\documentclass[12pt,a4paper]{amsart}
\usepackage{amscd,amsmath,amssymb,amsfonts,mathrsfs,times}
\input xy
\xyoption{all}
\numberwithin{equation}{subsection}
   \newtheorem{thm}{Theorem}[subsection]
   \newtheorem{lem}[thm]{Lemma}
   \newtheorem{prop}[thm]{Proposition}
   \newtheorem{cor}[thm]{Corollary}
   
   \newtheorem{defn}[thm]{Definition}

   \newtheorem{rem}[thm]{Remark}

\newenvironment{pf}{\medskip\noindent\emph{Proof.}}{\qed\par\medskip}
\newenvironment{pf*}[1]{\medskip\noindent\emph{#1.}}{\qed\par\medskip}
\begin{document}
\title[\'etale duality for constructible sheaves]{\'etale duality for constructible sheaves \endgraf on arithmetic schemes}
\author[U. Jannsen, S. Saito and K. Sato]{Uwe Jannsen,\, Shuji Saito \, and \, Kanetomo Sato}
\date{August, 2011}
\thanks{The authors are grateful to to the referee for offering numerous constructive comments to improve greatly the presentation of the paper. The third author carried out the research for this article during his stay at University of Southern California supported by JSPS Postdoctoral Fellowships for Research Abroad. He expresses his gratitude to Professors Wayne Raskind and Thomas Geisser for their great hospitality. Thanks are also due to Atsushi Shiho. The arguments for Theorems \ref{prop:appA} and \ref{thm:shiho} were inspired by discussions with him.}
\maketitle
%
%
%
%
%
%
%
%
%
%
\def\can{\text{\rm can}}
\def\cd{\text{\rm cd}}
\def\codim{\text{\rm codim}}
\def\ch{\text{\rm ch}}
\def\CH{\text{\rm CH}}
\def\cl{\text{\rm cl}}
\def\Coker{\text{\rm Coker}}
\def\cone{\text{\rm Cone}}
\def\cont{\text{\rm cont}}
\def\Cor{\text{\rm Cor}}
\def\DC{\mathcal E}
\def\dlog{d{\text{\rm log}}}
\def\et{\text{\rm \'{e}t}}
\def\ep{\epsilon}
\def\Ext{\text{\rm Ext}}
\def\Frac{\text{\rm Frac}}
\def\gys{\text{\rm Gys}}
\def\gysc{\gys^\circ}
\def\h{\text{\rm h}}
\def\H{H}
\def\Hom{\text{\rm Hom}}
\def\sHom{\hspace{0.5pt}{\mathscr H}\hspace{-1.2pt}om}
\def\id{\text{\rm id}}
\def\Image{\text{\rm Im}}
\def\ker{\text{\rm Ker}}
\def\L{\varLambda}
\def\lim{\mathop{\text{lim}}}
\def\loc{\text{\rm loc}}
\def\mod{\text{\rm mod}\ }
\def\op{\text{\rm op}}
\def\ord{\text{\rm ord}}
\def\qsimeq{\text{\rm qis.}}
\def\red{\text{\rm red}}
\def\res{\text{\rm res}}
\def\sh{\text{\rm sh}}
\def\Sh{\text{\rm Sh}}
\def\Spec{\text{\rm Spec}}
\def\Supp{\text{\rm Supp}}
\def\tor{\text{\rm tors}}
\def\tr{\text{\rm tr}}
\def\Tr{\text{\rm Tr}}
\def\val{\text{\rm val}}

\def\Sval{{S\text{-}\val}}

\def\II{I\hspace{-1.5pt}I}
\def\III{I\hspace{-1.5pt}I\hspace{-1.5pt}I}

\def\lwr{\wr\hspace{-1.5pt}}
\def\rwr{\hspace{-1.5pt}\wr}
\def\C{\mathbb C}
\def\N{\mathbb N}
\def\P{\mathbb P}
\def\Q{\mathbb Q}
\def\R{\mathbb R}
\def\Z{\mathbb Z}

\def\bA{\mathbb A}
\def\bC{\mathbb C}
\def\bF{\mathbb F}
\def\bG{\mathbb G}
\def\bH{\mathbb H}
\def\bK{\mathbb K}
\def\bL{\mathbb L}
\def\bN{\mathbb N}
\def\bP{\mathbb P}
\def\bQ{\mathbb Q}
\def\bR{\mathbb R}
\def\bT{\mathbb T}
\def\bZ{\mathbb Z}

\def\D{\mathscr D}
\def\K{\mathscr K}
\def\O{\mathscr O}

\def\cA{\mathscr A}
\def\cD{\mathscr D}
\def\cF{\mathscr F}
\def\cG{\mathscr G}
\def\cH{\mathscr H}
\def\cI{\mathscr I}
\def\cK{\mathscr K}
\def\cL{\mathscr L}
\def\cM{\mathscr M}
\def\cO{\mathscr O}
\def\cP{\mathscr P}
\def\cS{\mathscr S}
\def\cV{\mathscr V}

\def\fF{F}
\def\fO{\mathfrak o}
\def\fp{\mathfrak p}
\def\p{\mathfrak p}

\def\zn{\bZ/n}
\def\zpn{\bZ/p^n}
\def\zpr{\bZ/p^r}
\def\qz{\bQ/\bZ}
\def\qzp{\bQ_p/\bZ_p}
\def\mupn{\mu_{p^n}}

\def\dval{\partial^\val}
\def\dKyx{\dval_{y,x}}
\def\dKwy{\dval_{w,y}}
\def\dKwz{\dval_{w,z}}
\def\dKzx{\dval_{z,x}}

\def\dMwy{\partial_{w,y}}
\def\dMzx{\partial_{z,x}}

\def\ra{\to}
\def\lra{\longrightarrow}
\def\llllra{-\hspace{-5pt}-\hspace{-5pt}-\hspace{-6pt}\lra}
\def\Ra{\Rightarrow}
\def\Lra{\Longrightarrow}
\def\hra{\hookrightarrow}
\def\lmt{\longmapsto}
\def\isom{\hspace{9pt}{}^\sim\hspace{-16.5pt}\lra}
\def\lisom{\hspace{10pt}{}^\sim\hspace{-17.5pt}\longleftarrow}

\def\ssm{\smallsetminus}

\def\wt#1{\widetilde{#1}}
\def\wh#1{\widehat{#1}}

\def\ol#1{\overline{#1}}
\def\ul#1{\underline{#1}}

\def\bs#1{\boldsymbol{#1}}

\def\us#1#2{\underset{#1}{#2}}
\def\os#1#2{\overset{#1}{#2}}

\def\Gm{{\mathbb G}_{\hspace{-.6pt}\text{\rm m}}}
\def\Ga{{\mathbb G}_{\hspace{-.6pt}\text{\rm a}}}

\def\witt#1#2#3{W_{\hspace{-2pt}#2}{\hspace{1pt}}\Omega_{#1}^{#3}}
\def\logwitt#1#2#3{W_{\hspace{-2pt}#2}{\hspace{1pt}}\Omega_{{#1},{\log}}^{#3}}
\thispagestyle{empty}
\par
In this note we relate the following three topics for arithmetic schemes: a general duality for \'etale constructible torsion sheaves, a theory of \'etale homology, and the arithmetic complexes of Gersten-Bloch-Ogus type defined by K. Kato \cite{kk:hasse}.

In brief, there is an absolute duality using certain dualizing sheaves on these schemes, we describe and characterize the dualizing sheaves to some extent, relate them to symbol maps, define \'etale homology via the dualizing sheaves, and show that the niveau spectral sequence for the latter, constructed by the method of Bloch and Ogus \cite{Bloch-Ogus}, leads to the complexes defined by Kato. Some of these relations may have been expected by experts, and some have been used implicitly in the literature, although we do not know any explicit reference for statements or proofs. Moreover, the main results are used in a crucial way in a paper by two of us \cite{js}. So a major aim is to fill a gap in the literature, and a special emphasis is on precise formulations, including the determination of signs. But the general picture developed here may be of interest itself.
\smallskip
\subsection{Gersten-Bloch-Ogus-Kato complexes}\label{sect0-1}
For a scheme $X$ and a positive integer $n$ invertible on $X$, denote by $\zn(1) = \mu_n$ the \'etale sheaf on $X$ of $n$-th roots of unity, and let $\zn(r) = \mu_n^{\otimes r}$ be the $r$-fold Tate twist, defined for $r \in \Z$. As usual, we let
\[ \qzp(r) = \varinjlim{}_{n \ge 1} \ \zpn(r) \qquad \hbox{(for $p$ invertible on $X$).} \]
For a smooth variety $X$ over a perfect field of positive characteristic $p>0$ and integers $n>0$ and $r \ge 0$, $\logwitt X n r$ denotes the \'etale subsheaf of the logarithmic part of the $r$-th Hodge-Witt sheaf $\witt X n r$ (\cite{il} Chapter I 5.7), which are $\zpn$-sheaves. It is also noted $\nu^r_{n,X}$ in the literature. We denote
\[ \logwitt X \infty r = \varinjlim{}_{n \ge 1} \ \logwitt X n r \,, \]
where the transition maps $\logwitt X n r \to \logwitt X {n+1} r$ are given by factoring  the multiplication by $p$.
Let $X$ be a noetherian excellent scheme, and let $y$ and $x$ be points on $X$ such that $x$ has codimension $1$ in the closure $\ol {\{ y \}} \subset X$. Then for a prime number $p$, Kato (\cite{kk:hasse} \S1) defined `residue maps'
\begin{align}
& \H^{i+1}(y,\mupn^{\otimes r+1}) \lra \H^i(x,\mupn^{\otimes r}) & \hbox{(if $\ch(x)\not = p$)} \notag\\
& \H^i(y,\logwitt y n {r+1}) \lra \H^i(x,\logwitt x n r) & \hbox{(if $\ch(y)=\ch(x)=p$)} \label{eq-0.1.1} \\
& \H^{i+r+1}(y,\mupn^{\otimes r+1}) \lra \H^i(x,\logwitt x n r) \quad & \hbox{(if $\ch(y)=0$ and $\ch(x)=p$)}, \notag
\end{align}
where the maps of second and third type have non-zero target only for $i = 0,1$, and in case $i=1$ they are only defined if $[\kappa(x):\kappa(x)^p] \le p^r$. For a point $x \in X$, we wrote $H^*(x,-)$ for \'etale cohomology of $x = \Spec(\kappa(x))$, so this is just the Galois cohomology of $\kappa(x)$, the residue field at $x$. The sheaf $\logwitt x r n$ is the inverse image of $\logwitt U r n$, where $U$ is a dense smooth open subscheme of $\ol{\{x\}}$. These maps are defined via the Galois cohomology of discrete valuation fields, symbol maps on Milnor $K$-theory, and the valuation (see \S\ref{sect0-6} below). Therefore we will write $\dKyx$ for these maps, and denote sheafified variants in the same way. In particular, for $i=r=0$ the first and the last maps via Kummer theory correspond to the map
\[ \kappa(y)^\times/(\kappa(y)^\times)^{p^n} \lra \zpn \]
induced by the discrete valuations on the normalization of $\O_{\overline{\{y\}},x}$.

It has become customary to denote
\[ \zpn(r) := \logwitt X n r [-r]\]
for an (essentially) smooth scheme over a perfect field of characteristic $p$. With this notation, all maps above have the form
\[ \dKyx: H^{i+1}(y,\zpn(r+1)) \lra H^i(x,\zpn(r))\,. \]
Suppose now that $X$ is of finite type over a  noetherian regular excellent scheme of finite and pure dimension. Denote by $(X/S)_q$ the set of points $x \in X$ of `virtual dimension $q$ over $S$' (see \eqref{eq1-5-3'} below). When $S$ is the spectrum of a field then $(X/S)_q$ means $X_q$, the set of points on $X$ whose Zariski closure $\ol{\{x\}}$ has dimension $q$. In \cite{kk:hasse}, Kato showed that, for each triple of integers $i, j$ and $n>0$, the sequence
\begin{align*}
\dotsb \lra \bigoplus_{x\in (X/S)_r} H^{r+i}(x,\bZ/n(r+j)) & \lra \bigoplus_{x\in (X/S)_{r-1}}H^{r+i-1}(x,\bZ/n(r+j-1)) \lra  \\ \dotsb & \lra \bigoplus_{x\in (X/S)_0}H^i(x,\bZ/n(j))\,,
\end{align*}
whose maps have the components $\dKyx$\,, forms a complex $C_n^{i,j}(X)$. It was a major motivation for this paper
to understand the maps $\dKyx$ and these complexes in
terms of \'etale duality.

\smallskip
\subsection{\'Etale duality}\label{sect0-2}
A very general duality for constructible \'etale torsion sheaves has been established in \cite{sga4}. This is a {\it relative} duality, encoded in an adjunction
\begin{equation}\label{eq-0.2.1}
\Hom_X(\cF,Rf^!\cG) \cong \Hom_S(Rf_!\cF,\cG)
\end{equation}
for a separated morphism of finite type $f: X \to S$ and bounded complexes of \'etale torsion sheaves $\cF$ (on $X$) and $\cG$ (on $S$) with constructible cohomology sheaves (cf.\ \cite{sga4} XVIII 3.1.4.9). There is also a derived version, replacing $\Hom$ by $R\sHom$. To obtain an {\it absolute} duality for the cohomology groups of sheaves on $X$, in the spirit of Poincar\'e duality, one
needs an additional duality on the base scheme $S$. For arithmetic applications one is interested in schemes $X$ of finite type over $\Z$. Therefore we may assume that $S = \Spec(\fO_k)$, where $\fO_k$ is the ring of integers in a number field $k$. Here one has the Artin-Verdier duality
\begin{equation}\label{eq-0.2.2}
H^m_c(S,\cF) \times \Ext^{3-m}_S(\cF,\Gm) \lra H^3_c(S,\Gm) = \qz \,,
\end{equation}
where $H^m_c$ denotes the `cohomology with compact support' \cite{kk:hasse} which takes care of the archimedean places of $k$. But the figuring `dualizing sheaf' $\Gm$ is not torsion, so the relative duality above, for a scheme $X/S$, does not apply. Nevertheless, for such a higher-dimensional arithmetic scheme $X$, various absolute duality theorems have been obtained (cf.\ \cite{de2}, \cite{sp}, \cite{moser}, \cite{milne:adual}, \cite{Ge}), although always under some restrictions. For example $n$-torsion sheaves for $n$ invertible on $X$ have been considered, or $X$ was assumed to be smooth over $S$, or $X$ was assumed to be a scheme over a finite field.

Our approach is to introduce a complex of torsion sheaves $\qz(1)'_S$ on $S$ (see Definition \ref{def:DC-S}, \eqref{eq4-1-3}) so that one has a perfect duality as in \eqref{eq-0.2.2} when replacing $\Gm{}_{,S}$ by $\qz(1)'_S$. Next we define the dualizing `sheaf' (it is really a complex of sheaves) on $S$ as
\[ \D_S =  \qz(1)'_S[2]\,, \]
and on every separated $S$-scheme $X$ of finite type as
\[ \D_X = Rf^!\D_S \,, \]
where $f: X \to S$ is the structural morphism. Then, by using \eqref{eq-0.2.1}, \eqref{eq-0.2.2} and additional arguments, one gets a duality (cf.\ \S\ref{sect4})
\[ H^m_c(X,\cF) \times \Ext^{1-m}_X(\cF, \D_X) \lra H^1_c(X,\D_X) \lra \qz \,. \]
This is more or less formal, but we make the following three points. First, the duality is completely general: $X$ and the constructible complex $\cF$ can be arbitrary. Hence $X$ may be highly singular, and we may consider $p$-torsion sheaves even if $p$ is not invertible on $X$ (so in particular, if $X$ is an algebraic scheme over $\bF_p$), and the approach connects this `$p$-case' and the case `away from $p$' in a nice way. Secondly, we have a lot of information on the complex $\D_X$. Thirdly, it is this information that we need for the applications we have in mind, cf.\ \cite{js} 2.20, 2.21.

We describe the information on $\D_X$ separately for each $p$-primary part $\D_{X,p^\infty}$, where $p$ is a prime. Put \[ \Z/p^\infty := \qzp \quad\hbox{ and }\quad \mu_{p^\infty} := \bigcup_{n \ge 1} \  \mupn \,. \]
In the rest of this \S\ref{sect0-2}, suppose $n \in \bN \cup \{ \infty \}$. First we describe $\zpn(1)'_S$.
\par\medskip
(i)\; Let $S=\Spec(\fO_k)$ be as before. The complex $\zpn(1)'_S$ is, by definition, the mapping fiber of a morphism
\begin{equation}\label{eq-defS}
 \delta_S^\val = \delta^\val_{S,p^n} : Rj_*\mupn \lra i_* \zpn[-1] \,.
\end{equation}
Here $j: U= \Spec(\fO_k[p^{-1}]) \hra S$ is the open immersion, $i: Z = S \ssm U \hra S$ is the closed immersion of the complement, $\zpn$ is the constant sheaf on $Z$, and $\mupn$ is the sheaf of $p^n$-th roots of unity on $U$ (note that $p$ is invertible on $U$). One has $R^q j_*\mupn = 0$ for $q \ge 2$ (\cite{se} II.3.3\,(c)), and hence $\delta_S^\val$ is determined by the morphism $R^1j_*\mupn \to i_*\zpn$ it induces, and by adjunction and localization, this is in turn completely described by the induced morphisms
\[ \partial_x: k^\times/p^n = H^1(k,\zpn(1)) \lra H^0(x,\zpn) = \zpn \]
for each closed point $x\in Z=S \ssm U$. Then $\delta_S^\val$ is completely determined by defining $\partial_x$ to be the residue map \eqref{eq-0.1.1}, i.e., as $\ord_x\otimes \zpn$, where $\ord_x: k^\times \to \Z$ is the normalized discrete valuation corresponding to $x$.

Moreover, we will show that the mapping fiber of $\delta_S^\val$ is unique up to unique isomorphism in the derived category of sheaves on $S_\et$. In other words, $\zpn(1)'_S$ is the unique complex $\cF$ with $\cF|_U = \zpn(1)$, $Ri^!\cF = \zpn[-2]$, and for which the canonical morphism $Rj_* \cF|_U \to i_* Ri^!\cF [1]$ is the morphism $\delta_S^\val$ described above. See the remarks after Definition \ref{def:DC-S} for details.
%
\par\medskip
Now we list the properties of $\D_{X,p^n} = Rf^!\D_{S,p^n} = Rf^! \zpn(1)'_S[2]$ for $f: X \to S$ separated and of finite type.
\par\medskip
(ii)\; For $p$ invertible on $X$, $\D_{X,p^\infty}$ is the usual dualizing sheaf for the `prime-to-$p$ theory' over $\fO_K[p^{-1}]$. In particular, $\D_{X,p^\infty} = \qzp(d)[2d]$ if $X$ is regular of pure dimension $d$. Here we use the absolute purity due to Gabber \cite{fujiwara}.
\par\medskip
(iii)\; For $X$ of characteristic $p$, i.e., of finite type over the prime field $\bF_p$, and of dimension $d$, $\D_{X,p^\infty}$ is represented by the explicit complex
\begin{equation}\label{eq-defM}
 \cM_X \; : \; \bigoplus_{x\in X_d} \logwitt x \infty d \lra \bigoplus_{x\in X_{d-1}} \ \logwitt x \infty {d-1} \lra \dotsb \lra \bigoplus_{x\in X_0} \ \qzp
\end{equation}
introduced by Moser \cite{moser} p.\ 128 (except that we put the rightmost term in degree zero, while Moser rather considers the complex $\wt{\nu}{}^d_{\infty,X} := \cM_X[-d]$). In fact, we generalize Moser's duality over finite fields \[ H^m_c(X,\cF) \times \Ext^{d+1-m}_X(\cF,\cM_X[-d]) \to H^{d+1}_c(X,\cM_X[-d]) = \qzp \] in the following way: We extend the duality to arbitrary perfect ground fields $k$ of characteristic $p$, and show that $\cM_X$ is in fact $Rg^!\qzp$, where $g: X \to \Spec(\bF_p)$ is the structural morphism. Together with the well-known duality of finite fields, this immediately gives back Moser's theorem. By Gros and Suwa \cite{gs} 1.6, one has $\cM_X = \logwitt X \infty d [d]$, if $X$ is regular.
\par\medskip
(iv)\; Finally, for $X$ flat over $S = \Spec(\fO_K)$, consider the closed immersion
\[\xymatrix{ i : Y := X\otimes_{\Z} \bF_p \, \ar@{^{(}->}[r] & X }\]
and the open immersion
\[\xymatrix{ j: U := X[p^{-1}] \, \ar@{^{(}->}[r] & X }\]
of the complement. There is a morphism
\[ \delta_X^\Sval = \delta_{X,p^\infty}^\Sval : Rj_*\cD_{U,p^\infty} \lra i_* \cD_{Y,p^\infty}[1] \]
obtained from $\delta_S^\val$ (cf.\ (i)) via $Rf^!$, where $f : X \to S$ denotes the structural map. The source and target are studied in (ii) and (iii) above, respectively, and it is clear from the definitions that $\D_{X,p^\infty}$ is a mapping fiber of $\delta_X^\Sval$. In general, such a mapping fiber is not unique (for the lack of the unicity of isomorphisms), but one of our main results is the following:
$\cD_{X,p^\infty}$ is a unique mapping fiber of $\delta_X^\Sval$ up to unique isomorphism (cf.\ Theorem \ref{thm.cone} and Remark \ref{rem.cone})
 and moreover when $U$ is smooth, $\delta_X^\Sval$ is uniquely characterized by the property that, for every generic point $y \in Y$ and every generic point $\xi \in U$ which specializes to $y$, the induced map
\[ H^d(\xi,\qzp(d)) \lra H^{d-1}(y,\qzp(d-1)) = H^0(y,\logwitt y {\infty} {d-1}) \]
coincides with the residue map in \eqref{eq-0.1.1}, cf.\ Theorem \ref{Th.1-1}\,(3). When $X$ is proper (but $U$ arbitrary), we have a similar uniqueness property.
\par
There is another morphism
\[ \delta_{U,Y}^\loc(\cD_{X,p^\infty}) : Rj_*\cD_{U,p^\infty} \lra i_* \cD_{Y,p^\infty}[1]\,, \]
the connecting morphism of localization theory for $\cD_{X,p^\infty}$. We will also prove that this morphism agrees with $\delta_X^\Sval$ up to a sign, cf.\ \eqref{eq3-9-5}.

\smallskip
\subsection{\'Etale homology}\label{sect0-3}
Let $k$ be a perfect field, and let $X$ be a separated scheme of finite type over $k$. For integers $n>0$, $a$ and $b$, we define the \'etale homology of $X$ by
\[ H_a(X,\zn(b)) = H^{-a}(X,Rf^!\zn(-b))\,, \]
where $f: X \to \Spec(k)$ is the structural morphism. Note that for $\ch(k) = p > 0$, we have $\zpr(-b) = \logwitt k r {-b}[b]$, which is the constant sheaf $\zpr$ for $b=0$, zero for $b<0$ (because $k$ is perfect), and zero by definition for $b>0$. Therefore we will either assume that $n$ is invertible in $k$, or that $b=0$. These groups satisfy all properties of a (Borel-Moore type) homology theory, cf.\ \cite{Bloch-Ogus} 1.2, \cite{js} 2.1\,(a). Thus the method of Bloch and Ogus provides a converging niveau spectral sequence (\cite{Bloch-Ogus} 3.7)
\begin{equation}\label{eq-0.3.1}
E^1_{s,t}(X,\zn(b)) = \bigoplus_{x \in X_s} \ H_{s+t}(x,\zn(b)) \Lra  H_{s+t}(X,\zn(b)).
\end{equation}
Here we put
\[ H_a(x,\zn(b)) = \varinjlim_{V \subset \ol{\{x\}}} \ H_a(V,\zn(b)) \]
and the limit is taken over all non-empty open subvarieties $V\subset \overline{\{x\}}$. If $V$ is smooth of pure dimension $d$ over $k$, then one has a canonical purity isomorphism
\[ H_a(V,\zn(b)) \cong H^{2d-a}(V,\zn(d-b)) \]
between homology and cohomology. This is one of the main results of the Artin-Verdier duality \cite{sga4} in the case $n$ is invertible in $k$, and follows from our results in \S\ref{sect2} for the other case. As a consequence, one has canonical isomorphisms
\[ \bigoplus_{x\in X_s} \ H_{s+t}(x,\zn(b)) \cong \bigoplus_{x\in X_s} \ H^{s-t}(x,\zn(s-b))\,, \]
and the complex $E^1_{*,t}$ of $E^1$-terms of the spectral sequence can be identified with a complex
\begin{align}
\notag
\dotsb \lra \bigoplus_{x\in X_s} \ H^{s-t}(x,\bZ/n(s-b)) & \lra \bigoplus_{x\in X_{s-1}} \ H^{s-t-1}(x,\bZ/n(s-b-1)) \lra \\
\label{eq-0.3.2}
\dotsb & \lra \bigoplus_{x\in X_0} \ H^{-t}(x,\bZ/n(-b))\,,
\end{align}
where we place the last term in degree zero. Another main result of this paper is that this complex coincides with the Kato complex $C_n^{-t,-b}(X)$ mentioned in \S\ref{sect0-1}, up to well-defined signs. In \S\ref{sect1} we also give an absolute variant of this result, for the case that $X$ is a regular excellent noetherian scheme and $n$ is invertible on $X$.

Finally let $X$ be a separated scheme of finite type over $S = \Spec(\fO_K)$, where $K$ is a number field, and let $n$, $a$ and $b$ be integers. If $n$ is invertible on $X$, we define the \'etale homology as
\[ H_a(X,\zn(b)) = H^{-a}(X,Rf^!\zn(-b)), \]
where $f: X \to S[n^{-1}]$ is the structural morphism. If $n$ is not invertible on $X$, we just consider the case $b=-1$ and define
\[ H_a(X,\zn(-1)) = H^{-a}(X,Rf^!\zn(1)'_S)\,, \]
where $f: X \to S$ is the structural morphism, and $\zn(1)'_S$ has the $p$-primary components $\Z/p^{r_p}(1)'_S$ from (i) for $n = \prod \, p^{r_p}$. Again, in both cases this defines a homology theory in the sense of \cite{js} 2.1\,(a) (cf.\ \cite{Bloch-Ogus} 1.2), and one gets a niveau spectral sequence with exactly the same numbering as in \eqref{eq-0.3.1}. By the purity isomorphisms explained above, the complex of $E^1$-terms is identified with a complex
\begin{align}
\notag
\dotsb \to \bigoplus_{x\in (X/S)_s} \ H^{s-t-2}(x,\zn(s-b-1)) & \to \bigoplus_{x\in (X/S)_{s-1}} \ H^{s-t-3}(x,\zn(s-b-2))  \to  \\
\label{eq-0.3.3}
\dotsb & \to \bigoplus_{x\in (X/S)_0} \ H^{-t-2}(x,\zn(-b-1)),
\end{align}
cf.\ \cite{js}. The difference in numbering between \eqref{eq-0.3.2} and \eqref{eq-0.3.3} is explained by the purity results for the inclusion of the fibers $X_P \hra X$ over closed points $P\in S$. A third main result of this paper is that, also in this mixed characteristic case, this complex coincides with a Kato complex, viz., $C_n^{-t-2,-b-1}(X)$. In fact, this gives an alternative definition of the Kato complexes under consideration, which is very useful for working with them.

\smallskip
\subsection{Notations and conventions}\label{sect0-4}
For an abelian group $M$ and a positive integer $n$, $M/n$ (resp.\ ${}_n M$) denotes the cokernel (resp.\ the kernel) of the map $M \os{\times n}{\ra} M$. \par
In this paper, unless indicated otherwise, all cohomology groups of schemes are taken for the \'etale topology. \par
For a scheme $X$, we will use the following notation. For a point $x \in X$, $\kappa(x)$ denotes its residue field, and $\ol x$ denotes $\Spec(\ol {\kappa(x)})$, the spectrum of a separable closure of $\kappa(x)$. For a point $x \in X$ and an \'etale sheaf $\cF$ on $X$, we define
\[ \H^*_x(X,\cF) := \H^*_x(\Spec(\O_{X,x}),\cF). \]
For a non-negative integer $q$, we use the notation $X^q$ and $X_q$ (only) in the following cases.
If $X$ is pure-dimensional, then $X^q$ denotes the set of points on $X$ of codimension $q$.
If $X$ is a scheme of finite type over a field, then $X_q$ denotes the set of points on $X$ whose closure in $X$ has dimension $q$.
\par
\smallskip
\subsection{Connecting morphism of localization sequences}\label{sect0-4'}
Let $X$ be a scheme and let $n$ be a non-negative integer. Let $i:Z \hra X$ be a closed immersion, and let $j:U \hra X$ be the open complement $X \ssm Z$. For an object $\cK \in D^+(X_\et,\zn)$, we define the morphism
\[ \delta_{U,Z}^{\loc}(\cK) : Rj_*j^*\cK \lra i_*Ri^!\cK[1] \quad \hbox{ in } \; D^+(X_\et,\zn) \]
 as the connecting morphism associated with the semi-splitting short exact sequence of complexes
\[ 0 \lra i_*i^!I^\bullet \lra I^\bullet \lra j_*j^*I^\bullet \lra 0 \]
(\cite{sga4.5} Cat\'egories D\'eriv\'ees I.1.2.4), where $I^\bullet$ is a resolution of $\cK$ by injective $\bZ/n$-sheaves on $X_\et$. It induces the usual connecting morphisms
\[ \delta_{U,Z}^\loc(\cK) : R^qj_* j^*\cK \lra i_* R^{q+1}i^!\K\,,\]
or the connecting morphisms in the localization sequence for $(X,Z,U)$:
\[ \delta_{U,Z}^\loc(\cK) : H^q(U,j^* \cK) \lra H^{q+1}_Z(X,\cK)\,. \]
The morphism $\delta_{X,U}^\loc(\cK)$ is functorial in $\cK$, but does not commute with shift functors in general. In fact, we have
\begin{equation}\label{eq-0.4.1}
\delta_{U,Z}^{\loc}(\cK)[q] = (-1)^q \cdot \delta_{U,Z}^{\loc}(\cK[q]) \quad \hbox{for } \; q \in \bZ.
\end{equation}
By the convention in \cite{sga4} XVII.1.1.1 (which we follow and is usually taken, but which is opposite to the one in \cite{sga4.5} Cat\'egories D\'eriv\'ees I.1.2.1), the following triangle is distinguished in $D^+(X_\et,\zn)$:
\begin{equation}\label{eq-0.4.2}
\xymatrix{ i_*Ri^!\cK \ar[r]^-{i_*} & \cK \ar[r]^-{j^*} & Rj_*j^*\cK \ar[rr]^-{-\delta_{U,Z}^{\loc}(\cK)} && i_*Ri^!\cK[1]\,, }
\end{equation}
where the arrow $i_*$ (resp.\ $j^*$) denotes the adjunction morphism $i_*Ri^! \ra \id$ (resp.\ $\id \ra Rj_*j^*$). We generalize the above definition of connecting morphisms to the following situation. Let $x$ be a point on $X$ and let $i_x$ be the natural map $x \hra X$. We define a functor
\[ Ri_x^! : D^+(X_\et,\zn) \lra D^+(x_\et,\zn) \]
as $\iota_x^*Ri^!$, where $i$ denote the natural closed immersion $\ol {\{x\}} \hra X$ and $\iota_x$ denotes the natural map $x \hra \ol {\{x\}}$. Note that $Ri_x^!$ is not right adjoint to $i_{x*}$ unless $x$ is a closed point on $X$. Now let $y$ and $x$ be points on $X$ such that $x$ has codimension $1$ in the closure $T:=\ol {\{ y \}} \subset X$. Put $Y:=\Spec(\O_{T,x})$, and let $i_x$ (resp.\ $i_y$, $i_Y$, $i_T$, $\iota_Y$) be the natural map $x \ra X$ (resp.\ $y \ra X$, $Y \ra X$, $T \hra X$, $Y \hra T$). Then we define a connecting morphism
\[ \delta_{y,x}^{\loc}(\cK) : Ri_{y*}Ri_y^!\cK \lra Ri_{x*}Ri_x^!\cK[1] \quad \hbox{ in } \; D^+(X_\et,\zn) \]
as $Ri_{Y*}(\delta_{y,x}^{\loc}(Ri_Y^!\cK))$. Here we defined $Ri_Y^!$ as $\iota_Y^*Ri_T^!$.
\smallskip
\subsection{Derived categories}\label{sect0-5}
We shall often use the following facts. Let $\cA$ be an abelian category, and let $D^*(\cA)$ be its derived category with boundary condition $*\in \{\emptyset, +, -, b\}$.

\subsubsection{}\label{sect0-5-1}
A sequence $A \os{\alpha}\to B \os{\beta}\to C \os{\gamma}\to A[1]$ in $D^*(\cA)$ is a distinguished triangle if and only if $B \os{\beta}\to C \os{\gamma}\to A[1] \os{-\alpha[1]}\to B[1]$ is a distinguished triangle. (This is the axiom (TR2) for triangulated categories, \cite{sga4.5} Cat\'egories Deriv\'ees I.1.1.)

\subsubsection{}\label{sect0-5-2}
Given a diagram
\[\xymatrix{
A \ar[r]^a \ar[d]_f & B \ar[r]^{b} & C \ar[r] \ar[d]^h & A[1] \ar[d]^{f[1]} \\ A' \ar[r]^{a'} & B' \ar[r]^{b'} & C' \ar[r] & A'[1]
}\]
in which the rows are distinguished triangles and the last square commutes, there is a morphism $g: B \to B'$ making the remaining squares commutative, i.e., giving a morphism of distinguished triangles. Moreover one has
\addtocounter{thm}{2}
\begin{lem}\label{lem.derived}
The morphism $g$ is unique in the following three cases{\rm:}
\begin{enumerate}
\item[{\rm (1)}] $\Hom_{D(\cA)}(B,A') = 0$.
\item[{\rm (2)}] $\Hom_{D(\cA)}(C,B') = 0$.
\item[{\rm (3)}] $\Hom_{D(\cA)}(C,A') = 0$ \, and \, $\Hom^{-1}_{D(\cA)}(A,C') = 0$.
\end{enumerate} \end{lem}
\begin{pf}
There is an induced commutative diagram with exact rows and columns
\[\xymatrix{
& \Hom(C,A') \ar[r] \ar[d] & \Hom(B,A') \ar[d] \\ & \Hom(C,B') \ar[r]^{b^*} \ar[d]_{b'_*} & \Hom(B,B') \ar[r]^{a^*} \ar[d]_{b'_*} & \Hom(A,B') \\ \Hom^{-1}(A,C') \ar[r] & \Hom(C,C') \ar[r]^{b^*} & \Hom(B,C')\,.\hspace{-3pt}
}\]
Suppose $g_1$ and $g_2$ both make the previous diagram commutative. Then the element $g_1 - g_2 \in \Hom(B,B')$ is mapped
to zero in $\Hom(A,B')$ and $\Hom(B,C')$. Under conditions (1) and (2), either the right hand $b'_*$ or $a^*$ is injective, so the claim follows. Under condition (3), the left hand $b'_*$ and the lower $b^*$ are both injective, and again we get $g_1 - g_2 = 0$.
\end{pf}
\par
\stepcounter{subsubsection}
\subsubsection{}\label{sect0-5-4}
Let $q,r$ be integers, and let $M$ be an object in $D(\cA)$ which is concentrated in degrees $\leq r$. Let $N$ be an object in $D(\cA)$ which is concentrated in degrees $\ge 0$. Then we have
\[ \Hom_{D(\cA)}(M,N[-q]) =
 \begin{cases}
 \Hom_{\cA}(\cH^q(M),\cH^0(N)) & \quad \hbox{(if $q=r$)} \ \ \dotsb\dotsb (1)\\
 0 & \quad \hbox{(if $q>r$)} \ \ \dotsb\dotsb (2)
 \end{cases}\]
Here for $s \in \bZ$, $\cH^s(M)$ denotes the $s$-th cohomology object of $M$. These facts are well-known and easily proved, using \cite{pervers} 1.3.2 and \cite{sga4.5} Cat\'egories D\'eriv\'ees I.1.2.
\smallskip
\subsection{Kato's residue maps}\label{sect0-6}
We recall Kato's definition of the residue maps in \eqref{eq-0.1.1}. Consider a noetherian excellent scheme $X$ and points $x,y\in X$ such that $x$ lies in $Z=\overline{\{y\}}$ and has codimension 1 in $Z$. The construction only depends on $Z$ (with the reduced subscheme structure). Put $A:=\cO_{Z,x}$, a local domain of dimension $1$. We may further replace $Z$ with $\Spec(A)$.
\par
\medskip
(I)\; {\it Regular case}.\;
First consider the case that $A$ is regular, i.e., a discrete valuation ring. Then $K:= \kappa(y) = \Frac(A)$ is a discrete valuation field and $k:=\kappa(x)$ is the residue field of $A$, i.e., of the valuation. The residue map
\[ \dval = \dKyx : H^{i+1}(K,\zpn(r+1)) \lra H^i(k,\zpn(r)) \]
is obtained by restricting to the henselization $K^\h$ (which corresponds to restricting to the henselization $A^h = \cO_{Z,x}^\h$) and defining a map for the discrete valuation field $K^\h$ which has the same residue field $k$. Hence we may and will assume that $K$ is henselian (i.e., $A=A^\h$ and $K=K^\h$). Let $K^\sh$ be the maximal unramified extension of $K$ (corresponding to the strict henselization $A^\sh = \cO_{Z,\overline x}$).
\par
\medskip
(I.1)\;
If $p\neq \ch(k)$, we first have a map
\[ H^1(K,\zpn(1)) \lisom K^\times/(K^\times)^{p^n} \lra \zpn = H^0(k,\zpn)\,, \]
where the first arrow is the Kummer isomorphism, and the second is induced by the valuation. This is $\dval$ for $(i,r) =(0,0)$. In general $\dval$ is the composition
\[ H^{i+1}(K,\zpn(r+1)) \lra H^i(k,H^1(K^\sh,\zpn(r+1))) \lra H^i(k,\zpn(r))\,. \]
Here the first map is an edge morphism from the Hochschild-Serre sequence for $K^\sh/K$ (note that $\cd(K^\sh)=1$), and the second map is induced by (the Tate twist of) the previously defined map.
\par
\medskip
(I.2)\;
Now let $p=\ch(k)$ (and recall that $K$ is henselian). In this case $H^i(k,\zpn(r))$ $= 0$ for $i \ne r, r+1$. Assume that $i = r$. Then $\dval$ is defined by the commutativity of the diagram
\begin{equation}\label{eq-0.6.2}
\xymatrix{
H^{r+1}(K,\zpn(r+1)) \ar[r]^-{\dval} & H^r(k,\zpn(r)) \\
K^M_{r+1}(K)/p^n \ar[u]^{h^{r+1}}_{\rwr} \ar[r]^-{\partial} & K^M_r(k)/p^n \,. \ar[u]^{h^r}_{\rwr}
}\end{equation}
Here $K^M_r(F)$ is the $r$-th Milnor $K$-group of a field $F$, $h^r$ is the symbol map into Galois cohomology, and $\partial$ is the suitably normalized residue map for Milnor $K$-theory. By definition,
\[ h^r(\{a_1,\dotsc,a_r\}) = h^1(a_1) \cup \dotsb \cup h^1(a_r) \in H^r(F,\zpn(r))\,, \]
where $h^1: F^\times/p^n \to H^1(F,\zpn(1))$ is defined as follows: it is the Kummer isomorphism if $p$ is invertible in $F$, and it is the isomorphism $\dlog: F^\times/p^n \to H^0(F,\logwitt F n 1)$ if $\ch(F)=p$. It is known that, under our assumptions, the symbol maps $h^i$ in \eqref{eq-0.6.2} are isomorphisms (\cite{Bloch-Kato} \S2, \S5). Finally, if $\pi$ is a prime element for $K$, then $\partial$ is determined by the property that
\[ \partial(\{\pi,a_1,\dotsc,a_r\}) = \{\ol{a_1},\dotsc,\ol{a_r}\}\,, \]
for units $a_1,\ldots,a_r \in A^\times$, where $\overline{a_i}$ denotes the residue class of $a_i$ in the residue field $k$.
\par
\medskip
(I.3)\; Now let $i=r+1$. In this case we assume $[k:k^p] \le p^r$.
Then the residue map $\dval$ is defined as the composition ($\ol k$ denotes the separable closure of $k$)
\[ \xymatrix{
H^1(k,H^{r+1}(K^{\sh},\zpn(r+1)) \ar[r]^-{(**)} & H^1(k,H^r(\ol k,\zpn(r))) \\
H^{r+2}(K,\zpn(r+1)) \ar[u]^{(*)}_{\rwr} \ar@{.>}[r]^-{\dval} & H^{r+1}(k,\zpn(r)) \,. \ar[u]^{(*)}_{\rwr}
}\]
Here the isomorphisms $(*)$ come from the Hochschild-Serre spectral sequences and the fact that $\cd_p(k)\le 1$ and
\[ H^{j+1}(K^{\sh}, \zpn(r+1)) = 0 = H^j(\ol k,\zpn(r)) \quad \hbox{ for \; $j>r$.} \]
The map $(**)$ is induced by the map
\[ H^{r+1}(K^{\sh},\zpn(r+1)) \lra H^r(\ol k,\zpn(r)) \]
defined in (I.2). In \cite{kk:hasse} the completion $\wh K$ is used instead of the henselization $K^\h$, but this gives the same, because the map
\[ H^j(K^{\sh}, \zpn(r+1)) \lra H^j(\wh K^{\sh}, \zpn(r+1))\]
is an isomorphism (\cite{katokuzumaki} proof of Theorem 1). Indeed $\wh K^{\sh}/K^{\sh}$ is separable by excellency of $X$.
\par \medskip
(II)\; {\it General case}.\;
Now consider the case that $A$ is not necessarily regular. In this case let $Z' \to Z = \Spec(A)$ be the normalization. Note that $Z'$ is finite over $Z$ because the latter is excellent. Then we define
\[\dKyx(a) \, = \, \sum_{x'|x} \ \Cor_{\kappa(x')/\kappa(x)}(\dval_{y,x'}(a)) \qquad (a \in H^{i+1}(y,\zpn(r+1)))\]
where the sum is taken over all points $x' \in Z'$ lying over $x$,
\[ \dval_{y,x'}: H^{i+1}(y,\zpn(r+1)) \lra H^i(x',\zpn(r)) \]
is the residue map defined for the discrete valuation ring $\cO_{Z',x'}$, and
\begin{equation}\label{eq-0.6.5}
\Cor_{\kappa(x')/\kappa(x)}: H^i(x',\zpn(r)) \lra H^i(x,\zpn(r))
\end{equation}
is the corestriction map in Galois cohomology. For $p$ invertible in $\kappa(x)$ this last map is well-known. For $\kappa(x)$ of characteristic $p$ and $i=r$, this corestriction map is defined as the composition
\begin{align}
\label{eq-0.6.6}
\H^0(x',\logwitt {x'} n r) & \os{(h^r)^{-1}}{\isom} K^M_r(\kappa(x'))/p^n \os{N_{x'/x}}\lra K^M_r(\kappa(x))/p^n \\
\notag
& \;\, \os{h^r}\lra \; \H^0(x,\logwitt x n r) \,.
\end{align}
This implies that the diagram \eqref{eq-0.6.2} is also commutative in this case. For the remaining case $i=r+1$ we may proceed as follows. It is easy to see that the map \eqref{eq-0.6.6} is compatible with \'etale base-change in $\kappa(x)$. Therefore we get an induced corestriction or trace map
\begin{equation}\label{eq-0.6.7}
\tr_{x'/x}: \pi_* \logwitt {x'} n r \lra \logwitt x n r \,.
\end{equation}
Then we define the corestriction \eqref{eq-0.6.5} for $\ch(\kappa(x))=p$ and $i=r+1$ as
\begin{equation}\label{eq-0.6.8}
\tr_{x'/x} : H^1(x',\logwitt {x'} n r) \lra H^1(x,\logwitt x n r) \,,
\end{equation}
the map induced by \eqref{eq-0.6.7}. If $\kappa(x)$ is finitely generated over a perfect field $k$, the morphisms \eqref{eq-0.6.6}, \eqref{eq-0.6.7} and \eqref{eq-0.6.8} agree with the trace map in logarithmic Hodge-Witt cohomology defined by Gros \cite{gros:purity}. See the appendix for this and further compatibilities.
\smallskip
\subsection{Sheafified residue maps}\label{sect0-7}
We further define residue maps for sheaves by sheafifing the above residue maps of Galois cohomology groups. The setting remains as in \S\ref{sect0-6}. Consider a noetherian excellent scheme $X$ and points $x,y\in X$ such that $x$ lies in $Z=\overline{\{y\}}$ and has codimension 1 in $Z$. Let $i_x : x \hra X$ and $i_y : y \hra X$ be the natural maps. We would like to define homomorphisms of sheaves on $X_\et$
\begin{align}
& R^{j+1}i_{y*}\mu_{p^n}^{\otimes r+1} \lra R^ji_{x*}\mu_{p^n}^{\otimes r} & \hbox{(if $\ch(x)\not = p$)} \notag\\
& i_{y*}(\logwitt y n {r+1}) \lra i_{x*}(\logwitt x n r) & \hbox{(if $\ch(y)=\ch(x)=p$)} \label{eq-0.1.1'} \\
& R^{r+1}i_{y*}\mu_{p^n}^{\otimes r+1} \lra i_{x*}(\logwitt x n r) \quad & \hbox{(if $\ch(y)=0$ and $\ch(x)=p$)}. \notag
\end{align}
To define the first map, it is enough to construct a morphism
\begin{align}
& Ri_{y*}\mu_{p^n}^{\otimes r+1} \lra Ri_{x*}\mu_{p^n}^{\otimes r}[-1] \quad \hbox{ in } \; D^+(X_\et,\bZ/p^n).  \notag
\end{align}
By adjunction, to define the maps in \eqref{eq-0.1.1'} it is enough to construct
\begin{align}
& i_x^*Ri_{y*}\mu_{p^n}^{\otimes r+1} \lra \mu_{p^n}^{\otimes r}[-1]  \notag \\
& i_x^*i_{y*}(\logwitt y n {r+1}) \lra \logwitt x n r \label{eq-0.1.1''} \\
& i_x^*\,R^{r+1}i_{y*}\mu_{p^n}^{\otimes r+1} \lra \logwitt x n r, \notag
\end{align}
respectively, on $x_\et$. Let $A$ be the strict henselization of $\cO_{Z,x}$ at $\ol x$. Let $Z_1,\dotsc,Z_a$ be the distinct irreducible components of $\Spec(A)$.
Let $A_i$ be the affine ring of $Z_i$, which is a strict henselian local domain of dimension $1$ with residue field $\kappa(\ol x)$. Let $\eta_i$ be the generic point of $Z_i$.
Noting the fact $\cd_p(\eta_i) = 1$ in the first case (cf.\ \cite{sga5} I.5) and looking at stalks, it is enough to construct
\begin{align}
& \bigoplus_{i=1}^a \ H^1(\eta_i,\mu_{p^n}^{\otimes r+1}) \lra H^0(\ol x,\mu_{p^n}^{\otimes r}) \notag \\
& \bigoplus_{i=1}^a \ H^0(\eta_i,\logwitt y n {r+1}) \lra H^0(\ol x, \logwitt x n r) \label{eq-0.1.1'''} \\
& \bigoplus_{i=1}^a \ H^{r+1}(\eta_i,\mu_{p^n}^{\otimes r+1}) \lra H^0(\ol x,\logwitt x n r), \notag
\end{align}
where we have used the fact that $\cd_p(\eta_i)=1$ if $p \ne \ch(k)$. We define these maps as the sum of the maps in \eqref{eq-0.1.1} for the $Z_i$'s, which provide us with the maps in \eqref{eq-0.1.1'}.

%
%
\newpage
\section{Duality for arithmetic schemes}\label{sect4}
\medskip
The aim of this section is to prove a general duality for constructible sheaves on separated schemes of finite type over $\Z$. The main result of this section will be stated in Theorem \ref{thm:global-duality} below.
\par\smallskip
\subsection{Dualizing complex and higher-dimensional duality}
Let $k$ be a number field with ring of integers $\fO_k$, and let $S=\Spec(\fO_k)$. Let $p$ be a prime number, let $n \in \bN \cup \{\infty\}$ and put $U:= \Spec(\fO_k[p^{-1}])$. Let $j: U= \Spec(\fO_k[p^{-1}]) \hra S$ be the open immersion, and let $i: Z = S \ssm U \hra S$ be the closed immersion of the complement. Let $\zpn$ be the constant sheaf on $Z$, and let $\mupn$ be the sheaf of $p^n$-th roots of unity on $U$. We consider the following diagram:
\[\xymatrix{ U \; \ar@{^{(}->}[r]^-j & S & \ar@{_{(}->}[l]_-i \; Z \,. }\]
\begin{defn}\label{def:DC-S}
For each integer $n \ge 1$ define
\[ \zpn(1)'_S := \cone(\delta_S^\val: Rj_*\mupn \to i_* \zpn [-1])[-1] \in D^b(S_\et,\zpn), \]
the mapping fiber of the morphism $\delta_S^\val$ defined in \eqref{eq-defS}.
\end{defn}
\noindent
Here $D^b(S_\et,\bZ/p^\infty)$ denotes the derived category of bounded complexes of \'etale sheaves on $S$ whose sections are torsion of $p$-power
order.

\smallskip
In general, mapping cone or fiber of a morphism in a derived category is only well-defined up non-canonical isomorphism. However in our case it is well-defined up to a unique isomorphism, because we can apply the criterion of Lemma \ref{lem.derived}\,(1). Indeed the complex $Rj_*\mupn$ is concentrated in $[0,1]$, $A[1]=i_* \zpn [-1]$ is concentrated in degree 1, and $\delta_S$ induces a surjection $R^1j_*\mupn \twoheadrightarrow i_* \zpn$ so that the mapping fiber $B$ is concentrated in $[0,1]$ as well. Therefore $\Hom_{D(S,\zpn)}(B,A)=0$. (This argument should replace the reasoning in \cite{js} p.\ 497, where the criterion is misstated.)

By the above, there is a canonical exact triangle
\stepcounter{equation}
\begin{equation}\label{eq3-9-2}
i_*\zpn[-2] \os{g}\lra \zpn(1)'_S \os{t}\lra Rj_*\mupn \os{\delta_S^\val}\lra i_*\zpn[-1]\,,
\end{equation}
which induces canonical isomorphisms
\begin{equation}\label{isos.DC-S}
t: j^*(\zpn(1)'_S) \cong \mupn \quad \hbox{ and } \quad g: \zpn[-2] \cong Ri^!(\zpn(1)'_S)\,.
\end{equation}
We write
\begin{equation}\label{eq4-1-3'}
\qzp(1)'_S:= \bZ/p^\infty(1)' \in D^b(S_\et,\bZ/p^\infty),
\end{equation}
and often also regard this complex as an object of $D^b(S_\et)$. We further define
\begin{equation}\label{eq4-1-3}
\qz(1)'_S:= \bigoplus_{p} \ \qzp(1)' \in D^b(S_\et),
\end{equation}
where $p$ runs through all rational prime numbers. We will explain a version of Artin-Verdier duality using $\qz(1)'_S$ below.
\par\medskip
Now let $X$ be a separated scheme of finite type over $S$, with structural morphism $f: X \to S$. We define
\[ \cD_{X,p}:= Rf^!\qzp(1)'_S[2] \quad \hbox{ and } \quad \cD_{X}:= \bigoplus_p \ \cD_{X,p} \quad \hbox{(cf.\ \S\ref{sect0-2})}. \]
One of the main purposes of this paper is to analyze the object $\cD_{X,p}$ for each prime number $p$. We will describe $\cD_{X,p}$ over $X[p^{-1}]$ and $X \otimes_\bZ \bF_p$. The restriction of $\cD_{X,p}$ to $X[p^{-1}]$ is the usual \'etale dualizing complex. $\cD_{X,p}$ is obtained by glueing the \'etale dualizing complex on $X[p^{-1}]$ and Moser's complex \eqref{eq-defM} of $X \otimes_\bZ \bF_p$. We prove that there is a unique way to glue them when we impose some compatibilities concerning residue maps associated to specializations from a point in $X[p^{-1}]$ to a point in $X \otimes_\bZ \bF_p$ (see \S\ref{sect3-11} below). \par
\bigskip
For $\cL \in D^+(X_\et)$, we define the $m$-th \'etale cohomology group {\it with compact support} as
\[ \H^m_c(X,\cL) := \H^m_c(S,Rf_!\cL)\,, \]
where for an \'etale sheaf or a complex of \'etale sheaves $\cF$ on $S$, $\H^m_c(S,\cF)$ denotes the $m$-th \'etale cohomology group with compact support (see e.g., \cite{milne:adual} II.2, \cite{kk:hasse} \S3 for generalities). The main result of this section is the following duality (see also \cite{de2}, \cite{sp}):
\begin{thm}\label{thm:global-duality}
\begin{enumerate}
\renewcommand{\labelenumi}{(\arabic{enumi})}
\item There is a canonical trace map
\[ \tr_X:\H_c^1(X,\cD_X) \lra \qz\,. \]
\item For $\cL \in D^b(X_\et)$ with constructible torsion cohomology sheaves, the pairing
\[ \H_c^m(X,\cL) \times {\Ext}^{1-m}_X(\cL,\cD_X) \lra \H_c^1(X,\cD_X) \os{\tr_X}\lra \qz \]
induced by Yoneda pairing is a non-degenerate pairing of finite groups.
\end{enumerate}
\end{thm}
The proof of this theorem will cover the following subsections and will be completed in subsection 1.5.
\smallskip
\subsection{Artin-Verdier duality}\label{sect4-1}
We review the Artin-Verdier duality for number fields (cf.\ \cite{arver}, \cite{mazur}, \cite{milne:adual} II.2--3).
Let $\Gm:=\Gm{}_{,S}$ be the sheaf on $S_\et$ given by the multiplicative group. By global class field theory, we have
\begin{equation}\label{isom:cft}
  \H_c^m(S,\Gm) \cong
 \begin{cases}
  \qz \qquad & \hbox{($m=3$)}\\
  0 \qquad & \hbox{($m=2$ or $m \ge 4$)}.
 \end{cases}
\end{equation}
We normalize the isomorphism for $m=3$ as follows. For a closed point $y$ of $S$, let $G_y$ be the absolute Galois group of $\kappa(y)$, and let
\[ \tr_{y,\qz}:\H^1(y,\qz) \lra \qz \]
be its trace map, i.e., the unique homomorphism that evaluates a continuous character $\chi \in \Hom_{\cont}(G_y,\qz) =\H^1(y,\qz)$ at the arithmetic Frobenius substitution $\varphi_y \in G_y$.
Then for any closed point $i_y:y \hra S$ of $S$ the composition
\[\xymatrix{
\H^1(y,\qz) \ar[r]^-{\delta} & \H^2(y,\Z) \ar[rr]^-{\gys_{i_y,\Gm}} && \H_c^3(B,\Gm) \ar[r]^-{\eqref{isom:cft}} & \qz
}\]
coincides with $\tr_{y,\qz}$, where $\gys_{i_y,\Gm}$ denotes the Gysin map $\Z[-1] \ra Ri_y^!\Gm$ defined in \cite{sga4.5} Cycle 2.1.1 (see also Proposition \ref{prop:Kummer}\,(1) below), and the map $\delta$ is the connecting homomorphism associated with the short exact sequence
\[ 0 \lra \bZ \lra \bQ \lra \qz \lra 0. \]
The Artin-Verdier duality shows that for an integer $m$ and a constructible sheaf $\cF$ on $S_\et$, the pairing
\begin{equation}\label{dual_Gm}
\xymatrix{
\H_c^m(S,\cF) \times \Ext^{3-m}_S(\cF,\Gm) \ar[r] & \H_c^3(S,\Gm) \ar[r]^-{\eqref{isom:cft}} & \qz
}
\end{equation}
induced by Yoneda pairing is a non-degenerate pairing of finite groups.

\smallskip
\subsection{Artin-Verdier duality revisited}\label{sect4-2'}
We formulate a version of Artin-Verdier duality replacing $\Gm$ by a complex of torsion sheaves.
Let $n$ be a positive integer, let $p$ be a prime, and let $T$ be any scheme. For any \'etale sheaf $\cG$ on $T$ let
$$
\cG_n = \sHom_T(\zn,\cG) = \ker(\cG \os{n} \lra \cG)
$$
be the subsheaf of sections annihilated by
$n$, let $\cG_{p^\infty} = \varinjlim{}_{n \ge 1} \cG_{p^n}$ be the subsheaf of sections annihilated
by a power of $p$ (also called the $p$-primary torsion subsheaf of $\cG$), and let $\cG_\tor = \varinjlim{}_{n \ge 1} \ \cG_n
= \bigoplus_p \ \cG_{p^\infty}$ be the torsion subsheaf of $\cG$.

Denote by $\Sh(T_\et)$, $\Sh(T_\et,\zn)$, $\Sh(T_\et,p^\infty)$ and $\Sh(T_\et,\tor)$
the categories of \'etale sheaves, \'etale $\zn$-sheaves, \'etale $p$-primary torsion sheaves
and \'etale torsion sheaves on $T$, respectively. The exact inclusion functors
\begin{equation}\label{inclusion.functors}
\Sh(T_\et,\zn) \hookrightarrow \Sh(T_\et), \quad  \Sh(T_\et,p^\infty) \hookrightarrow \Sh(T_\et),
\quad \Sh(T_\et,\tor) \hookrightarrow  \Sh(T_\et)
\end{equation}
have the left exact right adjoints
\begin{equation}\label{torsion.sheaves}
\cG \mapsto \cG_n \quad , \quad \cG \mapsto \cG_{p^\infty} \quad , \quad  \cG \mapsto \cG_\tor \,,
\end{equation}
respectively, and these functors derive triangulated functors
\begin{equation}\label{torsion.sheaves.derived}
D^+(T_\et) \rightarrow D^+(T_\et,\zn), \quad D^+(T_\et) \rightarrow D^+(T_\et,p^\infty), \quad D^+(T_\et) \rightarrow D^+(T_\et,\tor)
\end{equation}
of the corresponding derived categories, which we denote by the same symbol.
The functors \eqref{torsion.sheaves} preserve injectives, because their left adjoints \eqref{inclusion.functors} are exact.
Hence we have
\begin{equation}\label{Rhom}
R\sHom_{S,\zpn}(\bullet,?_n) = R\sHom_S(\bullet,?)
\end{equation}
as bifunctors from $D^-(S_\et,\zn)^\op \times D^+(S_\et)$ to $D^+(S_\et,\zn)$ by \cite{sga4.5} Ca\-t\'e\-go\-ries D\'e\-ri\-v\'ees II.1.2.3\,(3).
The analogous result holds for the functors $?_{p^\infty}$ and $?_\tor$.
In particular (taking zero degree sections), the functors \eqref{torsion.sheaves.derived}
are right adjoint to the natural functors in the opposite
direction which are induced by the functors \eqref{inclusion.functors}.

\smallskip
We note that, for $\cG$ in $D^+(T_\et)$, the objects $\cG_n$, $\cG_{p^\infty}$ and $\cG_\tor$ can
also be regarded in $D^+(T_\et)$, and that the adjunctions gives canonical morphisms $\cG_n \os{\iota}\lra \cG$
in $D^+(T_\et)$ and canonical factorizations for them, for positive integers $n\mid n'$ and primes $p$,
$$
\cG_n \rightarrow \cG_{n'} \rightarrow \cG_\tor \rightarrow \cG \quad \mbox{ and } \quad \cG_{p^n} \rightarrow \cG_{p^\infty} \rightarrow \cG_\tor \rightarrow \cG\,,
$$
which induce isomorphisms $\cG_n \cong (\cG_\tor)_n \cong \cG_n$ and $\cG_{p^n} \cong (\cG_{p^\infty})_{p^n} \cong \cG_{p^n}$.

\begin{cor}\label{cor.Ext}
Assume that $T$ is noetherian. For $\cL \in D^-(T_\et)$ with constructible torsion cohomology sheaves
and $\cG$ in $D^+(T_\et)$, the morphism $\cG_\tor \rightarrow \cG$ induces a functorial isomorphism
\[ R\sHom_T(\cL,\cG_\tor) = R\sHom_T(\cL,\cG) \]
\end{cor}

\begin{pf}
We may assume that $\cL$ is a constructible torsion sheaf, by a standard argument using spectral sequences.
Then by the constructabilty we may further assume that $\cL$ is annihilated by some positive integer $n$,
and then the claim follows from \eqref{Rhom}.
\end{pf}

\begin{lem}\label{Kummer.sequence}
Let $\cG$ be an object in $D^+(T_\et)$.
\begin{itemize}
\item[(1)]
For any positive integer $n$ there is a canonical distinguished triangle
\[
\xymatrix{ \cG_{n} \ar[r]^-{\iota} & \cG \ar[r]^-{\times n} & \cG \ar[r]^-{\delta^\tr_{\cG}} & \cG_{n}[1]\,.}
\]
\item[(2)]
If $T$ is quasi-compact, then one has isomorphisms
\[
H^i(T,\cG_{p^\infty}) \cong \varinjlim{}_{n \ge 1} H^i(T,\cG_{p^n}) \quad \mbox{ and }
\quad H^i(T,\cG_\tor) \cong \varinjlim{}_{n \ge 1} H^i(T,\cG_n)\,.
\]
\end{itemize}
\end{lem}

\begin{pf}
The exact sequence of sheaves $0 \to \bZ \os{\times n}\lra \Z \to \zn \to 0$ induces a canonical distinguished triangle
\begin{equation}\label{KummerZ}
\xymatrix{ \bZ \ar[r]^-{\times n} & \bZ \ar[r]^-{\can} & \zn \ar[r]^-{\delta^\tr_\Z} & \bZ[1]\,. }
\end{equation}
Claim (1) folllows from applying the exact functor $R\sHom(-,\cG)$ to this triangle.
Alternatively, if $\cG$ is represented by a bounded below complex $I$ with injective components,
then $\cG_n$ is represented by the complex $I_n$, and the distinguished triangle in (1) is
represented by the exact sequence of complexes
\[
\xymatrix{ 0  \ar[r] & I_n \ar[r]^-{\iota} & I \ar[r]^-{\times n} & I \ar[r] &  0 \;.}
\]
In these terms, $\cG_\tor$ is represented by the complex $I_\tor$, and the second claim in (2) follows from the fact
that $I = \varinjlim{}_{n \ge 1} I_n$ and that the cohomology commutes with filtered direct limits for a quasi-compact
scheme $T$. The proof for the first claim in (2) is analogous.
\end{pf}

\smallskip
Applying these results to the sheaf $\Gm$ we get the following variant of Artin-Verdier duality.
To distinguish from the sheaves $(\Gm)_n$, $(\Gm)_{p^\infty}$ and $(\Gm)_\tor$, we write $R(\Gm)_n$,
$R(\Gm)_{p^\infty}$ and $R(\Gm)_\tor$ for the corresponding objects in the derived category.

\begin{cor}\label{cor:AVvariant}
\begin{enumerate}
\item[(1)]
There is a canonical trace isomorphism
\[ \tr_S: H^3_c(S,R(\Gm)_\tor) \isom \qz\,. \]
\item[(2)]
For $\cL \in D^b(S_\et)$ with constructible torsion cohomology sheaves, the pairing
\[ \H_c^m(S,\cL) \times \Ext^{3-m}_S(\cL,R(\Gm)_\tor) \lra H^3_c(S,R(\Gm)_\tor) \os{\tr_S}\lra \qz \]
induced by Yoneda pairing is a non-degenerate pairing of finite groups.
\end{enumerate}
\end{cor}
\begin{pf}
By Lemma \ref{Kummer.sequence}\,(1), we have a long exact sequence
\begin{align*}
\dotsb & \lra \H_c^m(S,R(\Gm)_n) \lra \H_c^m(S,\Gm) \os{\times n}{\lra} \H_c^m(S,\Gm) \\
 & \lra \H_c^{m+1}(S,R(\Gm)_n) \lra \dotsb .
\end{align*}
By \eqref{isom:cft} and this exact sequence, we obtain $\H_c^m(S,R(\Gm)_n)=0$ for $m \geq 4$ and a trace isomorphism
\[ \tr_{S,n} : \H_c^3(S,R(\Gm)_n) \isom \zn\,. \]
We get the trace isomorphism in (1) by passing to the limit on $n \ge 1$. In fact, \ref{Kummer.sequence} (2)
easily extends to the cohomology with compact support $H^i_c(S,-)$.
\par
The claim (2) follows from the non-degeneracy of \eqref{dual_Gm} and Proposition \ref{cor.Ext}.
\end{pf}

\smallskip
\subsection{Kummer theory}\label{sect4-2}
We establish canonical isomorphisms
\[ R(\Gm)_{p^n} \cong \zpn(1)' \qquad \hbox{($n \in \bN$ \, or \, $n=\infty$)} \]
where we let $\bZ/p^\infty(1)' = \qzp(1)'$, by definition.
\begin{prop}\label{prop:Kummer}
Let $p$ be a prime number and let $n$ be a positive integer or $\infty$. Let $i$ be the closed immersion $Y:= S \times_\Z \bF_p \hra S$,
and let $j$ be the open immersion of the complement $U:=S[p^{-1}] \hra S$. Finally let $R(\Gm)_{p^n}$ be as above. Then{\rm:}
\begin{enumerate}
\item[(1)]
There is a canonical isomorphism $\beta_U: \mupn \isom j^* R(\Gm)_{p^n}$ on $U_\et$.
\item[(2)]
For any closed subscheme $i_Z: Z \hra S$ of codimension $1$ there are canonical Gysin isomorphisms on $Z_\et$
\[ \gys_{i_Z,\Gm} : \bZ[-1] \isom Ri_Z^!\Gm \quad \text{and}\quad \gys_{i_Z,p^n} : \zpn[-2] \isom Ri_Z^! R(\Gm)_{p^n}\,. \]
\item[(3)] There is a unique isomorphism $\beta: \zpn(1)'_S \isom R(\Gm)_{p^n}$ completing the following diagram to an isomorphism of distinguished triangles
\stepcounter{equation}
\begin{equation}\label{eq4-2-1}
 \xymatrix{
i_*\zpn[-2] \ar[r]^-g \ar[d]_{i_*(\gys_{i,p^n})}^{\rwr}
& \zpn(1)'_S \ar[r]^-t \ar[d]_\beta^{\rwr}
& Rj_*\mupn \ar[rr]^-{\delta_S^\val} \ar[d]_{Rj_*(\beta_U)}^{\rwr}
&& i_*\zpn[-1] \ar[d]_{i_*(\gys_{i,p^n})[1]}^{\rwr} \\
i_* Ri^! R(\Gm)_{p^n} \ar[r]^-{i_*}
& R(\Gm)_{p^n} \ar[r]^-{j^*} & Rj_* j^* R(\Gm)_{p^n} \ar[rr]^-{-\delta^\loc(R(\Gm)_{p^n})}
&& i_* Ri^! R(\Gm)_{p^n}[1]\,.\hspace{-5pt}}
\end{equation}
Here the top triangle comes from the definition of $\zpn(1)'_S$, and the bottom triangle from the localization sequence \eqref{eq-0.4.2}
for $R(\Gm)_{p^n}$.
\item[(4)] For $n<\infty$ there is a canonical distinguished triangle {\rm(}in $D^b(S_\et)${\rm)}
\[ \zpn(1)'_S \os{\gamma}\lra \Gm \os{\times p^n}\lra \Gm \lra \zpn(1)'_S[1] \,. \]
\end{enumerate}
\end{prop}
\begin{pf}
First let $n < \infty$. Applying Lemma \ref{Kummer.sequence} (1) to $\cG = \Gm$, we get a canonical distinguished triangle
\begin{equation}\label{KummerGm}
\xymatrix{ R(\Gm)_{p^n} \ar[r]^-{\iota} & \Gm \ar[r]^-{\times p^n} & \Gm \ar[r]^-{\delta^\tr_{\Gm}} & R(\Gm)_{p^n}[1]\,.}
\end{equation}
On the other hand, since $p$ is invertible on $U$, there is an exact sequence
\begin{equation}\label{KummerGmU}
\xymatrix{ 0 \ar[r] & \mupn \ar[r] & \Gm{}_{,U} \ar[r]^-{\times p^n} & \Gm{}_{,U} \ar[r] & 0 \,. }
\end{equation}
This gives canonical isomorphisms $j^* R(\Gm)_{p^n}=R\sHom_U(\zpn,\Gm{}_{,U}) \cong \mupn$, as claimed.
Since we have sheaves on both sides we also get a canonical isomorphism $j^\ast R(\Gm)_{p^\infty} \cong \mu_{p^\infty}$
by passing to the direct limit.

(2) First one notes that
\[ R^mi_Z^!\Gm \cong \begin{cases}
                            \Z \qquad & \hbox{($m=1$)}\\
                             0 \qquad & \hbox{($m \ne 1$)}
                         \end{cases} \]
(cf.\ e.g., \cite{milne:adual} p.\ 185, bottom). Therefore
\[ \Hom_Z(\bZ[-1],Ri_Z^!\Gm) = H^1_Z(S,\Gm) \cong \bigoplus_{z \in Z} \ \bZ \,, \]
and to get a canonical isomorphism $\gys_{i_Z,\Gm}$ it suffices to replace $Z$ by a point $z\in Z$ and to find a canonical generator of $H^1_z(S,\Gm)$. This is done by the localization sequence
\[ \cO_{S,z}^\times \lra k^\times \os{\delta}\lra H^1_z(S,\Gm) \lra H^1(k,\Gm)=0 \]
for the discrete valuation ring $\cO_{S,z}$. Now we take $\delta(\pi)$ as a generator for any prime element $\pi$ of $\cO_{S,z}$.
\par
As for the second Gysin isomorphism in (2), consider a diagram on $Z_\et$
\begin{equation}\label{betaY}
\xymatrix{
 \bZ[-1] \ar[rr]^-{\times p^n} \ar[d]_{\gys}^{\rwr} \ar@{}[rrd]|{\text{(}*\text{)}}
&& \bZ[-1] \ar[rr]^-{\can} \ar[d]_{\gys}^{\rwr}
&& \zpn[-1] \ar[rr]^-{-\delta^\tr_\Z\,[-1]} \ar@{.>}[d]_{\beta_Z}
&& \bZ \ar[d]_{\gys[1]}^{\rwr} \\
Ri_Z^!\Gm \ar[rr]^-{\times p^n}
&& Ri_Z^!\Gm \ar[rr]^-{-Ri_Z^!(\delta^\tr_{\Gm})}
&& Ri_Z^!R(\Gm)_{p^n} [1] \ar[rr]^-{Ri_Z^!(\iota)[1]}
&& Ri_Z^!\Gm \,,\hspace{-5pt}}
\end{equation}
where $\gys$ denotes $\gys_{i_Z,\Gm}$. The top sequence is a distinguished triangle by \eqref{KummerZ} and the rule recalled in \S\ref{sect0-5-1}. The bottom distinguished triangle is obtained by applying $Ri^!$ to \eqref{KummerGm} and shifting suitably.
Now the commutativity of the square $(*)$ implies the existence of a morphism $\beta_Z$ making the diagram commutative (cf.\ \S\ref{sect0-5-2}), which then necessarily is an isomorphism. Moreover, since
\[ \Hom_{D(S_\et)}(\zpn[-1],Ri_Z^!\Gm) \cong \Hom_{D(S_\et)}(\zpn,\Z[-1]) = 0
 \quad \hbox{(cf.\ \S\ref{sect0-5-4}\,(2)),} \]
such $\beta_Z$ is unique by Lemma \ref{lem.derived}\,(1). So $\gys_{i_Z,p^n} := \beta_Z$ gives the desired canonical isomorphism.
\par (The sign $-1$ on $Ri_Z^!(\delta^\tr_{\Gm})$ is motivated by the fact that the restriction $(-\delta^\tr_{\Gm})|_U$ is the connecting morphism $\Gm{}_{,U} \to \mupn[1]$ associated with the short exact sequence \eqref{KummerGmU}, which appears in the definition of Deligne's cycle class \cite{sga4.5} Cycle. In particular, by our choice, $\gys_{i_Z,p^n}$ agrees with the Gysin morphism in \S\ref{sect1-1} when $Z$ is contained in $U$.)
\par
(3) Since $\zpn(1)'_S$ is concentrated in $[0,1]$, we have
\[ \Hom_{D(S_\et)}(\zpn(1)'_S,i_* Ri^!R(\Gm)_{p^n}) \cong \Hom_{D(S_\et)}(\zpn(1)'_S,i_*\zpn[-2]) = 0\,. \]
In view of Lemma \ref{lem.derived}\,(1) and the fact that $\beta_U$ and $\gys_{i,p^n}$ are isomorphisms, our task is to show that the right hand square of \eqref{eq4-2-1} is commutative. Since this comes down to morphisms of sheaves (see \ref{eq-defS}), it suffices to show this for $n<\infty$,
and the case $n=\infty$ follows by passing to the inductive limit.
For $n<\infty$ there is a commutative diagram of distinguished triangles
\begin{equation}\label{9diagram}
\xymatrix{
i_* Ri^!R(\Gm)_{p^n} \ar[r]^-{i_*} \ar[d]_{\iota}
& R(\Gm)_{p^n} \ar[r]^-{j^*} \ar[d]_{\iota}
& Rj_* j^*R(\Gm)_{p^n} \ar[rr]^-{-\delta^\loc(R(\Gm)_{p^n})} \ar[d]_{\iota}
&& \\
i_* Ri^!\Gm \ar[r]^-{i_*} \ar[d]_{\times p^n}
& \Gm \ar[r]^-{j^*} \ar[d]_{\times p^n}
& Rj_* j^*\Gm \ar[rr]^-{-\delta^\loc(\Gm)} \ar[d]_{\times p^n}
&& \\
i_* Ri^!\Gm \ar[r]^-{i_*} \ar[d]_{\delta^\tr_{\Gm}}
& \Gm \ar[r]^-{j^*} \ar[d]_{\delta^\tr_{\Gm}}
& Rj_* j^*\Gm \ar[rr]^-{-\delta^\loc(\Gm)} \ar[d]_{\delta^\tr_{\Gm}}
&&\,, \\
& & &&
}
\end{equation}
where the columns are the distinguished triangles coming from \eqref{KummerGm}, and the rows are localization triangles. We now obtain the following diagram of sheaves on $S_\et$:
\[ \xymatrix{
R^1j_*\mupn \ar[rrrr]^{\delta^\val_S} \ar[dr]^{\beta_U} \ar@{}[rrrrd]|{\text{\scriptsize(a)}}
& & & & i_*\zpn \ar[dl]_{\gys_{i,p^n}} \\
& R^1j^* R(\Gm)_{p^n} \ar[rr]^{\delta^\loc(R(\Gm)_{p^n})} \ar@{}[rrd]|{\text{\scriptsize(c)}}
& & i^* R^2i^! R(\Gm)_{p^n} & \\
& j_*j^*\Gm \ar[u]^{-\delta^\tr_{\Gm}} \ar[rr]_{\delta^\loc(\Gm)}_{\raisebox{-0,8cm}{\text{\scriptsize(e)}}}
& & i_* R^1i^! \Gm \ar[u]_{-\delta^\tr_{\Gm}}_{\hspace{1,2cm}\raisebox{-0,25cm}{\text{\scriptsize(d)}}} \\
j_*\Gm{}_{,U} \ar@{->>}[uuu]^{-\delta^\tr_{\Gm}}_-{\hspace{.7cm}\raisebox{-0,25cm}{\text{\scriptsize(b)}}} \ar@{=}[ur] \ar[rrrr]_{\ord} & & & & i_*\bZ\,. \hspace{-5pt} \ar[uuu]_{\can}\ar[ul]^{\gys_{i,\Gm}}
}\]
The middle square (c) with the four $\delta$'s comes from diagram \eqref{9diagram} and anti-commutes, because $\delta^\loc$ is functorial for the morphism $\delta^\tr_{\Gm}: \Gm \to R(\Gm)_{p^n}[1]$ and we have
\[ \delta^\loc(R(\Gm)_{p^n}[1]) = - \delta^\loc(R(\Gm)_{p^n})[1] \qquad \hbox{(cf.\ \eqref{eq-0.4.1})}. \]
The top arrow $\delta_S^\val$ is induced by residue maps, so the outer square of the diagram commutes by the remark after the proof of (2). The diagram (b) commutes by the definition of $\beta_U$, and the diagram (d) commutes by the definition of $\gys_{i,p^n}$, i.e., by the commutativity of the diagram \eqref{betaY}. The bottom arrow is induced by the normalized discrete valuations for the points $y \in Y$, and the diagram (e) commutes by the definition of the Gysin map $\gys_{i,\Gm}$.
Consequently the diagram (a) anti-commutes, and the right hand square of \eqref{eq4-2-1} commutes by \S\ref{sect0-5-4}\,(1).
Thus we obtain (3).

Finally (4) follows from \eqref{KummerGm} and the isomorphism $\beta$ in (3) by letting $\gamma := \iota \circ \beta$.
\end{pf}

\subsection{Proof of the duality}\label{sect4-5}

By using the canonical isomorphisms $R(\Gm)_{p^\infty}\cong \qzp(1)'$ from Proposition \ref{prop:Kummer}
and the deduced isomorphism
$$
R(\Gm)_\tor=\oplus_p R(\Gm)_{p^\infty} \cong \oplus_p \qzp(1)' = \qz(1)'\,,
$$
we immediately get the following from Corollary \ref{cor:AVvariant}

\begin{cor}\label{cor:AVvariant2}
\begin{enumerate}
\item[(1)]
There is a canonical trace isomorphism
\[ \tr_S: H^3_c(S,\qz(1)') \isom \qz\,. \]
\item[(2)]
For $\cL \in D^b(S_\et)$ with constructible torsion cohomology sheaves, the pairing
\[ \H_c^m(S,\cL) \times \Ext^{3-m}_S(\cL,\qz(1)') \lra H^3_c(S,\qz(1)') \os{\tr_S}\lra \qz \]
induced by Yoneda pairing is a non-degenerate pairing of finite groups.
\end{enumerate}
\end{cor}

With this we are now ready to prove Theorem \ref{thm:global-duality}:

(1) Let $\tr_f: Rf_!\cD_X=Rf_!Rf^!\cD_S \to \cD_S = \qz(1)'[2]$ be the canonical trace map, i.e., the adjunction morphism for the adjunction between $Rf^!$ and $Rf_!$ (\cite{sga4} XVIII.3.1.4). We then define the trace map $\tr_X$ as the composite
\[ \tr_X: \H_c^1(X,\cD_X) \os{\tr_f}\lra \H_c^1(S,\cD_S) = H^3_c(S,\qz(1)'_S) \os{\tr_S}\isom \qz\,. \]
(2) There is a commutative diagram of Yoneda pairings
\[ \xymatrix{
H^m_c(X,\cL) \ar@<19pt>@{=}[d] \hspace{-40pt}
& \times \; \Ext^{1-m}_X(\cF,Rf^!\D_S) \ar@<7pt>@{=}[d] \ar[r]
& H^1_c(X,Rf^!\D_S) \ar[d]^{\tr_f} \\
H^m_c(B,Rf_!\cL) \hspace{-28pt}
& \times \; \Ext^{1-m}_S(Rf_!\cL,\D_S) \ar[r]
& H^1_c(S,\D_S)\,,
}\]
and the assertion follows from Corollary \ref{cor:AVvariant2}.

\newpage
\section{The case where $p$ is invertible on the scheme $X$}\label{sect1}
\medskip
In this section, we work in the following setting. Let $X$ be a noetherian excellent regular scheme, let $n$ be a positive integer invertible on $X$ and put $\L:=\bZ/n$. For integers $q \in \bZ$, put $\L(q):=\mu_n^{\otimes q}$, the $q$-fold Tate twist of the \'etale sheaf $\L$ on $X$ or $X$-schemes. Let $Z \subset X$ be a regular closed subscheme of pure codimension $c$. By Gabber's construction \cite{fujiwara} 1.1.2, there is a cycle class $\cl_X(Z)$ in the \'etale cohomology group $\H^{2c}_Z(X,\L(c))$ (without using the absolute purity), which satisfies the following three properties:
\begin{enumerate}
\item[(G1)] For an \'etale morphism $X' \ra X$ and $Z':=Z \times_X X'$, the pull-back of $\cl_X(Z)$ to $\H^{2c}_{Z'}(X',\L(c))$ agrees with $\cl_{X'}(Z')$.
\item[(G2)] For regular closed subschemes $Z \subset Y \subset X$, we have $\cl_X(Y) \cap \cl_{Y}(Z) = \cl_X(Z)$ in $\H^{2c}_{Z}(X,\L(c))$.
\item[(G3)] The image of $\cl_X(Z)$ in $\H^0(Z,R^{2c}i^!\L(c))$ coincides with Deligne's cycle class \cite{sga4.5} Cycle \S2.2. Here $i=i_Z$ denotes the closed immersion $Z \hra X$.
\end{enumerate}
\par
\smallskip
\subsection{Gysin maps and compatibility}\label{sect1-1}
For $q,r \in \bZ$, one defines the Gysin map $\gys_i$ as
\[ \gys_i: \H^q(Z,\L(r)) \lra \H^{q+2c}_Z(X,\L(r+c)), \quad  \alpha \mapsto \cl_X(Z) \cup \alpha\,. \]
The main aim of this section is the following compatibility result:
\begin{thm}\label{prop:appA}
Let $c$ be a positive integer, and let $i_x:x \hra X$ and $i_y: y \hra X$ be points on $X$ of codimension $c$ and $c-1$, respectively, with $x \in \ol {\{y\}}$. Then the following square commutes for integers $q, r \ge 0${\rm:}
\stepcounter{equation}
\begin{equation}\label{CD:app1}
\xymatrix{
\H^{q+1}(y,\L(r+1)) \ar[rr]^-{-\dKyx} \ar[d]_{\gys_{i_y}} && \H^q(x,\L(r)) \ar[d]^{\gys_{i_x}} \\
\H^{q+2c-1}_y(X,\L(r+c)) \ar[rr]^-{\delta^\loc_{y,x}(\L(r+c)_X)} && \H^{q+2c}_x(X,\L(r+c))\,.\hspace{-5pt}
}\end{equation}
\end{thm}
To prove the theorem, we may assume that
\begin{equation}\label{eq-local-ell}
\hbox{ $X$ is local with closed point $x$.}
\end{equation}
Put $Z:=\ol {\{ y \}} \subset X$, which has dimension $1$ and consists of two points $\{ y, x \}$. Let $i_Z$ (resp.\ $\iota_x$) be the closed immersion $Z \hra X$ (resp.\ $x \hra Z$). The proof will be finished in \S\ref{sect1-4} below.

\smallskip
\subsection{Regular case}\label{sect1-2}
We first prove Theorem \ref{prop:appA}, assuming that $Z$ is regular. In this case $Z$ is the spectrum of a discrete valuation ring $A$, and we have the cycle classes
\[ \cl_X(Z) \in \H^{2(c-1)}_Z(X,\L(c-1)) \qquad \hbox { and } \qquad \cl_Z(x) \in \H^{2}_x(Z,\L(1)) \]
by Gabber's construction, where $\cl_Z(x)$ agrees with Deligne's construction in \cite{sga4.5} Cycle \S2.2 by (G3). There is a  diagram of boundary maps
\[\xymatrix{
 & \H^{q+1}(y,\L(r+1)) \ar[r]^-{\gys_{i_y}} \ar[ld]_{-\dKyx} \ar[d]^{\delta^{\loc}_{y,x}(\L(r+1)_Z)}  & \H^{q+2c-1}_y(X,\L(r+c)) \ar[d]^{\delta^{\loc}_{y,x}(\L(r+c)_X)} \\
\H^q(x,\L(r)) \ar[r]_-{\gys_{\iota_x}} & \H^{q+2}_x(Z,\L(r+1)) \ar[r]_-{\gys_{i_Z}} & \H^{q+2c}_x(X,\L(r+c))\,.
}\]
Here $\gys_{i_y}$ is the map taking the cup-product with $\cl_X(Z) \vert_{\Spec(\O_{X,y})}$, by the property (G1), and hence the right square commutes by the naturality of cup products and \eqref{eq-0.4.1}. The composite of the bottom row agrees with $\gys_{i_x}$ by (G2). Thus we obtain the commutativity of the diagram \eqref{CD:app1}, once we show the left triangle commutes. But this commuting follows from (G3) and \cite{sga4.5} Cycle 2.1.3. Indeed, by noting that
\[ R^q j_*\L(r+1)=0 \quad \hbox{ for $q \ge 2$, where $j : y \hra Z$,} \]
the left triangle is induced by the following square in $D^b(x_\et,\L)$:
\[\xymatrix{
\iota_x^*R^1 j_*\L(r+1)[-1] \ar[d]_{-\dval[-1]} && \tau_{\leq 1}\,\iota_x^*R j_*\L(r+1) \ar[ll]_-{\text{canonical}}  \ar[d]^{\delta^{\loc}_{y,x}(\L(r+1)_Z)} \\
\L(r) [-1] \ar[rr]^-{\gys_{\iota_x}} && R \iota_x^!\L(r+1)[1]\,,
}\]
where $\dval: \iota_x^*R^1 j_*\L(r+1) \to \L(r)$ denotes a map of sheaves on $x_\et$ induced by the valuation of $A$. We note $R i^!\L(r+1)[1]$ is concentrated in degree $1$. Therefore it suffices to show its commutativity after taking the cohomology sheaves $\cH^1(-)$ in degree 1, so that we are reduced to showing the commutativity of the diagram
\[ \xymatrix{
\iota_x^*R^1 j_*\L(r+1) \ar[d]_{-\dval} & \iota_x^*R^1j_*\L(r+1) \ar@{=}[l]  \ar[d]^{\delta^{\loc}_{y,x}(\L(r+1)_Z)} \\
\L(r) \ar[r]^-{\gys_{\iota_x}} & R^2 \iota_x^!\L(r+1)\,,
}\]
where $\gys_{\iota_x} : \L(r) \to R^2 i^!\L(r+1)$ is given by $a \mapsto \cl_Z(x)\cup a$. By looking at the stalks, we are now reduced to the case that $A$ is strictly henselian and to showing the anti-commutativity of
\[\xymatrix{
\kappa(y)^\times/n \ar[d]_{\ord_A} \ar[rr]^-{h^1}_-{\sim} && H^1(y,\L(1)) \ar[d]^{\delta^{\loc}_{y,x}(\L(r+1)_Z)} \\
\L \ar[rr]^-{\gys_{\iota_x}}_-{1\,\lmt \,\cl_Z(x)} && H^2_x(Z,\L(1))\,,
}\]
which is a consequence of \cite{sga4.5} Cycle 2.1.3 ($h^1$ is the Kummer isomorphism). Note that we have $\dval = \ord_A \circ (h^1)^{-1}$ by \S\ref{sect0-6}\,(I.1), and that $\gys_{\iota_x} \circ \ord_A$ sends a prime element $\pi$ of $A$ to $\cl_Z(x)$ and hence agrees with the map induced the composition
\[ \kappa(y)^\times \os{\delta}\lra H^1_x(Z,\Gm) \os{\delta}\lra H^2_x(Z,\L(1))  \qquad \hbox{(loc.\ cit.\ 2.1.2)}. \]

\smallskip
\subsection{General case}\label{sect1-3}
We prove Theorem \ref{prop:appA} in the case that $Z$ is not regular. Take the normalization $f:T \ra Z$. Put $\varSigma:=f^{-1}(x) \subset T$ with the reduced subscheme structure, which is the set of all closed points on $T$ by the assumption \eqref{eq-local-ell}. Note that $T$ is regular and that $f$ is finite by the excellence of $Z$. Since a finite morphism is projective (\cite{ega2} 6.1.11), the map $T \ra X$ factors as $T \hra \P^d_X \ra X$ for some integer $d \ge 1$. Let $\iota_y$ be the composite map $y \hra T \hra \P^d_X$.
There is a commutative diagram of schemes
\begin{equation}\label{CD:schemes}
\xymatrix{
\varSigma \; \ar@<-1pt>@{^{(}->}[r] \ar[d]_h \ar@/^7mm/[rr]^{i_{\varSigma}}
& T \; \ar@<-1pt>@{^{(}->}[r]^-{i_{T}} \ar[d]_f
& \P:=\P^d_X \ar[d]_g \\
x \; \ar@<-1pt>@{^{(}->}[r]^-{\iota_x} \ar@/_7mm/[rr]_{i_x}
& Z \; \ar@<-1pt>@{^{(}->}[r]^-{i_Z}
& X \,, \hspace{-5pt}
}
\end{equation}
where $i_x$, $i_{\varSigma}$ and all horizontal arrows are closed immersions cf.\ \eqref{eq-local-ell}. Now put $c':=c+d=\codim_{\P}(\varSigma)$, and consider the diagram in Figure 1.
\par
\begin{figure}[htp]
\setlength{\unitlength}{.6mm} {\scriptsize
\begin{picture}
(220,90)(-131,-45) \thinlines
       \put(-99,-35){\vector(2,1){24}}
       \put(67,-35){\vector(-2,1){24}}
       \put(67,36){\vector(-2,-1){24}}
       \put(-99,36){\vector(2,-1){24}}
       \put(-119,34){\line(0,-1){68}}
       \put(-118,34){\line(0,-1){68}}
       \put(84,34){\vector(0,-1){68}}
       \put(-97,40){\vector(1,0){162}}
       \put(-97,-39){\vector(1,0){162}}
       \put(-65,13){\vector(0,-1){26}}
       \put(32,13){\vector(0,-1){26}}
       \put(-35,18){\vector(1,0){36}}
       \put(-35,-18){\vector(1,0){36}}
       \put(-84,-32){{\tiny $\gys_{i_y}$}}
       \put(40,-32){{\tiny $\gys_{i_x}$}}
       \put(36,31){{\tiny $\gys_{i_{\varSigma}}$}}
       \put(-84,31){{\tiny $\gys_{\iota_y}$}}
       \put(86,0){{\tiny $\tr_{h}$}}
       \put(-22,44){{\tiny $-\dval_{y,\varSigma}$}}
       \put(-22,-45){{\tiny $-\dKyx$}}
       \put(-63,0){{\tiny $\alpha$}}
       \put(27,0){{\tiny $\beta$}}
       \put(-32,12){{\tiny $\delta^{\loc}_{y,\varSigma}(\L(r+c'))$}}
       \put(-32,-14){{\tiny $\delta^{\loc}_{y,x}(\L(r+c))$}}
       \put(-18,28){\tiny {\text(1)}}
       \put(-95,0){\tiny {\text(2)}}
       \put(-18,0){\tiny {\text(3)}}
       \put(58,0){\tiny {\text(4)}}
       \put(-18,-30){\tiny {\text(5)}}
\put(70,39){$\H^q(\varSigma,\L(r))$}
\put(-142,39){$\H^{q+1}(y,\L(r+1))$}
\put(-142,-41){$\H^{q+1}(y,\L(r+1))$}
\put(70,-41){$\H^q(x,\L(r))$}
\put(8,17){$\H^{q+2c'}_{\varSigma}(\P,\L(r+c'))$}
\put(-92,17){$\H^{q+2c'-1}_y(\P,\L(r+c'))$}
\put(-92,-20){$\H^{q+2c-1}_y(X,\L(r+c))$}
\put(8,-20){$\H^{q+2c}_x(X,\L(r+c))$}
\end{picture}
\caption{A diagram for the proof in the general case} }
\end{figure}
In the diagram, the arrows $\alpha$ and $\beta$ are induced by the composite morphism
\[ \xymatrix{
\gamma : Rf_*R i_{T}^!\L(r+c')_{\P}[2d] \ar[r]^-{(*)} & Ri_{Z}^!R g_* \L(r+c')_{\P}[2d] \ar[rr]^-{Ri_Z^!(\tr_{g})} && R i_Z^!\L(r+c)_X
}\]
in $D^+(Z_\et,\L)$, where $(*)$ is the cobase-change morphism (cf.\ \cite{sga4} XVIII.3.1.13.2) for the right square of \eqref{CD:schemes}. More precisely, $\alpha$ is obtained by restricting $\gamma$ to $y$, and $\beta$ is defined as the composite
\[ \xymatrix{
\beta : Rh_*Ri_{\varSigma}^!\L(r+c')_{\P}[2d] \ar[r]^-{(**)} & R\iota_x^!Rf_*R i_{T}^!\L(r+c')_{\P}[2d] \ar[r]^-{R\iota_x^!(\gamma)} & R i_x^!\L(r+c)_X\,,
}\]
where $(**)$ is the cobase-change morphism for the left square of \eqref{CD:schemes}. Therefore the square (3) is commutative by \eqref{eq-0.4.1}. On the other hand, the diagram (1) is commutative by the regular case \S\ref{sect1-2}, and moreover, the outer large square of Figure 1 commutes by the definition of $\dKyx$ (cf.\ \S\ref{sect0-6} (II)) and the fact that the trace map $\tr_{h}$ coincides with the corestriction map of Galois cohomology groups (cf.\ \cite{sga4} XVIII.2.9 (Var 4)). Therefore, once we show that the diagrams (2) and (4) commute, we will have obtained the commutativity of (5), i.e., Theorem \ref{prop:appA}. One can easily deduce the commutativity of (2) and (4) from a compatibility result of Riou \cite{riou} 2.34. However we include a simple proof of the commutativity of (4) here for the convenience of the reader (the diagram (2) is simpler and left to the reader).
\smallskip
\subsection{Commutativity of (4)}\label{sect1-4}
Put \[ Q:=\P^d_x \] and take a section $s : X \hra \P$. Let $t:x \hra Q$ be the restriction of $s$ to $x$, and let $i_Q$ (resp.\ $\iota_{\varSigma}$) be the closed immersion $Q \hra \P$ (resp.\ $\varSigma \hra Q$). Consider the following diagram:
\begin{equation}\label{CD:appA2}
\xymatrix{
& \H^q(x,\L(r)) \ar[r]^-{\gys_{i_x}} \ar[d]_{\gys_{t}} \ar@{}[rd]|{\text{(6)}} & \H^{q+2c}_x(X,\L(r+c)) \ar[d]_{\gys_{s}} \\
\H^q(\varSigma,\L(r)) \ar[r]^-{\gys_{\iota_\varSigma}} \ar[rd]_{\tr_h}^{\quad\;\;\text{(7)}} & \H^{q+2d}(Q,\L(r+d)) \ar[r]^-{\gys_{i_Q}} \ar[d]_{\tr_{Q/x}} \ar@{}[rd]|{\text{(8)}} & \H^{q+2c'}_{Q}(\P,\L(r+c')) \ar[d]_{\beta'} \\
& \H^q(x,\L(r)) \ar[r]^-{\gys_{i_x}} & \H^{q+2c}_x(X,\L(r+c))\,.
}\end{equation}
Here the square (6) commutes by (G2). The arrow $\beta'$ is a trace map defined in a similar way as for $\beta$. The diagram (4) in question is related to the large tetragon (7)$+$(8) in \eqref{CD:appA2} by a diagram
\begin{equation}\label{CD:appA2.2}
\xymatrix{
\H^q(\varSigma,\L(r)) \ar[rr]^-{\gys_{i_\varSigma}} \ar[rrd]^{\hspace{-2pt}\varphi}_{\text{(7)}+\text{(8)}\qquad\qquad}  \ar[d]_{\gys_{i_x} \circ \tr_h} && \H^{q+2c'}_{\varSigma}(\P,\L(r+c')) \ar[d]^{\iota_{\Sigma*}} \\ \H^{q+2c}_x(X,\L(r+c)) && \ar[ll]_{\beta'} \H^{q+2c'}_{Q}(\P,\L(r+c'))\,. \hspace{-5pt}
}\end{equation}
Here the arrow $\varphi$ denotes the composite of the middle row of \eqref{CD:appA2}, and the upper right triangle of \eqref{CD:appA2.2} commutes obviously. The composition $\beta'\circ \iota_{\Sigma*}$ is $\beta$, and the square agrees with the diagram (4) in Figure 1. To prove the commutativity of (4), it thus suffices to check that of the lower left triangle of \eqref{CD:appA2.2}, i.e., the tetragon (7)$+$(8) in \eqref{CD:appA2}. To prove this, it suffices to check the following claims concerning the diagram \eqref{CD:appA2}:
\addtocounter{thm}{2}
\begin{lem}\label{lem1-1}
\begin{enumerate}
\item[{\rm (a)}] The triangle {\rm (7)} is commutative.
\item[{\rm (b)}] The triangle {\rm (8)} is commutative.
\end{enumerate}
\end{lem}
\begin{pf}
The claim (a) follows from standard arguments using \cite{sga4.5} Cycle 2.3.8\,(i),\,(ii). To prove prove (b), we use the notation fixed in \eqref{CD:schemes} By the absolute purity \cite{fujiwara}, we have
\[ H^{q+2c}_x(X,\L(r+c)) \simeq H^q(x,R^{2c}i_x^!\L(r+c)), \]
and the problem is reduced to the case that $q=0$ and $\kappa(x)$ is separably closed. Then the Gysin map
\[ \H^0(x,\L(r)) \lra \H^{2d}(Q,\L(r+d)) \]
is bijective. Therefore the diagram (8) commutes by the following facts concerning the diagram \eqref{CD:appA2}:
\begin{enumerate}
\item[{\rm (c)}] the composite of the right vertical column is the identity map,
\item[{\rm (d)}] the composite of the middle vertical column is the identity map,
\end{enumerate}
which are consequences of \cite{sga4.5} Cycle 2.3.8\,(i),\,(ii).
\end{pf}

\noindent
This completes the proof of Theorem \ref{prop:appA}.

\smallskip
\subsection{Bloch-Ogus complexes and Kato complexes}\label{sect1-5}
We use Theorem \ref{prop:appA} to identify the Kato complexes with those defined via the method of Bloch and Ogus. Keep the assumptions as in the beginning: $X$ is a noetherian excellent {\it regular} scheme, $n$ is a positive integer invertible on $X$, and $\L(q)$\,($q \in \bZ$) is the $q$-fold Tate twist of the \'etale sheaf $\L$ on $X$ or $X$-schemes. Assume that $X$ is has pure dimension $d < \infty$, and let $Z \subset X$ be a closed subscheme. For a non-negative integer $q \ge 0$, we define
\begin{equation}
\label{eq1-5-0} (Z/X)_q:=\{x \in Z \ | \ \codim_X(x)=d-q \}. 
\end{equation}
By a standard argument using localization and an exact couple (cf.\ \cite{Bloch-Ogus} 3.10),
there is a localization spectral sequence
\begin{equation}\label{eq1-5-1}
 E^1_{s,t}(Z/X,\L(b)) = \bigoplus_{x\in (Z/X)_s} \ H^{-s-t}_x(X,\L(-b)) \Lra H^{-s-t}_Z(X,\L(-b)),
\end{equation}
which induces a filteration on $H^{-s-t}_Z(X,\L(-b))$ with respect to the dimension of support. It is regarded as the niveau spectral sequence on $Z$ (cf.\ \cite{Bloch-Ogus} (3.7) and \cite{js} 2.7), for the homology theory which is defined on all subschemes $V$ on $X$ by
\[ H_a(V/X,\L(b)) := H^{-a}_V(U,\L(-b))\,, \]
if $V$ is a closed subscheme of an open subscheme $U\subset X$.
By definition, \[ H^a_x(X,\L(-b)) = \varinjlim_{x \in U \subset X} \ H^a_{\ol{\{x\}}\cap U}(U,\L(-b))\,, \] where the limit is taken over all open subsets $U \subset X$ containing $x$.
Hence we have 
\begin{equation}\label{eq1-5-3}
 E^1_{s,t}(Z/X,\L(b)) \cong \bigoplus_{x\in (Z/X)_s} \ H^{s-t-2d}(k(x),\L(s-d-b))
\end{equation}
by the absolute purity \cite{fujiwara}. For a complex of abelian groups $C^*$ denote by $(C^*)^{(-)}$ the complex with the same components, but with the differentials multiplied by $-1$.
\addtocounter{thm}{3}
\begin{thm}\label{thm1-5-1}
The Bloch-Ogus complex $E^1_{*,t}(Z/X,\L(b))$ agrees with the sign-modified Kato complex $C_n^{-t-2d,-d-b}(Z)^{(-)}$ via the Gysin isomorphisms \eqref{eq1-5-3}.
\end{thm}
\begin{pf}
By the construction of the spectral sequence \eqref{eq1-5-1}, its $d^1$-differentials have the components
\[ \delta_{y,x}^{\loc}: H^{-s-t}_y(X,\L(-b)) \lra H^{-s+1-t}_x(X,\L(-b)) \]
for $y\in X_s$ and $x\in X_{s-1}$ with $x\in \overline{\{y\}}$ (cf.\ \cite{js} Remarks 2.8). Therefore the claim directly follows from Theorem \ref{prop:appA}.
\end{pf}
\begin{rem}
For $Z=X$ it is often customary in literature to renumber the spectral sequence \eqref{eq1-5-1} into a cohomological coniveau spectral sequence {\rm(}with $c=-b${\rm)}
\[ E_1^{p,q}(X,\L(c)) = \bigoplus_{x\in X^p} \ H^{q-p}(k(x),\L(c-p)) \Lra H^{p+q}(X,\L(c))\,. \]
See \S{\rm\ref{sect0-4}} for the definition of $X^q$. This does not change the differentials, and so the $E_1$-terms compare in a similar way to the Kato complexes. More precisely one obtains that $E_1^{-*,q}(X,\L(c))$ coincides with $C_n^{q-d,c-d}(X)^{(-)}$.
\end{rem}

By this method, we only get Bloch-Ogus complexes for schemes $Z$ which can be globally embedded in a regular scheme $X$. But the following slight variant covers all the cases considered in \cite{js} p.494, case (i) of 2.B.

Let $S$ be a noetherian excellent regular base scheme of pure dimension $d$, let $n$ be invertible on $S$, and let $b$ be an integer. Similarly as in \S\ref{sect0-3} we define a homology theory (in the sense of \cite{js} 2.1) for all separated $S$-schemes of finite type $f: X \to S$ by defining
\[ H_a(X/S,\L(b)) := H^{-a}(X,Rf^!\L(-b))\,. \]
For a non-negative integer $q \ge 0$, we define a set $(X/S)_q$ as
\addtocounter{equation}{2}
\def\trdeg{\text{tr.deg}}
\begin{equation}\label{eq1-5-3'}
 (X/S)_q:=\{ x \in X \ | \ \codim_S(f(x))-\trdeg(x/f(x)) = d - q  \},
\end{equation}
where for points $x \in X$ and $y \in S$ with $y=f(x)$, we wrote $\trdeg(x/y)$ for the transcendetal degree of $\kappa(x)$ over $\kappa(y)$. The set $(Z/S)_q$  for a closed subscheme $Z \subset S$ agrees with that in the sense of \eqref{eq1-5-0}. When $S$ is the spectrum of a field, $(X/S)_q$ agrees with $X_q$ in the sense of \S\ref{sect0-4}. Under this definition, one gets a niveau spectral sequence
\begin{equation}\label{eq1-5-4}
 E^1_{s,t}(X/S,\L(b)) = \bigoplus_{x\in (X/S)_s} H_{s+t}(x/S,\L(b)) \Lra H_{s+t}(X/S,\L(b))\,,
\end{equation}
where $H_a(x/S,\L(b))$ is defined as the inductive limit of $H_a(V/S,\L(b))$ over all non-empty open subschemes $V\subset \overline{\{x\}}$. Since $\overline{\{x\}}$, being of finite type over $S$, is again excellent, there is a non-empty open subset $V$ which is regular. Then, for all non-empty open $V' \subset V$ one has a canonical isomorphism due to the absolute purity
\begin{equation}\label{eq1-5-5}
 H_a(V'/S,\L(b)) \cong H^{2(s-d)-a}(V',\L(s-b-d))
\end{equation}
with $s=\trdeg(x/f(x)) + d - \codim_S(f(x))$, by the construction in \cite{fujiwara} p.\ 157. This induces an isomorphism
\begin{equation}\label{eq1-5-6}
E^1_{s,t}(X/S,\L(b)) \cong \bigoplus_{x\in (X/S)_s} \ H^{s-t-2d}(x,\L(s-b-d))
\end{equation}
The following theorem generalizes Theorem \ref{thm1-5-1} (which is the case $X=S$).
\addtocounter{thm}{4}
\begin{thm}\label{thm:Katocomplexprimetop}
The Bloch-Ogus complex $E^1_{*,t}(X/S,\L(b))$ agrees with the sign-modified Kato complex $C_n^{-t-2d,-b-d}(X)^{(-)}$ via the isomorphisms \eqref{eq1-5-6}.
\end{thm}
\begin{pf}
The question is local on $X$ and $S$. Thus we can assume that there is a factorization $f=p\circ i$, where $p: P \to S$ is a smooth morphism of relative dimension $N$ (e.g., $P = \bA^N_S$) and $i: X \hra P$ is a closed immersion. Then there is a canonical isomorphism from Poincar\'e duality (\cite{sga4} XVIII 3.2.5) \[ Rp^!\L(-b) \cong \L(N-b)[2N], \] which induces an isomorphism
\begin{align*}
H_a(X/S,\L(b)) &  =  H^{-a}(X,Rf^!\L(-b)) \\
& \cong H^{-a}(X,Ri^!\L(N-b)[2N]) = H^{2N-a}_X(P,\L(N-b))\,.
\end{align*}
Similarly, for a locally closed subset $V\subset P$, say a closed immersion $i_V: V \hra U$ with $j_U: U \hra P$ open, there is an isomorphism
\stepcounter{equation}
\begin{align}
H_a(V/S,\L(b)) & = H^{-a}(V,R(f\circ j_U\circ i_V)^!\L(-b))  \notag \\
 & \cong H^{-a}(V,R(i_V)^!\L(N-b)[2N]) \label{eq1-5-8} \\
 &  = H^{2N-a}_V(U,\L(N-b))\,. \notag
\end{align}
Moreover, this is compatible with localization sequences. If $V$ is regular and of dimension $s$ (hence of codimension $d+N-s$ in $P$), then by definition, the isomorphism \eqref{eq1-5-3} is the composition of this map with the inverse of the Gysin isomorphism
\[ \gys_{i_V}: H^{2s-2d-a}(V,\L(s-d-b)) \isom H^{2N-a}_V(U,\L(N-b))\,. \]
This shows the following: Via the maps \eqref{eq1-5-8}, we get an isomorphism between the homology theory $H_*(-/S,\L(*'))$, restricted to subschemes of $X$, and the homology theory $H_{*-2N}(-/P,\L(*'+N))$ from \eqref{eq1-5-1}, restricted to subschemes of $X$, and therefore an isomorphism of the corresponding spectral sequence. Moreover, via this isomorphisms, the isomorphisms \eqref{eq1-5-3} and \eqref{eq1-5-6} correspond. Therefore the claim follows from Theorem \ref{thm1-5-1}.
\end{pf}

\newpage
\section{The case of $p$-torsion over a perfect field of characteristic $p$}\label{sect2}
\medskip
Throughout this section, $k$ always denotes a perfect field of positive characteristic $p$ and $n$ denotes a positive integer. We will often write $s$ for $\Spec(k)$.
\smallskip
\subsection{Gros' Gysin map}\label{sect2-1}
Let us recall that Gros has defined Gysin morphisms
\[ \gys_f: Rf_*\logwitt Y n r \lra \logwitt X n {r+c}[c] \]
for any proper morphism $f: Y \to X$ of smooth equidimensional varieties over $k$, where $c = \dim(X) - \dim(Y)$ (\cite{gros:purity} II.1). These induce maps
\[ \gys_f: H^q(Y,\logwitt Y n r) \lra H^{q+c}(X,\logwitt X n {r+c})\,. \]
If $i: Y \hra X$ is a closed immersion of smooth $k$-schemes, it also induces Gysin maps
\[ \gys_i: \H^q(Y,\logwitt Y n r) \lra \H^{q+c}_Y(X,\logwitt X n {r+c})\,, \]
where $c$ is now the codimension of $Y$ in $X$. The following result is a $p$-analogue of Theorem \ref{prop:appA}, cf. Remark \ref{rem:smooth}.

\begin{thm}\label{thm:shiho}
Let $X$ be a smooth variety over $k$. Let $n$ and $c$ be positive integers. Let $i_x:x \hra X$ and $i_y:y \hra X$ be points on $X$ of codimension $c$ and $c-1$, respectively, with $x \in \ol {\{ y \}}$. Then the following diagram commutes{\rm:}
\[\xymatrix{
\H^0(y,\logwitt y n {r-c+1}) \ar[rr]^-{(-1)^r \, \dKyx} \ar[d]_{\gys_{i_y}}
&& \H^0(x,\logwitt x n {r-c}) {\ar[d]^{\gys_{i_x}}}\\
\H^{c-1}_y(X,\logwitt X n r) \ar[rr]^-{\delta^{\loc}_{y,x}(\logwitt X n r)}
&& \H^{c}_x(X,\logwitt X n r)\,.\hspace{-5pt} }\]
\end{thm}
In \cite{shm} 5.4, Shiho proved this compatibility property assuming $n=1$, but in a more general situation. The proof of Theorem \ref{thm:shiho} given below relies on the following properties of the Gysin maps:
\begin{enumerate}
\item[(P1)] Local description of Gysin maps. See \cite{gros:purity} II.3.3.9, but we will only need the case of a regular prime divisors, where one can give a simpler proof.
\item[(P2)] Transitivity of Gysin maps \cite{gros:purity} II.2.1.1.
\item[(P3)] For a finite map $h : z \to x$ of spectra of fields which are finitely generated over $k$, the Gysin map $\gys_h: \H^0(z,\logwitt z n r) \to \H^0(x,\logwitt x n r)$ agrees with the corestriction map \eqref{eq-0.6.5}, cf.\ Lemma \ref{lem:norm0} in the appendix.
\end{enumerate}
To prove the theorem, replacing $X$ with $\Spec(\O_{X,x})$, we suppose that
\stepcounter{equation}
\begin{equation}\label{eq-local-p}
\hbox{ $X$ is local with closed point $x$.}
\end{equation}
The proof proceeds in three steps, which will be finished in \S\ref{sect2-4} below.

\smallskip
\subsection{Divisor case}\label{sect2-2}
We first prove Theorem \ref{thm:shiho} assuming $c=1$. In this case, $A:=\cO_{X,x}$ is a discrete valuation ring. Let $\pi$ be a prime element of $A$, and put $K:=\Frac(A)=\kappa(y)$ and $F:=\kappa(x)$. By the Bloch-Gabber-Kato theorem \cite{Bloch-Kato} 2.1, the group $\H^0(y,\logwitt y n r)$ is generated by elements of the forms
\[
\hbox{(i)} \ \ \dlog(\ul{f_1}) \cdot \, \dotsb \, \cdot \dlog(\ul{f_r})
\quad \hbox{ and } \quad
\hbox{(ii)} \ \ \dlog(\ul{\pi}) \cdot \dlog(\ul{f_1}) \cdot \, \dotsb \, \cdot \dlog(\ul{f_{r-1}})\,,
\]
where each $f_j$ belongs to $A^\times$, and for $a \in A$, $\ul a \in W_n(A)$ denotes its Teichm\"uller representative. The diagram in question commutes for elements of the form (i) obviously. We consider the element
\[ \alpha :=  \dlog(\ul {\pi}) \cdot \dlog(\ul {f_1}) \cdot \ \dotsb \ \cdot \dlog(\ul {f_{r-1}}) \]
with each $f_j \in A^\times$, in what follows. By \cite{CTSS} p.\ 779 Lemma 2, we have
\[ \witt X n r/\logwitt X n r \cong \witt X n r/dV^{n-1}\Omega^{r-1}_X \,, \]
which is a finitely successive extension of (locally) free $\cO_X$-modules by \cite{il} I.3.9. Hence the natural map
\[ \varrho : \H^1_x(X,\logwitt X n r) \lra \H^1_x(X,\witt X n r) \]
is injective, and our task is to show the equality
\begin{equation}\label{xy}
(-1)^r \big(\varrho \circ \gys_{i_x} \circ \dKyx\big) (\alpha) = \big(\varrho \circ \delta^{\loc}_{y,x}(\logwitt X n r)\big) (\alpha)
\end{equation}
in $\H^1_x(X,\witt X n r)$. We regard the complex
\[ \witt A n r \os{i_y^*}{\lra} \witt K n r \]
as a representative of $R\Gamma_x(X, \witt X n r)$, where $\witt A r n$ is placed in degree $0$, cf.\ \cite{gros:purity} II.3.3.3. This identification induces an isomorphism
\[ \varphi : \witt K n r / \witt A n r \isom \H^1_x(X,\witt X n r). \]
Now consider a commutative diagram
\[ \xymatrix{ \witt F n {r-1} \ar[rr]^-{\omega \, \mapsto \, \wt{\omega} \cdot \dlog(\ul{\pi})} \ar@{}[rrd]|{\text{(1)}} \ar@{=}[d]
&& \witt K n r / \witt A n r \ar@{}[rrd]|{\text{(2)}} \ar[d]_{\varphi}^{\rwr}
&& \ar@{->>}[ll]_-{(-1) \times \, \text{natural projection}} \witt K n r \ar@{=}[d] \\
\H^0(x,\witt x n {r-1}) \ar[rr]^-{\gys'_{i_x}} \ar@{}[rrd]|{\text{(3)}}
&& \H^1_x(X,\witt X n r) \ar@{}[rrd]|{\text{(4)}}
&& \ar[ll]_-{\delta^{\loc}_{y,x}(\witt X n r)} \H^0(y,\witt y n r) \\
\H^0(x,\logwitt x n {r-1}) \ar[rr]^-{\gys_{i_x}} \ar@{^{(}->}[u]
&& \H^1_x(X,\logwitt X n r) \ar@{^{(}->}[u]^{\varrho}
&& \ar[ll]_-{\delta^{\loc}_{y,x}(\logwitt X n r)} \ar@{^{(}->}[u] \H^0(y,\logwitt y n r)\,,\hspace{-5pt}
}\]
where for $\omega \in \witt F n {r-1}$, $\wt{\omega} \in \witt A n {r-1}$ denotes a lift of $\omega$. The square (1) commutes by the property (P1) mentioned before. The square (2) commutes by a simple (but careful) computation of boundary maps, cf.\ \cite{shm} p.\ 612. By these commutative squares we have
\begin{align*}
 \hbox{RHS of \eqref{xy}} & = (-1)^{r-1} \big(\varrho \circ \delta^{\loc}_{y,x}(\logwitt X n r)\big) \big(\dlog(\ul{f_1}) \cdot \, \dotsb \, \cdot \dlog(\ul{f_{r-1}}) \cdot \dlog(\ul{\pi}) \big) \\
 & \os{\text{(4)}}{=} (-1)^{r-1} \big(\delta^{\loc}_{y,x}(\witt X n r)\big) \big(\dlog(\ul{f_1}) \cdot \, \dotsb \, \cdot \dlog(\ul{f_{r-1}}) \cdot \dlog(\ul{\pi}) \big) \\
 & \hspace{-8.0pt}\os{\text{(2)}+\text{(1)}}{=} (-1)^r \, \gys'_{i_x} \big(\dlog(g_1) \cdot \ \dotsb \ \cdot \dlog(g_{r-1}) \big) \os{\text{(3)}}{=} \hbox{LHS of \eqref{xy}},
\end{align*}
where $g_j \in W_n(F)^{\times}$ denotes the residue class of $\ul{f_j}$ for each $j$. We thus obtain Theorem \ref{thm:shiho} in the case $c=1$.

\smallskip
\subsection{Regular case}\label{sect2-3}
We next treat the case that $c$ is arbitrary but the closure $Z:=\ol {\{ y \}} \subset X$ is regular at $x$. Let $\iota_x: x \hra Z$ and $i : Z \hra X$ be the natural closed immersions. Let us consider the following diagram:
\[\xymatrix{
& \H^0(y,\logwitt y n {r-c+1}) \ar[r]^-{\gys_{i_y}} \ar[ld]_{(-1)^r \, \dKyx \;} \ar[d]^{(-1)^{c-1}\delta}
& \H^{c-1}_y(X,\logwitt X n r) \ar[d]^{\delta^{\loc}_{y,x}(\logwitt X n r)} \\
\H^{0}(x,\logwitt x n {r-c}) \ar[r]_-{\gys_{\iota_x}}
& \H^1_x(Z,\logwitt Z n {r-c+1}) \ar[r]_-{\gys_i}
& \H^c_x(X,\logwitt X n r)\,, \hspace{-5pt}
}\]
where we put $\delta:=\delta^{\loc}_{y,x}(\logwitt Z n {r-c+1})$ for simplicity. The right upper arrow and the right lower arrow are induced by the Gysin morphism for $i$,
so the right square commutes by \eqref{eq-0.4.1}. The left triangle commutes by the previous case. The composite of the bottom row coincides with $\gys_{i_x}$ by (P2). Hence the assertion follows in this case.

\smallskip
\subsection{General case}\label{sect2-4}
We finally consider the general case. The arguments here proceed similarly as for \S\ref{sect1-3}. Let $Z=\ol {\{ y \}} \subset X$ be as in the previous step. We assume that $Z$ is not regular. Take the normalization $f:T \ra Z$. Put $\varSigma:=f^{-1}(x) \subset T$ with the reduced subscheme structure, which is the set of all closed points on $T$ by the assumption \eqref{eq-local-p}. Note that $T$ is regular and that $f$ is finite by the excellence of $Z$. Since a finite morphism is projective (\cite{ega2} 6.1.11), the map $T \ra X$ factors as $T \hra \P^e_X \ra X$ for some integer $e \geq 1$. Let $\iota_y$ be the composite map $y \ra T \ra \P^e_X$.  There is a commutative diagram of schemes
\begin{equation}\label{CD:schemes2}
\xymatrix{
\varSigma \; \ar@<-1pt>@{^{(}->}[r] \ar[d]_h \ar@/^7mm/[rr]^{i_{\varSigma}}
& T \; \ar@<-1pt>@{^{(}->}[r]^-{i_{T}} \ar[d]_f
& \P:=\P^e_X \ar[d]_g \\
x \; \ar@<-1pt>@{^{(}->}[r]^-{\iota_x} \ar@/_7mm/[rr]_{i_x}
& Z \; \ar@<-1pt>@{^{(}->}[r]^-{i_Z}
& X \,, \hspace{-5pt}
}
\end{equation}
where $i_x$, $i_{\varSigma}$ and all horizontal arrows are closed immersions cf.\ \eqref{eq-local-p}.
 Now put $c':=c+e=\codim_{\P}(\varSigma)$, and consider the diagram in Figure 2.
\begin{figure}[htp]
\setlength{\unitlength}{.6mm} {\scriptsize
\begin{picture}
(220,90)(-127,-45) \thinlines
       \put(-99,-35){\vector(2,1){24}}
       \put(67,-35){\vector(-2,1){24}}
       \put(67,36){\vector(-2,-1){24}}
       \put(-99,36){\vector(2,-1){24}}
       \put(-119,34){\line(0,-1){68}}
       \put(-118,34){\line(0,-1){68}}
       \put(84,34){\vector(0,-1){68}}
       \put(-97,40){\vector(1,0){161}}
       \put(-97,-39){\vector(1,0){161}}
       \put(-65,13){\vector(0,-1){26}}
       \put(32,13){\vector(0,-1){26}}
       \put(-42,18){\vector(1,0){51}}
       \put(-42,-18){\vector(1,0){51}}
       \put(-84,-32){{\tiny $\gys_{i_y}$}}
       \put(40,-32){{\tiny $\gys_{i_x}$}}
       \put(36,31){{\tiny $\gys_{i_{\varSigma}}$}}
       \put(-84,31){{\tiny $\gys_{\iota_y}$}}
        \put(86,0){{\tiny $\gys_{h}$}}
        \put(-22,44){{\tiny $(-1)^r\,\dval_{y,\varSigma}$}}
        \put(-22,-45){{\tiny $(-1)^r\,\dKyx$}}
        \put(-63,0){{\tiny $\gys_g$}}
        \put(20,0){{\tiny $\gys_g$}}
        \put(-40,12){{\tiny $(-1)^e\,\delta^{\loc}_{y,\varSigma}(\logwitt {\P} n {r+e})$}}
        \put(-34,-14){{\tiny $\delta^{\loc}_{y,x}(\logwitt X n r)$}}
       \put(-18,28){\tiny {\text(1)}}
       \put(-95,0){\tiny {\text(2)}}
       \put(-18,0){\tiny {\text(3)}}
       \put(58,0){\tiny {\text(4)}}
       \put(-18,-30){\tiny {\text(5)}}
\put(66,39){$\H^0(\varSigma,\logwitt {\varSigma} n {r-c})$}
\put(-140,39){$\H^0(y,\logwitt y n {r-c+1})$}
\put(-140,-41){$\H^0(y,\logwitt y n {r-c+1})$}
\put(66,-41){$\H^0(x,\logwitt x n {r-c})$}
\put(12,17){$\H^{c'}_{\varSigma}(\P,\logwitt {\P} n {r+e})$}
\put(-89,17){$\H^{c'-1}_y(\P,\logwitt {\P} n {r+e})$}
\put(-89,-20){$\H^{c-1}_y(X,\logwitt X n r)$}
\put(12,-20){$\H^{c}_x(X,\logwitt X n r)$}
\end{picture}
\caption{A diagram for the proof of the general case}}
\end{figure}
The square (3) commutes by \eqref{eq-0.4.1}. Moreover, the diagrams (2) and (4) commute by the transitivity property (P2). On the other hand, the diagram (1) is commutative by the result in the previous case. Finally, the outer large square of Figure 2 commutes by Lemma \ref{lem:norm0} in the appendix and the definition of $\dKyx$, cf.\ \S\ref{sect0-6}\,(II). Thus the diagram (5) commutes, i.e., Theorem \ref{thm:shiho}. \hfill \qed
\begin{cor}\label{cor:shiho}
Let $X$ be a smooth variety of pure dimension $d$ over $k$, and let $c$ be a positive integer. Let $i_x:x \ra X$ and $i_y:y \ra X$ be points on $X$ of codimension $c$ and $c-1$, respectively, with $x \in \ol {\{ y \}}$. Then the following diagram commutes{\rm:}
\[ \xymatrix{
\H^1(y,\logwitt y n {d-c+1}) \ar[rr]^-{(-1)^d\,\dKyx} \ar[d]_{\gys_{i_y}}
&& \H^1(x,\logwitt x n {d-c}) \ar[d]^{\gys_{i_x}} \\
\H^c_y(X,\logwitt X n d) \ar[rr]^-{\delta^{\loc}_{y,x}(\logwitt X n d)}
&& \H^{c+1}_x(X,\logwitt X n d)\,.\hspace{-5pt}
}\]
\end{cor}
\begin{pf}
First of all we note that $[\kappa(x):\kappa(x)^p] = p^{d-c}$, because $\kappa(x)$ has transcendence degree $d-c$ over the perfect field $k$. Therefore the upper map is well-defined. For the prove of the corollary we just have to consider the case $c=1$. In fact, the reduction to this case works as in \S\ref{sect2-3} and \S\ref{sect2-4} for Theorem \ref{thm:shiho}; we only have to consider the case $n=c$, and to raise the degrees of all cohomology groups by $1$. Furthermore we have to replace Lemma \ref{lem:norm0} by Lemma \ref{cor:norm1}. \par
In the case $c=1$ we again may replace $X$ by the spectrum of the discrete valuation ring $A=\cO_{X,x}$. By the definition of Kato's residue maps (cf.\ \S\ref{sect0-6}\,(I)), and since
\[ H^2_x(X,\logwitt X n d) \isom H^2_x(X^\h,\logwitt {X^\h} n d) \]
for the henselization $X^\h$ of $X$ at $x$, we may furthermore replace $X$ by $X^\h$. Then $y=\Spec(K)$ for a henselian discrete valuation field $K$ with residue field $\kappa(x)$. Let $y'=\Spec(K^{\sh})$, where $K^{\sh}$ is the strict henselization of $K$. Put $\ol x := \Spec(\ol{\kappa(x)})$ for the separable closure $\ol{\kappa(x)}$ of $\kappa(x)$. Then we get a diagram in Figure 3.\par
\begin{figure}[htp]
\setlength{\unitlength}{.6mm} {\scriptsize
\begin{picture}
(220,90)(-126,-45) \thinlines
       \put(-99,-35){\vector(2,1){24}}
       \put(67,-35){\vector(-2,1){24}}
       \put(67,36){\vector(-2,-1){24}}
       \put(-99,36){\vector(2,-1){24}}
       \put(-117,34){\line(0,-1){68}}
       \put(-116,34){\line(0,-1){68}}
       \put(84,34){\vector(0,-1){68}}
       \put(-97,40){\vector(1,0){161}}
       \put(-97,-39){\vector(1,0){161}}
       \put(-67,13){\line(0,-1){26}}
       \put(-68,13){\line(0,-1){26}}
       \put(38,13){\vector(0,-1){26}}
       \put(-39,18){\vector(1,0){46}}
       \put(-39,-18){\vector(1,0){46}}
       \put(-84,-31){\tiny $a$}
       \put(-94,-28){\tiny $\simeq$}
       \put(49,-31){\tiny $c$}
       \put(59,-28){\tiny $\simeq$}
       \put(49,30){\tiny $b$}
       \put(59,28){\tiny $\simeq$}
       \put(-84,30){\tiny $a$}
       \put(-94,28){\tiny $\simeq$}
       \put(86,0){\tiny $\gys$}
       \put(-26,44){\tiny $(-1)^d\,\dval$}
       \put(-31,-45){\tiny $\delta^\loc(\logwitt X n d)$}
       \put(13,0){\tiny $H^1(x,\gys)$}
       \put(-26,12){\tiny $(-1)^d\,\dval$}
       \put(-40,-14){\tiny $H^1(x,\delta^\loc(\logwitt X n d))$}
       \put(-18,28){\tiny {\text(1)}}
       \put(-93,0){\tiny {\text(2)}}
       \put(-18,0){\tiny {\text(3)}}
       \put(58,0){\tiny {\text(4)}}
       \put(-18,-30){\tiny {\text(5)}}
\put(66,39){$H^1(x,\logwitt x n {d-1})$}
\put(-137,39){$H^1(y,\logwitt y n d)$}
\put(-137,-41){$H^1(y,\logwitt y n d)$}
\put(66,-41){$H^2_x(X,\logwitt X n d)$}
\put(9,17){$H^1(x,H^0(\ol x,\logwitt x n {d-1}))$}
\put(-97,17){$H^1(x,H^0(y',\logwitt y n d))$}
\put(-97,-20){$H^1(x,H^0(y',\logwitt y n d))$}
\put(9,-20){$H^1(x,H^1_{\ol x}(X^\sh,\logwitt X n d))$}
\end{picture}
\caption{A diagram for the proof of Corollary \ref{cor:shiho}}
}
\end{figure}
\noindent
Here the isomorphisms $a, b$ and $c$ come from Hochschild-Serre spectral sequences for the pro-\'etale covering $X^\sh \to X\,(=X^\h)$ given by the strict henselization of $X$. See \S\ref{sect0-6}\,(I.3) for $a$ and $b$, and note the isomorphism
\[  H^{i+1}_{\ol x}(X^\sh,\logwitt X n d) \cong H^i(\ol x , \logwitt x n {d-1}) \; \hbox{($=0$ for $i \ne 0$)} \]
for $c$, cf.\ \cite{moser} Corollary to 2.4. Then the diagram (1) commutes by the definition of Kato's residue maps. The diagram (2) commutes trivially, and the diagrams (4) and (5) commute, because the vertical maps and the two lower horizontal maps are induced by morphisms of sheaves ($\gys$ and $\delta^\loc$), and hence are compatible with the Hochschild-Serre spectral sequence. Finally, it follows from Theorem \ref{thm:shiho} that the square (3) commutes. This implies the commutativity of the outer square and hence the corollary.
\end{pf}
%
%
\smallskip
\subsection{The complex $\bs{\cM^\bullet_{n,X}}$}\label{sect2-5}
Building on work of Moser \cite{moser}, and motivated by Theorem \ref{thm:shiho}, we introduce a complex of \'etale
sheaves and prove a duality result for it (cf.\ \S\ref{sect2-6} below).
\begin{defn}
Let $X$ be a scheme of finite type over $s$. For a point $x$ on $X$, let $i_x$ be the canonical map $x \hra X$.
We define the complex $\cM^\bullet_{n,X}$ of \'etale sheaves on $X$ as
\[ \cM^\bullet_{n,X} :=
 \left( \left\{ \bigoplus{}_{x \in X_{-q}} \ i_{x*} \logwitt x n {-q} \right\} {}_q, \left\{ -\partial^{-q} \right\}_q \right), \]
where $\partial^{-q}$ has the components $\dKyx$ with $y \in X_{-q}$ and $x \in X_{-q-1}$ {\rm(}cf.\ {\rm\S\ref{sect0-7}\rm)}. We often write $\cM_{n,X}$ for the image of $\cM^\bullet_{n,X}$ in $D^b(X_\et,\zpn)$. See Remark \ref{rem:smooth} below for the reason of the sign of the differentials.
\end{defn}
The complex $\cM^\bullet_{n,X}$ coincides with the complex $\wt{\nu}{}_{r,X}$ defined in \cite{moser} up to signs of boundary operators and a shift. If $X$ is smooth over $s$ of pure dimension $d$, then, by a theorem of Gros and Suwa \cite{gs}, the embedding $\logwitt X n d \hra \bigoplus_{x \in X_d} \, i_{x*} \logwitt x n d$ induces a canonical quasi-isomorphism
\stepcounter{equation}
\begin{equation}\label{same}
\logwitt X n d [d] \os{\qsimeq}\lra \cM^\bullet_{n,X}\,.
\end{equation}
Note also the following simple facts: For a closed immersion $i:Z \hra X$ of schemes of finite type over $s$, there is a natural map of complexes
\begin{equation}\label{gysin:DC}
 i_*\cM^\bullet_{n,Z} \lra \cM^\bullet_{n,X}\,.
\end{equation}
If $X$ and $Z$ are smooth of pure dimension $d$ and $d'$, respectively, then this map induces a morphism
\begin{equation}\label{eq2-5-5}
 \gysc_i : i_*\logwitt Z n {d'}[d'] \lra \logwitt X n d [d] \quad \hbox{ in }\; D^b(X_\et,\zpn)
\end{equation}
via \eqref{same} for $X$ and $Z$, which we call the {\it modified Gysin morphism} for $i$.
\addtocounter{thm}{3}
\begin{rem}\label{rem:smooth}
The reason we put the sign $-1$ on the differentials of $\cM^\bullet_{n,X}$ is as follows. Because of these signs, the modified Gysin map \eqref{eq2-5-5} agrees with Gros' Gysin map $\gys_i$ only up to the sign $(-1)^{d-d'}$, cf.\ \cite{sato1} {\rm 2.3.1}. However by this fact, if we define $\zpn(r) := \logwitt X n r [-r]$ for {\rm(}essentially{\rm)} smooth schemes $X$ over $k$ and note property \eqref{eq-0.4.1},
Theorem {\rm\ref{thm:shiho}} for $r=d$ and Corollary {\rm\ref{cor:shiho}} become a commutative diagram
\[\xymatrix{
\H^{m-c+1}(y,\zpn(d-c+1)) \ar[rr]^-{-\dKyx} \ar[d]_{\gysc_{i_y}} && \H^{m-c}(x,\zpn(d-c)) \ar[d]^{\gysc_{i_x}} \\
\H^{m+c-1}_y(X,\zpn(d)) \ar[rr]^-{\delta^{\loc}_{y,x}(\zpn(d))} && \H^{m+c}_x(X,\zpn(d))\,.\hspace{-5pt}
}\]
Note that the groups in the top row are only non-zero for $m=d, d+1$. This shows the perfect analogy with Theorem {\rm\ref{prop:appA}}.
\end{rem}
\noindent
The following lemma shows that the complex $\cM^\bullet_{n,X}$ is suitable for cohomological operations:
\begin{lem}\label{lem:moser}
Let $x$ be a point on $X$ of dimension $q \geq 0$. Then{\rm:}
\begin{enumerate}
\item[(1)]
The sheaf $\logwitt x n q$ on $x_\et$ is $i_{x*}$-acyclic.
\item[(2)]
For a closed immersion $i:Z \hra X$, the sheaf $i_{x*} \logwitt x n q$ on $X_\et$ is $i^!$-acyclic.
\item[(3)]
For an $s$-morphism $f:X \ra Y$, the sheaf $i_{x*} \logwitt x n q$ on $X_\et$ is $f_*$-acyclic.
\end{enumerate}
\end{lem}
\begin{pf}
For (1) and (2), see \cite{moser} 2.3 and 2.4. We prove (3). For a point $y \in Y$, we have
\begin{align*}
\big(R^mf_*(i_{x*} \logwitt x n q)\big)_{\ol y} & \cong \H^m(X \times_Y \Spec(\cO_{Y,\ol y}^{\sh}),i_{x*} \logwitt x n q) \\
 & \os{\text{(1)}}{\cong} \H^m(x \times_Y \Spec(\O_{Y,\ol y}^\sh), \logwitt x n q)
\end{align*}
and the last group is zero for $m>0$ by the same argument as in loc.\ cit.\ 2.5.
\end{pf}
\begin{cor}[cf.\ \cite{moser} Corollary to Theorem 2.4]
\label{cor:moser} For a closed immersion $i : Z \hra X$, the map \eqref{gysin:DC} induces an isomorphism
\[ \gysc_i : \cM_{n,Z} \isom Ri^!\cM_{n,X} \quad \hbox{ in } \; D^+(Z_\et,\zpn). \]
\end{cor}

\smallskip
\subsection{Relative duality theory}\label{sect2-6}
Let $\cV_s$ be the category of schemes separated of finite type over $s$ and separated $s$-morphisms of finite type.
\begin{thm}\label{thm:moser2}
Suppose that there exists an assignment of morphisms
\[ \Tr: ( f : Y \ra X \hbox{ in } \cV_s ) \lmt (\Tr_f : Rf_!\cM_{n,Y} \ra \cM_{n,X} \hbox{ in } D^+(X_\et,\zpn)) \]
which satisfy the following three conditions {\rm(i)\,--\,(iii):}
\begin{enumerate}
\item[(i)]
If $f$ is \'etale, then $\Tr_f$ agrees with the composite morphism
\[ Rf_!\cM_{n,Y} = Rf_!f^*\cM_{n,X} \os{f_!}\lra \cM_{n,X}\,, \]
where the arrow $f_!$ denotes the adjunction morphism $Rf_!f^*=Rf_!Rf^! \to \id$ {\rm(}cf.\ \cite{sga4} {\rm XVIII.3}, the equality $\cM_{n,Y} = f^*\cM_{n,X}$ is straight-forward{\rm)}.
\item[(ii)]
If $f$ is a closed immersion, then $\Tr_f$ agrees with the composite morphism
\[\xymatrix{
Rf_*\cM_{n,Y} \ar[rr]^-{Rf_*(\gysc_f)} && Rf_*Rf^!\cM_{n,X} \ar[r]^-{f_*} & \cM_{n,X}\,,
}\]
where the arrow $f_*$ denotes the adjunction morphism $Rf_*Rf^!=Rf_!Rf^! \ra \id$.
\item[(iii)]
For morphisms $g:Z \ra Y$ and $f:Y \ra X$ with $h:=f \circ g$, $\Tr_h$ agrees with the composition
\[\xymatrix{
Rh_!\cM_{n,Z} = Rf_!Rg_!\cM_{n,Z} \ar[rr]^-{Rf_!(\Tr_g)} && Rf_!\cM_{n,Y} \ar[r]^-{\Tr_f} & \cM_{n,X}\,.
}\]
\end{enumerate}
Then for a map $f:Y \ra X$ in $\cV_s$, the adjoint morphism $\Tr^f$ of $\Tr_f$ is an isomorphism{\rm:}
\[ \Tr^f: \cM_{n,Y} \isom R f^!\cM_{n,X} \quad \hbox{ in } \; D^+(Y_\et,\zpn). \]
\end{thm}
\noindent
This theorem is a variant of Moser's duality \cite{moser} 5.6 (which itself generalizes Milne's duality for smooth projective varieties \cite{milne:mot}). However, because Theorem \ref{thm:moser2} looks quite different from Moser's formulation, we outline a proof of our statement below in \S\ref{sect2-7}. The main result of this section is the following theorem:
\begin{thm}\label{thm:jss}
There exists a unique assignment of morphisms
\[ \tr: ( f:Y \ra X \hbox{ in } \cV_s ) \lmt ( \tr_f : Rf_!\cM_{n,Y} \ra \cM_{n,X} \hbox{ in } D^+(X_\et,\zpn) ) \]
that satisfies the conditions {\rm(i)\,--\,(iii)} in Theorem {\rm\ref{thm:moser2}} with $\Tr:=\tr$. Consequently, for a map $f:Y \ra X$ in $\cV_s$, the morphism $\tr^f:\cM_{n,Y} \ra Rf^!\cM_{n,X}$ adjoint to $\tr_f$ is an isomorphism.
\end{thm}
\noindent
We will prove Theorem \ref{thm:jss} in \S\S\ref{sect2-8}--\ref{sect2-9} below.

\smallskip
\subsection{Proof of Theorem \ref{thm:moser2}}\label{sect2-7}
By the transitivity property (iii) of $\Tr$, the assertion is reduced to the case of a structure morphism $f:X \ra s$, and moreover, by the property (i) of $\Tr$, we may suppose that $s= \ol s$ (i.e., $k$ is algebraically closed) and that $f$ is proper. In this situation, we claim the following:
\begin{thm}\label{thm:moser}
Let $X$ be a proper scheme of finite type over the algebraically closed field $k$ of characteristic $p>0$, with structural morphism $f: X \to \Spec(k)$. Then, for any constructible $\zpn$-sheaf $\cF$ on $X_\et$ and any integer $m$, the pairing
\[ \alpha_X(m,\cF):\H^m(X,\cF) \times \Ext^{-m}_{X,\zpn}(\cF,\cM_{n,X}) \lra \H^0(X,\cM_{n,X}) \os {\Tr_f}\lra \zpn \]
induced by Yoneda pairing is a non-degenerate pairing of finite groups.
\end{thm}
\def\jzp{j_!\hspace{1pt}\zpn}
We first prove Theorem \ref{thm:moser2}, admitting Theorem \ref{thm:moser}: Applying \ref{thm:moser} to $\cF=\jzp$ with $j:U \ra X$ \'etale, and noting the isomorphisms
\[ \Ext^{-m}_{X,\zpn}(\jzp,\cM_{n,X}) \cong \Ext^{-m}_{U,\zpn}(\zpn,\cM_{n,U})\cong H^{-m}(U,\cM_{n,U})\,, \]
we obtain isomorphisms
\[ \H^{-m}(U,\cM_{n,U}) \os{a}\cong \Hom_{\,\zpn\text{-mod}}(\H^{m}(X,\jzp),\zpn) \os{b}\cong \H^{-m}(U,Rg^!\zpn) \]
for any $m \in \bZ$, where $g= f \circ j$, the first isomorphism comes from the pairing, and the second isomorphism comes from the adjunction between $Rg^!$ and $Rg_!$ and the fact that $\zpn$ is an injective object in the category of $\zpn$-modules. We verify that this composite map agrees with that coming from $\Tr^f$ --- then the morphism $\Tr^f$ is bijective on cohomology sheaves, and we obtain Theorem \ref{thm:moser2}. Indeed, by the definition of the pairing, the map $a$ sends $x \in H^{-m}(U,\cM_{n,U})=\Hom_{D(X,\zpn)}(\jzp,\cM_{n,X}[-m])$ to the composition \[ \H^m(X,\jzp) \os{x}\lra H^0(X,\cM_{n,X})=H^0(s,Rf_*\cM_{n,X}) \os{\Tr_f}\lra \zpn, \] which, by \ref{thm:moser2}\,(i) and (iii), coincides with the map induced by
\[ Rg_!\hspace{1pt}\zpn[m] \os{x}\lra Rg_!\cM_{n,U} \os{\Tr_g}\lra \zpn\,. \]
By definition (and functoriality of adjunction), the map $b$ sends this to the composition
\[ \zpn[m] \os{x}\lra \cM_{n,U} \os{\Tr^{\hspace{1pt}g}}\lra Rg^!\zpn\,, \]
which shows the claim, again by \ref{thm:moser2}\,(i) and (iii). \qed
\par\medskip\smallskip
As for Theorem \ref{thm:moser}, it follows from the arguments in \cite{moser} 5.6. More precisely, it follows from the properties (i)\,--\,(iii) of $\Tr$, the steps (a)\,--\,(c), (f)\,--\,(k) of loc.\ cit.\ 5.6, and the following lemma:
\begin{lem}\label{lem:moser2}
Let $f : X \ra  s (=\ol s)$ be a proper smooth morphism with $X$ connected. Then for an integer $m$ and a positive integer $t \leq r$, the pairing $\alpha_X(m,\Z/p^t)$ {\rm(}cf.\ Theorem {\rm\ref{thm:moser})} is a non-degenerate pairing of finite groups.
\end{lem}
\begin{pf*}{\it Proof of Lemma \ref{lem:moser2}}
The problem is reduced to the case $t=r$ by \eqref{same} and \cite{moser} 5.4. Now we note that Milne duality \cite{milne:mot} 1.11 gives an isomorphism of finite groups in our case. Indeed, with the notation of \cite{milne:mot} p.\ 305, the unipotent part of the group scheme $\ul{H}^m(X,\zpn)$ is trivial, which follows from the short exact sequence
\[ 0 \lra \ul{H}^m(X,\zpn) \lra \ul{H}^m(X,W_n\cO_X) \os{1-F}\lra \ul{H}^m(X,W_n\cO_X) \lra 0 \]
and the fact that $1-F$ is \'etale. Therefore it is enough to show that the composite map
\[ \Tr'_f : \H^d(X,\logwitt X n d) \os{\eqref{same}}\cong \H^0(X,\cM_{n,X}) \os{\Tr_f}\lra \zpn \]
with $d:=\dim(X)$ coincides with the trace map $\eta_n$ in \cite{milne:mot} p.\ 308, up to a sign. But, for a closed point $i_x:x \hra X$, $\Tr'_f$ sends the cycle class $\gysc_{i_x}(1) \in \H^d(X,\logwitt X n d)$ (= the image of $1$ under $\gysc_{i_x}$) to $1$ by the properties (ii), (iii) of $\Tr$, and hence $\Tr'_f=(-1)^d \cdot \eta_n$ by Remark \ref{rem:smooth}\,(1). This completes the proof of Lemma \ref{lem:moser2}, Theorem \ref{thm:moser} and Theorem \ref{thm:moser2}.
\end{pf*}
\begin{rem}\label{rem:moser}
Note that step {\rm (j)} of \cite{moser} {\rm5.6} uses de Jong's alteration theorem \cite{dj} {\rm4.1}.
\end{rem}
\begin{cor}\label{cor:jss}
Suppose that there exist two assignments $\sigma:f \mapsto \sigma_f$ and $\tau : f \mapsto \tau_f$ satisfying {\rm(i)\,--\,(iii)} in Theorem {\rm\ref{thm:moser2}} with $\Tr:=\sigma$ and $\tau$, respectively. Then we have $\sigma=\tau$.
\end{cor}
\vspace{-5pt}
\begin{pf}
Let $f:Y \ra X$ be a morphism in $\cV_s$. We show that $\sigma_f=\tau_f$ as morphisms $Rf_!\cM_{n,Y} \ra \cM_{n,X}$ in $D^+(X_\et,\zpn)$, in two steps. We first prove the case $X=s$ (hence $\cM_{n,X}=\zpn$). By the properties (i) and (iii), we may suppose that $f$ is proper. Then $Rf_!\cM_{n,Y}=Rf_*\cM_{n,Y}$ is computed by the complex $f_*\cM_{n,Y}^\bullet$ by Lemma \ref{lem:moser}\,(3), and the morphisms $\sigma_f$ and $\tau_f$ are determined by maps $f_*\cM_{n,Y}^0 \ra \zpn$ of sheaves on $s_\et$ by \S\ref{sect0-5-4}\,(1). Hence in view of the properties (ii) and (iii) and the assumption that $s$ is perfect, the problem is reduced to the case where $f$ is \'etale, and we obtain $\sigma_f=\tau_f$ by the property (i). This completes the first step. Next we prove the general case. Let $\sigma^f$ and $\tau^f$ be the adjoint morphisms of $\sigma_f$ and $\tau_f$, respectively. By adjunction, it is enough to show $\sigma^f=\tau^f$ as morphisms $\cM_{n,Y} \ra Rf^!\cM_{n,X}$ in $D^+(Y_\et,\zpn)$. Let $g:X \ra s$ be the structure map and put $h:=g \circ f$. By the first step and the property (iii), we have
\[ Rf^!(\sigma^g) \circ \sigma^f = \sigma^h = \tau^h = Rf^!(\tau^g) \circ \tau^f \]
as morphisms $\cM_{n,Y} \ra Rh^!\zpn$. On the other hand, we have $Rf^!(\sigma^g)=Rf^!(\tau^g)$ by the first step, and these are isomorphisms in $D^+(Y_\et,\zpn)$ by Theorem \ref{thm:moser2}. Hence we have $\sigma^f = \tau^f$. This completes the proof of Corollary \ref{cor:jss}.
\end{pf}

\smallskip
\subsection{Covariant functoriality}\label{sect2-8}
In this subsection, we prove Lemma \ref{lem:jss2} stated below (cf.\ \cite{moser} 4.1), which is a key ingredient of Theorem \ref{thm:jss}. Let $f:Y \ra X$ be a morphism in $\cV_s$. Let $q$ be a non-negative integer and let $x$ (resp.\ $y$) be a point on $X$ (resp.\ on $Y$) of dimension $q$. Let $f_y$ (resp.\ $i_x$) be the composite map $y \ra Y \ra X$ (resp.\ $x \ra X$). We define a map of sheaves on $X_\et$
\[ \tr_{f,(y,x)}:f_{y*} \logwitt y n q \lra i_{x*}\logwitt x n q \]
as Gros' Gysin map for $y \to x$ (\cite{gros:purity} II.1.2.7), if $y$ is finite over $x$ via $f$. We define $\tr_{f,(y,x)}$ as zero otherwise. Collecting this map for points on $Y$ and $X$, we obtain a map of {\it graded abelian sheaves} on $X_\et$
\[ \tr_f^\bullet : f_*\cM^\bullet_{n,Y} \lra \cM^\bullet_{n,X}\,. \]
By definition and \cite{gros:purity} II.2.1.1, this map of graded sheaves satisfies transitivity, that is, for morphisms $g:Z \ra Y$ and $f:Y \ra X$ in $\cV_s$, we have the equality
\begin{equation}\label{lem:jss2:transitive}
\tr_f^\bullet \circ f_*(\tr_g^\bullet) = \tr_{f \circ g}^\bullet
\end{equation}
of maps of graded sheaves on $Y_\et$. We prove the following lemma:
\stepcounter{thm}
\begin{lem}\label{lem:jss2}
Suppose that $f$ is proper. Then $\tr_f^\bullet$ is a map of complexes. Consequently, $\tr_f^\bullet$ induces a morphism
\[ \tr_f: Rf_*\cM_{n,Y} \lra \cM_{n,X} \quad \hbox{ in } \; D^b(X_\et,\zpn) \]
by Lemma {\rm\ref{lem:moser}\,(3)}.
\end{lem}
\begin{pf}
We have to show the commutativity of the following diagram for each $q \le -1$:
\begin{equation}\notag
\xymatrix{
f_*\cM_{n,Y}^{-q} \ar[r]^{\partial^{-q}} \ar[d]_{\tr_f^{-q}} & f_*\cM_{n,Y}^{-q-1} \ar[d]^-{\tr_f^{-q-1}} \\
\cM_{n,X}^{-q} \ar[r]^{\partial^{-q}} & \cM_{n,X}^{-q-1} \,,
}
\end{equation}
Fix a negative integer $q$, and a point $w \in Y_{-q}$. Put $z:=f(w)$. There happen the following three cases:
\begin{enumerate}
\item[(i)] $z \in X_{-q}$ \qquad\qquad (ii) $z \in X_{-q-1}$ \qquad\qquad (iii) otherwise.
\end{enumerate}
In the third case, we see that all points $y$ in the closure $\ol{\{w\}}$ map to points on $X$ of dimension $\le -q-2$. Hence the maps $\partial^{-q} \circ \tr_f^{-q-1}$ and $\tr_f^{-q} \circ \partial^{-q}$ are both zero on the direct summand $i_{w*}\logwitt w {-q} n$ of $f_*\cM_{n,Y}^{-q}$. We next treat the case (ii). In this case, we are ought to show that the following sequence is a complex:
\[f_{w*}\logwitt w n {-q} \lra \bigoplus_{y \in Y_{-q-1} \, \cap \, \ol {\{ w \} } \, \cap \, f^{-1}(z)} \ f_{y*}\logwitt y n {-q-1} \lra i_{z*}\logwitt z n {-q-1}. \]
We put $C:=f^{-1}(z) \cap \ol {\{ w \} }$ endowed with reduced subscheme structure. Since $C$ is proper over $z$ by the properness of $f$, $C$ is a proper curve over $z$ with generic point $w$ (\cite{hartshorne} III.9.6).
By Lemma \ref{lem:norm0}, the assertion is reduced to the case that $C$ is a projective line over $z$. The assertion then follows from {\it Claim} in the proof of Lemma \ref{lem:norm0}.

We finally prove the case (i). Fix an arbitrary point $x \in X_{-q-1} \cap \ol {\{z\}}$. Our task is to show the commutativity of the diagram
\stepcounter{equation}
\begin{equation}\label{lem:jss2:commute2}
\xymatrix{
\displaystyle  f_{w*}\logwitt w n {-q}
\ar[rr]^-{\bigoplus_{y} \dKwy } \ar[d]_-{\tr_{f,(w,z)}}
&&
\displaystyle \bigoplus_{y \in Y_{-q-1} \cap f^{-1}(\ol {\{ x \} })}^{\phantom{|}} \ f_{y*}\logwitt y n {-q-1} \ar[d]^-{\sum_y \, \tr_{f,(y,x)}} \\ i_{z*}\logwitt z n {-q} \ar[rr]^-{\dKzx}
&& i_{x*}\logwitt x n {-q-1}\,.\hspace{-5pt}
}
\end{equation}
Let $T$ be the localization of $\ol{ \{ z \} }$ at $x$, and put
\[ Z_w := \ol { \{ w \} } \times_X T \;(=\ol { \{ w \} } \cap f^{-1}(T))\,. \]
If $y \in Y_{-q-1} \cap f^{-1}(\ol {\{z\} })$ is away from $Z_w$, then $y$ is outside of $\ol {\{ w \} } \cap f^{-1}(x)$ and hence at least one of $\dKwy$ and $\tr_{f,(y,x)}$ is zero. Thus the commutativity of \eqref{lem:jss2:commute2} is reduced to that of the following diagram:
\vspace{-8pt}
\begin{equation}\label{lem:jss2:commute3}
\xymatrix{
f_{w*}\logwitt w n {-q} \ar[rr]^-{\bigoplus_y\,\dKwy} \ar[d]_-{\tr_{f,(w,z)}}
&& \displaystyle \bigoplus_{y \in Y_{-q-1} \cap Z_w}^{\phantom{|}}\ f_{y*} \logwitt y n {-q-1} \ar[d]^-{\sum_y \, \tr_{f,(y,x)}} \\
i_{z*}\logwitt z n {-q} \ar[rr]^-{\dKzx} && i_{x*}\logwitt x n {-q-1}\,.\hspace{-5pt}
}
\end{equation}
We claim here the following:
\par
\vspace{8pt}
\noindent
{\it Claim. The canonical morphism $f_T:Z_w \ra T$ is finite, and $Y_{-q-1} \cap Z_w$ agrees with the set of all closed points on $Z_w$.}
\par
\vspace{8pt}
\noindent
{\it Proof of Claim.}
By the properness of $f_T$ and \cite{ega3} 4.4.2, it suffices to show that $f_T$ is quasi-finite. Note that $\kappa(w)$ is a finite field extension of $\kappa(z)$. Let $\nu:U \ra T$ be the normalization of $T$ in $\kappa(w)$. Then $\nu$ is finite (cf.\ \cite{hartshorne} I.3.9A) and $U$ is the spectrum of a Dedekind ring, which imply that $\nu$ factors as $U \ra Z_w \ra T$ by the valuative criterion for proper morphisms (cf.\ loc.\ cit.\ II.4.7). Here the map $U \ra Z_w$ is surjective, because it is proper and dominant. Hence $f_T$ is quasi-finite by the finiteness of $\nu$ and we obtain the assertion. The second assertion immediately follows from the finiteness of $f_T$. \qed
\par
\vspace{8pt}
We turn to the proof of Lemma \ref{lem:jss2} and prove the commutativity of \eqref{lem:jss2:commute3}. Since the problem is \'etale local at $x \in T$, we assume that $T$ and $Z_w$ are strictly henselian by replacing them with $\Spec(\cO_{T,\ol x}^\sh)$ and a connected component of $Z_w \times_T \Spec(\cO_{T,\ol x}^\sh)$, respectively. Then by the Bloch-Gabber-Kato theorem (\cite{Bloch-Kato} 2.1) and Lemma \ref{lem:norm0} in the appendix, we are reduced to the commutativity of residue maps of Milnor $K$-groups (\S\ref{sect0-6}\,(I.2)) via norm maps due to Kato \cite{kk:res}, Lemma 3
(which assumes the domains concerned are normal, but is easily generalized to our situation by a standard argument using normalization):
\vspace{-6pt}
\[\xymatrix{
K^M_{-q}(\kappa(w))/p^n \ar[rr]^-{\bigoplus_y\,\dMwy} \ar[d]_-{N_{\kappa(w)/\kappa(z)}}
&&\displaystyle\bigoplus_{y \in (Z_w)_0}^{\phantom{|}}\ K^M_{-q-1}(\kappa(y))/p^n \ar[d]^-{\sum_y \, N_{\kappa(y)/\kappa(x)}}\\
K^M_{-q}(\kappa(z))/p^n \ar[rr]^-{\dMzx} && K^M_{-q-1}(\kappa(x))/p^n\,.\hspace{-5pt}
}\]
Hence \eqref{lem:jss2:commute3} commutes in this case. This completes the proof of Lemma \ref{lem:jss2}.
\end{pf}
\addtocounter{thm}{2}
\begin{rem}\label{rem:jss}
If $f$ is finite and \'etale, then $\tr_f^\bullet$ coincides with the adjunction map
\[ f_* : f_*\cM^\bullet_{n,Y} = f_*f^*\cM^\bullet_{n,X} \lra \cM^\bullet_{n,X}\,. \]
Indeed, the claim is reduced to the case of a finite separable extension of a point, which follows from a standard base-change argument and \cite{milne} {\rm V.1.12}.
\end{rem}

\smallskip
\subsection{Proof of Theorem \ref{thm:jss}}\label{sect2-9}
By Corollary \ref{cor:jss}, it remains to show the existence of a desired assignment. For a map $f:Y \ra X$ in $\cV_s$, we define the morphism $\tr_f$ as follows. If $f$ is proper, then we define $\tr_f$ as that constructed in Lemma \ref{lem:jss2}. Next suppose that $f$ is not proper. Take a compactification of $f$, i.e., an open immersion $ j : Y \hra Z$ and a proper map $ g : Z \ra X$ with $f=g \circ j$, and define $\tr_{f,(Z,j,g)}$ as the composite morphism
\begin{equation}\label{comp2}
\tr_{f,(Z,j,g)}: Rf_! \cM_{n,Y} = Rg_*Rj_!j^*\cM_{n,Z} \os{j_!}\lra Rg_*\cM_{n,Z} \os{\tr_g}\lra \cM_{n,X}\,,
\end{equation}
where the arrow $j_!$ is defined by the adjunction morphism $Rj_!j^*\cM_{n,Z} \to \cM_{n,Z}$. We are going to define
\[ \tr_f:=\tr_{f,(Z,j,g)}\,. \]
To verify the well-definedness, it suffices to show the following:
\stepcounter{thm}
\begin{lem}\label{lem:jss}
Let $Y \os{\varrho}{\hra} V \os{h}{\ra} X$ be another compactification of $f$. Then we have
\[ \tr_{f,(Z,j,g)}=\tr_{f,(V,\varrho,h)}\,. \]
\end{lem}
\vspace{-8pt}
\begin{pf}
Replacing $Z$ by the closure of the image of $Y  \os{(j,\varrho)}{\lra} Z \times_X V$, we may suppose that there exists a proper morphism $\pi:Z \ra V$ fitting into a commutative diagram with cartesian square
\[\xymatrix{ Y \, \ar@{^{(}->}[r]^j \ar@<-1.5pt>@{=}[d] \ar@{}[rd]|{\square} & Z \ar[rd]^g \ar[d]_\pi \\
 Y \, \ar@{^{(}->}[r]^\varrho  & V \ar[r]^h & X\,.\hspace{-5pt} }\]
Since $\tr_g= \tr_h \circ Rh_*(\tr_\pi)$ by \eqref{lem:jss2:transitive} and Lemma \ref{lem:moser}\,(3), it remains to show that $\tr_{\varrho,(Z,j,\pi)}$ agrees with the adjunction morphism
\[ \varrho_! \; : \; \varrho_!\cM_{n,Y} = \varrho_!\varrho^*\cM_{n,V} \lra \cM_{n,V}\,. \]
Indeed, since $\pi^{-1}(Y)=Y$, $\varrho^*(\tr_{\varrho,(Z,j,h)})$ is the identity morphism of $\cM_{n,Y}$ (cf.\ \eqref{comp2}), which implies that $\tr_{\varrho,(Z,j,h)} = \varrho_!$. This completes the proof of the lemma.
\end{pf}
\medskip
Thus we obtained a well-defined assignment $f \mapsto \tr_f$. We show that this satisfies the conditions (i)--(iii) in Theorem \ref{thm:moser2}. The condition (ii) holds obviously by definition (cf.\ \eqref{gysin:DC}, \S\ref{sect2-8}). We next show the condition (i). Suppose that $f:Y \ra X$ is \'etale. Take an open immersion $j:Y \hra Z$ and a finite map $g:Z \ra X$ with $f=g \circ j$ (cf.\ \cite{milne} 1.8). We claim that the morphism $\tr_{f,(Z,j,g)}$ coincides with the adjunction $f_!$, which implies (i).
Indeed, since $f$ is \'etale and $g$ is finite, $\tr_{f,(Z,j,g)}$ is represented by the composite map of complexes
\begin{equation}\notag
\xymatrix{
f_{!}\cM^\bullet_{n,Y} = g_* j_!j^*\cM^\bullet_{n,Z} \ar[r]^-{g_*(j_!)} & g_*\cM^\bullet_{n,Z} \ar[r]^-{\tr_g^\bullet} & \cM^\bullet_{n,X} \,,
}
\end{equation}
which agrees with $f_!$ by a similar argument as for Remark \ref{rem:jss} (see also \cite{milne} II.3.18).
We finally show the condition (iii), that is, for two maps $g:Z \ra Y$ and $f:Y \ra X$ in $\cV_s$, we prove
\stepcounter{equation}
\begin{equation}\label{eq2-9-3} \tr_h=\tr_f \circ Rf_!(\tr_g) \quad \hbox{ with \, $h:=f \circ g$.} \end{equation}
If $f$ and $g$ are open immersions, \eqref{eq2-9-3} follows from the property (i) and the transitivity of adjunction maps for open immersions. If $f$ and $g$ are proper, \eqref{eq2-9-3} follows from \eqref{lem:jss2:transitive} and Lemma \ref{lem:moser}\,(3). Hence, if $g$ is an open immersion {\it or} $f$ is proper, then we obtain \eqref{eq2-9-3} by the previous two cases. We show the general case. Take compactifications of $f$ and $g$ as follows:
\[ \xymatrix{ & T \ar[rd]^q & & V \ar[rd]^{\pi} & \\
 Z \ar[rr]_g \ar@{^{(}->}[ru]^{j} && Y \ar[rr]_f \ar@{^{(}->}[ru]^{\varrho} && X\,,\hspace{-5pt}} \]
where $j$ and $\varrho$ are open immersions and $q$ and $\pi$ are proper maps which make the triangles commutative.
Because we already know, by the previous cases, that
\[ \tr_h= \tr_{\pi} \circ R\pi_*(\tr_{\alpha}) \circ R(f \circ q)_!(\tr_j) \quad \hbox{ with \, $\alpha:=\varrho \circ q$,} \]
 it remains to show the following composite morphism agrees with $\tr_{\alpha}$:
\[\xymatrix{
R\alpha_!\cM_{n,T} \ar[rr]^-{R \varrho_!(\tr_q)} && R\varrho_! \cM_{n,Y} \ar[r]^-{\varrho_!} & \cM_{n,V}\,.
}\]
We prove this equality. Take an open immersion $\beta:T \hra W$ with {\it dense image} and a proper map $\gamma:W \ra V$ satisfying $\alpha= \gamma \circ \beta$. Then one can easily check that the square
\[\xymatrix{
T \, \ar@{^{(}->}[r]^\beta \ar@<-1.5pt>[d]_q & W \ar[d]^\gamma \\ Y \, \ar@{^{(}->}[r]^\varrho & V
}\]
is cartesian. Hence we have
\begin{align*}
 \varrho^*(R\gamma_*(\beta_!)) &=\id_{Rq_*}  \quad & \hbox{($\beta_!$ denotes $R\beta_!\beta^* \to \id$)} \\
 \varrho^*(\tr_{\gamma})& =\tr_q  \qquad & \hbox{(cf.\ \S\ref{sect2-8})}
\end{align*}
and thus $\varrho^*(\tr_{\alpha,(W,\beta,\gamma)})=\tr_q$, which implies $\varrho_! \circ R \varrho_!(\tr_q)=\tr_{\alpha,(W,\beta,\gamma)} = \tr_{\alpha}$. This completes the proof of Theorem \ref{thm:jss}. \qed
\smallskip
\subsection{Purity for logarithmic Hodge-Witt sheaves}\label{sect2-10}
Theorem \ref{thm:jss} implies the following purity result, whose special case was needed in \cite{js} 2.14.
\begin{cor}\label{cor:p-purity}
Let $f : X \to Y$ be a morphism of smooth varieties of pure dimension $d$ and $e$, respectively over $s=\Spec(k)$.  Then there is a canonical Gysin isomorphism
\[ \tr^f : \logwitt X n d [d] \isom Rf^! \logwitt Y n e [e] \,. \]
Especially, for $g : X \to s$ smooth of dimension $d$ we get a canonical isomorphism
\[ \tr^g : \logwitt X n d [d] \isom Rg^!\zpn \,. \]
\end{cor}
\begin{pf}
The first claim follows from the isomorphisms
\[ \logwitt X n d [d] \; \us{\;\eqref{same}}\isom \; \cM_{n,X} \; \us{\;\text{\ref{thm:jss}}}\isom \;
 Rf^!\cM_{n,Y} \; \us{\;\eqref{same}}\lisom \; Rf^!\logwitt Y n e [e] \,. \]
For the special case note that $\logwitt s n 0 = \zpn$.
\end{pf}
\begin{rem}\label{rem:p-purity}
With the notation $\zpn(r)_X := \logwitt X n r [-r]$ the purity isomorphism in Corollary {\rm\ref{cor:p-purity}} becomes
\addtocounter{equation}{2}
\begin{equation}\label{eq2-10-3}
\tr^f : \zpn(d)_X[2d] \isom Rf^!\zpn(e)_Y[2e] \,.
\end{equation}
When $f$ is a closed immersion, $\tr^f$ is adjoint to the modified Gysin morphism \eqref{eq2-5-5}. When $f$ is proper, $\tr^f$ is adjoint to Gros' Gysin morphism $\gys_f$ only up to the sign $(-1)^{d-e}$, cf.\ Remark {\rm\ref{rem:smooth}}.
\end{rem}

\smallskip
\subsection{Bloch-Ogus complexes and Kato complexes}\label{sect2-11}
Finally we have the following application to Kato complexes, which is analogous to Theorem \ref{thm:Katocomplexprimetop}. Let $S$ be a scheme which is smooth of finite type over $k$ and of pure dimension $d$. (Most interesting is the case $S=\Spec(k)$, $d=0$, which was needed in \cite{js} 2.14.) For a separated scheme of finite type over $S$, $f: X \to S$, define its homology with coefficients in $\zpn(-d)$ by
\begin{equation}\label{eq2-11-1}
H_a(X/S,\zpn(-d)) := H^{-a}(X, Rf^!\zpn(d)_S)\,.
\end{equation}
These groups define a homology theory on the category of separated $S$-schemes of finite type, in the sense of \cite{js} 2.1 (cf.\ loc.\ cit.\ 2.2), and in a standard way one obtains a niveau spectral sequence
\begin{equation}\label{eq2-11-2}
E^1_{q,t}(X/S,\zpn(-d)) = \bigoplus_{x\in X_q} H_{q+t}(x/S,\zpn(-d)) \Ra H_{q+t}(X/S,\zpn(-d)) \hspace{-10pt} \end{equation}
for $X$ as above (cf.\ \S\ref{sect1-5} and \cite{js} 2.7). Note that $X_q$ agrees with $(X/S)_q$ in the sense of \eqref{eq1-5-4}, because $S$ is of finite type over $k$.
\addtocounter{thm}{2}
\begin{thm}\label{thm:Katocomplexp-case}
Let $X$ be a separated $S$-scheme of finite type.
\begin{enumerate}
\item[(1)] There is a canonical isomorphism
\[ E^1_{q,t}(X/S,\zpn(-d)) \cong \bigoplus_{x\in X_q} \ H^{q-t-2d}(x,\zpn(q)) =
 \bigoplus_{x\in X_q} \ H^{-t-2d}(x, \logwitt x n q)\,. \]
\item[(2)] Via these isomorphisms, the Bloch-Ogus complex $E^1_{*,t}(X/S,\zpn)$ coincides with the sign-modified modified Kato complex $C_{p^n}^{-t-2d,0}(X)^{(-)}$.
\item[(3)] Especially, for a separated $k$-scheme $X$ of finite type, purity induces an isomorphism $E^1_{*,t}(X/k,\zpn)\cong C_{p^n}^{-t,0}(X)^{(-)}$.
\end{enumerate}
\end{thm}
\begin{pf}
(1) follows from the purity isomorphism
\stepcounter{equation}
\begin{align}
\notag
H_a(V/S,\zpn(-d)) & = H^{-a}(V,Rf^!\zpn(d)_S) \\
\label{eq2-11-3}
& \hspace{-8.2pt}\us{\text{\eqref{eq2-10-3}}}\cong H^{-a+2q-2d}(V,\zpn(q)_V)
\end{align}
for $f: V \to S$ with $V$ smooth of pure dimension $q$.
\par
Since (3) is a special case of (2), we prove (2) in what follows, by similar arguments as in the proof of Theorem \ref{thm:Katocomplexprimetop}. The question is local in $S$ and $X$. Therefore we may assume that $f : X \to S$ factors as follows:
\[\xymatrix{ X \; \ar@{^{(}->}[r]^i & P \ar[r]^\pi & S\,, }\]
where $\pi$ is a smooth morphism of pure relative dimension $N$ and $i$ is a closed immersion. The Gysin isomorphism $\zpn(d+N)[2N] \simeq R\pi^!\zpn(d)$ from \eqref{eq2-10-3} induces an isomorphism of homology theories
\[ \gamma: H_{*-2N}(-/P,\zpn(-d-N)) \isom H_*(-/S,\zpn(-d)) \]
on all subschemes of $P$, and therefore an isomorphism between the corresponding spectral sequences. Moreover, for an open subscheme $j_U: U \hra P$ and a closed subscheme $i_V: V \hra U$ of dimension $q$, the purity isomorphism \eqref{eq2-10-3} for the composition
\[\xymatrix{ g = \pi \circ j_U \circ i_V : V \; \ar@{^{(}->}[r]^-{i_V} & U \; \ar@{^{(}->}[r]^-{j_U} & P \ar[r]^-{\pi} & S }\]
factors as
\[ \begin{CD}
\zpn(q)_V[2q] @>{\quad\tr^{i_V}\quad}>> Ri_V^!\zpn(d+N)_U[2(d+N)] \\
 @>{Ri_V^!j_U^*(\tr^\pi)}>> Ri_V^!j_U^* R\pi^!\,\zpn(d)_S[2d] \\
 @= Rg^!\zpn(d)_S[2d] \,.
\end{CD}\]
The first morphism here induces the modified Gysin map
\[ \gysc_{i_V} : H^{m+2q}(V,\zpn(q)_V) \lra H^{m+2(d+N)}_V(U,\zpn(d+N)_U) \]
in \eqref{eq2-5-5}. Thus the compatibility facts in Remark \ref{rem:smooth} implies the claim.
\end{pf}

\newpage
\section{The case of $p$-torsion in mixed characteristic $(0,p)$}\label{sect3}
\setcounter{thm}{0}
\setcounter{equation}{0}
\medskip
Let $S$ be the spectrum of a henselian discrete valuation ring $A$ with fraction field $K$ of characteristic zero and {\it perfect} residue field $k$ of characteristic $p>0$. Consider a diagram with cartesian squares
\[\xymatrix{
X_\eta \, \ar@{^{(}->}[r]^-{j_X} \ar[d]_{f_\eta} \ar@{}[rd]|{\square} & X \ar[d]_{f\hspace{-1pt}} \ar@{}[rd]|{\square} &  \ar@{_{(}->}[l]_-{i_X} \ar[d]_{f_s \hspace{-2pt}}  \,X_s \\ \eta \, \ar@{^{(}->}[r]^j & S & \ar@{_{(}->}[l]_i \,s\,,\hspace{-5pt}
} \]
where $\eta$ (resp.\ $s$) is the generic point (resp.\ closed point) of $S$, and $f$ is separated of finite type. Let $n$ be a positive integer, and let $\DC_\eta$ (resp.\ $\DC_s$) be the \'etale sheaf $\mupn$ on $\eta$ (resp.\ the constant \'etale sheaf $\zpn$ on $s$). We define
\begin{align*}
\DC_{X_\eta} & := Rf_\eta^! \DC_\eta \in D^+(X_{\eta,\et},\zpn),\\
\DC_{X_s} &:= Rf_s^! \DC_s \in D^+(X_{s,\et},\zpn).
\end{align*}
We recall some standard facts on $\DC_{X_\eta}$ (compare Theorems \ref{thm:jss} and \ref{thm:shiho} for $\DC_{X_s}$)
\subsubsection{}\label{rel:dual}
If $X_\eta$ is smooth over $\eta$ of pure dimension $d$, then there is a canonical isomorphism
\[ \tr^{f_\eta}: \mupn^{\otimes d+1}[2d] \isom \DC_{X_\eta} \]
in $D^+(X_{\eta,\et},\zpn)$ by the relative Poincar\'e duality \cite{sga4} XVIII.3.25.

\subsubsection{}\label{rem:lam}
For points $y \in (X_\eta)_q$ and $x \in (X_\eta)_{q-1}$ with $x \in \ol { \{ y \} }\subset X_\eta$, there is a commutative diagram
\[\xymatrix{
\iota_x^*R\iota_{y*}\mupn^{\otimes q+1}[2q] \ar[rr]^-{-\dKyx} \ar[d]_{\iota_x^*R\iota_{y*}(\tau_y)}^{\rwr}
&& \mupn^{\otimes q}[2q-1] \ar[d]_{\lwr}^{\tau_x[1]} \\
\iota_x^*R\iota_{y*}R\iota_y^!\DC_{X_\eta} \ar[rr]^-{\iota_x^*\{\delta^\loc_{y,x}(\DC_{X_\eta})\}}
&&  R\iota_x^!\DC_{X_\eta}[1]
}\]
in $D^+(x_\et,\zpn)$. Here for a point $v \in (X_\eta)_m$, $\iota_v$ denotes the canonical map $v \hra X_\eta$ and $\tau_v$ denotes the canonical isomorphism $\mupn^{\otimes m+1}[2m] \cong R\iota_v^!\DC_{X_\eta}$ obtained from \S\ref{rel:dual} for a smooth dense open subset of $\ol {\{ v \}}$. See \eqref{eq-0.1.1''} for the top arrow.
One can check this commutativity in the following way. Localizing and embedding $X_\eta$ into an affine space, we may suppose that $X_\eta$ is smooth. Because $R\iota_x^!\DC_{X_\eta}[1]$ (resp.\ $\iota_x^*R\iota_{y*}\mupn^{\otimes q+1}[2q]$) is concentrated in degree $-2q+1$ (resp.\ $\leq -2q+1$), the problem is reduced, by \S\ref{sect0-5-4}\,(1), to the commutativity at the $(-2q+1)$-st cohomology sheaves, which follows from Theorem {\rm\ref{prop:appA}} and \cite{sga4.5} Cycle, 2.3.8\,(i).

\smallskip
\subsection{Condition $\bs{\bK_q}$}\label{sect3-1}
The complexes $\DC_{X_\eta}$ and $\DC_{X_s}$ are important for the theory of duality and homology over $\eta$ and $s$, as we have seen in \S\ref{sect1} and \S\ref{sect2}. For working over $S$, we study morphisms
\[ Rj_{X*}\DC_{X_\eta} \lra i_{X*}\DC_{X_{s}}[-1]\,, \] see \S\ref{sect3-2} and \S\ref{sect3-6} below.
In particular, we want to investigate local conditions. For a point $v \in X$, let $i_v$ be the canonical map $v \hra X$. Let $q$ be a non-negative integer, and take points $y\in (X_\eta)_q$ and $x\in (X_s)_q$ with $x\in \overline{\{y\}} \subset X$. Put $Y:=\Spec(\O_{\ol { \{ y \} },x})$ and $x':=Y \times_X X_s$ and let $\pi : Y \hra X$ be the natural map. Then we have cartesian squares
\[\xymatrix{
y \, \ar@{^{(}->}[r]^-{j_Y} \ar[d]_{\iota_y} \ar@{}[rd]|{\square}
& Y \ar[d]_{\pi\hspace{-1pt}} \ar@{}[rd]|{\square}
& \ar@{_{(}->}[l]_-{i_Y} \ar[d]_{\ep_{x'} \hspace{-2pt}}  \,x' \\
X_\eta \, \ar@{^{(}->}[r]^{j_X}
& X
& \ar@{_{(}->}[l]_{i_X} \, X_s\,,\hspace{-5pt}
}\]
and a canonical nilpotent closed immersion $x \hra x'$.
Now let
\[ \delta_X:Rj_{X*}\DC_{X_\eta} \lra i_{X*}\DC_{X_s}[-1] \]
be a morphism in $D^+(X_\et,\zpn)$. Applying $R\pi_* R\pi^!$ to $\delta_X$, we obtain a morphism
\begin{equation}\label{mor1}
R\pi_{*}R\pi^!(\delta_X):Ri_{y*} R\iota_y^!\DC_{X_\eta} \lra Ri_{x*} R\ep_x^! \DC_{X_s}[-1]\,,
\end{equation}
where $\ep_x$ denotes the canonical map $x \hra X_s$, and we have used base-change isomorphisms
\[ R\pi^!Rj_{X*}=Rj_{Y*}R\iota_y^!\quad \hbox{ and } \quad R\pi^!i_{X*}=i_{Y*}R\ep_{x'}^!\,, \]
and the isomorphism
\[ R\ep_{x'*}R\ep_{x'}^!=R\ep_{x*}R\ep_{x}^! \]
by the invariance of \'etale topology. Furthermore, we have $R\iota_y^!\DC_{X_\eta} \cong \mupn^{\otimes q+1}[2q]$ by \S\ref{rel:dual}, and we have $R\iota_x^! \DC_{X_s} \cong \logwitt x n q [q]$ by Theorem \ref{thm:jss}. Therefore the morphism \eqref{mor1} is identified with a morphism $Ri_{y*} \mupn^{\otimes q+1}[2q] \ra  Ri_{x*} \logwitt x n q [q-1]$, which induces a map of cohomology sheaves in degree $-q+1$:
\[ \delta_X(y,x): R^{q+1}i_{y*} \mupn^{\otimes q+1} \lra i_{x*} \logwitt x n q \,. \]
We are going to compare this map of sheaves on $X_\et$ with Kato's residue map (cf.\ \eqref{eq-0.1.1'}):
\[ \dKyx : R^{q+1}i_{y*} \mupn^{\otimes q+1} \lra i_{x*} \logwitt x n q \,. \]
\stepcounter{thm}
\begin{defn}\label{def:K_q}
We say that $\delta_X$ satisfies $\bK_q$ if the induced map $\delta_X(y,x)$ agrees with $\dKyx$ for all points $y\in (X_\eta)_q$ and $x\in (X_s)_q$ with $x\in \overline{\{y\}}$.
\end{defn}
\begin{rem}
In view of \S{\rm\ref{sect0-5-4}\,(1)}, the morphism \eqref{mor1} is determined by $\delta_X(y,x)$. In fact, we have $R^mi_{y*} \mupn^{\otimes q+1}=0$ for any $m>q+1$ by a similar argument as for Lemma {\rm\ref{lem1}} below.
\end{rem}

\smallskip
\subsection{Functoriality of Kato's residue maps}\label{sect3-2}
Let
\[ \delta_S^\val:Rj_*\DC_\eta  \lra i_*\DC_s[-1]\]
be the composite morphism
\[ Rj_*\DC_\eta \cong \tau_{\leq 1} Rj_*\DC_\eta \lra R^1j_* \DC_\eta [-1] \lra i_*\DC_s[-1] \]
in $D^b(S_\et,\zpn)$, where the first isomorphism follows from a theorem of Lang: $\cd({\eta})=1$ (cf.\ Lemma \ref{lem1} below) and the last morphism is induced by Kummer theory and the normalized valuation $v_A$ on $K^{\times}$, i.e., Kato's residue map (cf.\ \S\ref{sect0-1}). By the base-change isomorphisms
\[ Rf^! Rj_* \DC_\eta = Rj_{X*}\DC_{X_\eta} \quad \hbox{ and } \quad Rf^! i_* \DC_s = i_{X*}\DC_{X_s}\,, \]
we obtain a morphism
\[ \delta_X^\Sval:=Rf^!(\delta_S^\val):Rj_{X*} \DC_{X_\eta} \lra i_{X*}\DC_{X_s}[-1] \quad\hbox{ in }\; D^b(X_\et,\zpn).\]
The first main result of this section is the following theorem:
\begin{thm}\label{Th.1-1}
\begin{enumerate}
\item[(1)]
The morphism $\delta_X^\Sval$ satisfies $\bK_q$ for all $q\ge 0$.
\item[(2)]
$\delta_X^\Sval$ is the only morphism that satisfies $\bK_q$ for all $q\ge 0$.
\item[(3)]
If $X_\eta$ is smooth of pure dimension $d$, $\delta_X^\Sval$ is the only morphism satisfying $\bK_d$.
\end{enumerate}
\end{thm}
\noindent
The proof of this result will be finished in \S\ref{sect3-8} below.
\smallskip
\subsection{First reductions}\label{sect3-3}
We first note that, to prove Theorem \ref{Th.1-1}, we may assume that $X$ is reduced and the closure of $X_\eta$. In fact, let $X' \subset X$ be the closure of $X_\eta$ with the reduced subscheme structure. Then we get cartesian squares
\[\xymatrix{
 (X_\eta)_{\text{red}} \, \ar@{^{(}->}[r]^-{j_{X'}} \ar[d]_{\kappa_\eta\hspace{-1pt}} \ar@{}[rd]|{\square}
& X' \ar[d]_{\kappa\hspace{-1pt}} \ar@{}[rd]|{\square}
& \ar@{_{(}->}[l]_-{i_{X'}} \ar[d]_{\kappa_s\hspace{-1pt}} \, X'_s  \\
X_\eta \, \ar@{^{(}->}[r]^-{j_X}
& X
& \ar@{_{(}->}[l]_-{i_X} \, X_s\,,\hspace{-5pt}
}\]
where $\kappa$ is the closed immersion. They induce a commutative diagram
\[\xymatrix{
\kappa_* Rj_{X'*}\DC_{X'_\eta} \ar[rr]^-{\kappa_*R\kappa^!(\delta_X)} \ar[d]_{\kappa_{\eta*}}^{\rwr}
&& \kappa_* i_{X'*}\DC_{X_s'} \ar[d]^{\kappa_{s*}} \\
Rj_{X*}\DC_{X_\eta} \ar[rr]^-{\delta_X}
&& i_{X*}\DC_{X_s}\,,\hspace{-5pt}
}\]
for any given morphism $\delta_X$ at the bottom. The left adjunction map is an isomorphism by topological invariance of \'etale cohomology. Moreover, $R\kappa^!(\delta_X^\Sval) = \delta_{X'}^\Sval$, and evidently $\delta_X$ satisfies $\bK_q$ if and only if $R\kappa^!(\delta_X)$ does. This shows that the claims of Theorem \ref{Th.1-1} hold for $X$ if and only if they hold for $X'$. We also note the following reduction:
\begin{lem}\label{Cl.1-3}
A morphism $\delta_X$ satisfies $\bK_q$ if and only if for all integral closed subschemes $\iota_Z : Z \hra X$ of dimension $q+1$ the morphism $R\iota_Z^!(\delta_X)$ satisfies $\bK_q$. In particular, Theorem {\rm\ref{Th.1-1}\,(1)} holds for $X$ if and only if $\delta_Z^\Sval$ satisfies $\bK_q$ for all integral subschemes $Z\subset X$ of dimension $d$.
\end{lem}
\begin{pf}
Let $X$ be arbitrary. Take a point $y \in (X_\eta)_q$ ($0 \le q \le d:=\dim(X_\eta)$), let $Z$ be its closure in $X$, and take an $x \in (X_s)_q$ with $x \in Z$. Let $\iota_Z: Z\hra X$ be the natural inclusion. We have base-change isomorphisms
\[ R\iota_Z^!Rj_{X*}\DC_{X_\eta}=Rj_{Z*}\DC_{Z_\eta}\quad \hbox{ and } \quad
 R\iota_Z^!i_{X*}\DC_{X_s}=i_{Z*}\DC_{Z_s}, \]
and it follows from the definitions in \S\ref{sect3-1} that $\delta_X(y,x) = (R\iota_Z^!(\delta_X))(y,x)$, if we regard these as maps $\iota_x^* R^{q+1}\iota_{y*} \mu^{\otimes{q+1}}_{p^n} \to \logwitt x n q$. This shows the first claim. The second claim follows, because $R\iota_Z^!(\delta_X^\Sval) = \delta_Z^\Sval$.
\end{pf}
Finally we note:
\begin{rem}\label{rem1}
To prove that $\delta_X^\Sval$ satisfies $\bK_q$ it suffices to assume that $f:X \to S$ is proper by taking a compactification of $f$.
\end{rem}

\smallskip
\subsection{Criterion in the proper case}\label{sect3-4}
Suppose that we are given two morphisms
\begin{align*}
\delta_S & :  Rj_*\DC_\eta \lra i_*\DC_s[-1] \quad \hbox{ in } \; D^+(S_\et,\zpn), \\
\delta_X & :  Rj_{X*}\DC_{X_\eta} \lra i_{X*}\DC_{X_s}[-1] \quad \hbox{ in } \; D^+(X_\et,\zpn).
\end{align*}
Assuming that $f$ is proper, we give a simple criterion as to when $\delta_X$ agrees with $Rf^!(\delta_S)$, that is, as to when the following diagram commutes in $D^+(X_\et,\zpn)$:
\begin{equation}\label{dg:proper}
\xymatrix{
Rj_{X*}\DC_{X_\eta} \ar[rr]^-{\delta_X} \ar@{=}[d] && i_{X*}\DC_{X_s}[-1] \ar@{=}[d] \\
Rf^!Rj_*\DC_\eta  \ar[rr]^-{Rf^!(\delta_S)} && Rf^!i_*\DC_s[-1]\,,\hspace{-5pt}
}
\end{equation}
where the equalities mean the identifications by base-change isomorphisms.
\stepcounter{thm}
\begin{prop}\label{Lem.1-1}
Suppose that $f$ is proper. Then the diagram \eqref{dg:proper} commutes if and only if the following diagram is commutative$:$
\stepcounter{equation}
\begin{equation}\label{dg:proper2}
\xymatrix{
\H^1(X_{\wt{\eta}},\DC_{X_\eta}) \ar[rr]^-{\delta_X} \ar[d]_{f_*} && \H^0(X_{\ol s},\DC_{X_s}) \ar[d]^{f_*} \\
\H^1(\wt{\eta},\DC_\eta)  \ar[rr]^-{\delta_S} && \H^0(\ol s,\DC_s)\,,\hspace{-5pt}
}
\end{equation}
where $\wt{\eta}$ denotes the generic point of the maximal unramified extension $\wt{S}$ of $S$ {\rm(}$\ol s$ is the closed point of $\wt{S}${\rm);} the vertical maps are defined by the adjunction map $Rf_!Rf^! \ra \id$ and the properness of $f$, that is, $Rf_!=Rf_*$.
\end{prop}
\begin{pf}
By the adjointness between $Rf^!$ and $Rf_!$, we have the adjunction maps $f^! : \id \ra Rf^!Rf_!$ and $f_! : Rf_!Rf^! \ra \id$,
which satisfy the relation that the composite
\[ Rf^! \os{f^!}\lra Rf^!Rf_!Rf^! \os{f_!}\lra Rf^! \]
is the identity map. By these facts, it is easy to see that the commutativity of \eqref{dg:proper} is equivalent to that of the following diagram in $D^+(S_\et,\zpn)$:
\begin{equation}\label{dg:proper3}
\xymatrix{
Rf_!Rj_{X*}\DC_{X_\eta} \ar[rr]^-{Rf_!(\delta_X)} \ar[d]_\alpha && Rf_!i_{X*}\DC_{X_s}[-1] \ar[d]^\beta \\
Rj_*\DC_\eta  \ar[rr]^-{\delta_S} && i_*\DC_s[-1]\,,\hspace{-5pt}
}
\end{equation}
where $\alpha$ is defined as the composite
\[\xymatrix{
\alpha : Rf_!Rj_{X*}\DC_{X_\eta} \ar@{=}[rr]^-{\text{base-change}} && Rf_!Rf^!Rj_*\DC_\eta \ar[r]^-{f_!} & Rj_*\DC_\eta
}\]
and $\beta$ is defined in a similar way (note that we do not need the properness of $f$ for this equivalence). We prove that the commutativity of \eqref{dg:proper3} is equivalent to that of \eqref{dg:proper2}.
For this, we first show the following:

\par\medskip\noindent
{\it Claim. $i^*Rf_!Rj_{X*}\DC_{X_\eta}=i^*Rf_*Rj_{X*}\DC_{X_\eta}$ is concentrated in degrees $\le 1$.}

\par\medskip\noindent
{\it Proof of Claim}.
Because the stalk at $\ol s$ of the $m$-th cohomology sheaf is
\[ \cH^m(i^*Rf_*Rj_{X*} \DC_{X_\eta})_{\ol s} \cong \H^m(X_{\wt{\eta}},\DC_{X_\eta}) \]
by the properness of $f$, it suffices to show that the group on the right hand side is zero for $m>1$. Take an open subset $U_\eta\subset X_\eta$ which is smooth over $\eta$ of pure dimension $d:=\dim(X_\eta)$ and such that $\dim(Z_\eta)<d$, where $Z_\eta$ denotes the closed complement $X_\eta \ssm U_\eta$. By \eqref{eq-0.4.2}, there is a localization exact sequence
\begin{equation}\label{exact:local}
\dotsb \lra \H^m(Z_{\wt{\eta}},\DC_{Z_\eta}) \lra \H^m(X_{\wt{\eta}},\DC_{X_\eta}) \lra \H^m(U_{\wt{\eta}},\DC_{U_\eta}) \lra \dotsb.
\end{equation}
Now we have $\DC_{U_\eta} \cong \mupn^{\otimes d+1} [2d]$ by \S\ref{rel:dual} for $U_\eta$, so that
\[ \H^m(U_{\wt{\eta}},\DC_{U_\eta}) \cong \H^{m+2d}(U_{\wt{\eta}},\mupn^{\otimes d+1})\,, \]
which vanishes for $m>1$ because $\cd(U_{\wt{\eta}}) \leq 2d+1$, cf.\ the proof of Lemma \ref{lem1} below. Thus the vanishing of $\H^m(X_{\wt{\eta}},\DC_{X_\eta})$ for $m>1$ is shown by induction on $\dim(X_\eta)$ and we obtain the claim.
\qed
\par\medskip \noindent
We turn to the proof of Proposition \ref{Lem.1-1}. By the above claim and \S\ref{sect0-5-4}\,(1), a morphism $i^*Rf_!Rj_{X*}\DC_{X_\eta} \ra \DC_s[-1]$ is determined by the map of the $1$st cohomology sheaves, and thus determined by the associated map of their stalks at $\ol s$. Hence by the adjointness between $i_*$ and $i^*$, the diagram \eqref{dg:proper3} commutes if and only if the diagram \eqref{dg:proper2} does. This completes the proof of Proposition \ref{Lem.1-1}.
\end{pf}

\smallskip
\subsection{Result for smooth generic fiber}\label{sect3-5}
In Proposition \ref{Cl.1-2} below we obtain a first step towards part (3) of Theorem \ref{Th.1-1} which will also be used for the other parts. We first show:
\begin{lem}\label{lem1}
Let $\cF$ be a torsion sheaf on $(X_\eta)_\et$. Then $R^mj_{X*}\cF=0$ for any $m > \dim(X_\eta)+1$.
\end{lem}
\begin{pf}
Clearly $R^mj_{X*}\cF$ is trivial on $X_\eta$ if $m>0$. Hence the problem is \'etale local on $X_s$ and we may suppose that $s = \ol s$. Let $x$ be a point on $X_s$. The stalk of $R^mj_{X*}\cF$ at $\ol x$ is isomorphic to $\H^m(\Spec(\O_{X,\ol x}^{\sh}[p^{-1}]),\cF)$, where $\Spec(\O_{X,\ol x}^{\sh}[p^{-1}])$ is written as a projective limit of affine varieties over $\eta$ of dimension $\leq \dim(X_\eta)$.
Hence the assertion follows from the affine Lefschetz theorem (\cite{sga4} XIV.3.2) and Lang's theorem: $\cd(\eta)=1$ (\cite{se} II.3.3).
\end{pf}
\begin{prop}\label{Cl.1-2}
If $X_\eta$ is smooth of pure dimension $d$, then there exists a unique morphism $\delta_X:Rj_{X*}\DC_{X_\eta} \to i_{X*}\DC_{X_s}[-1]$ satisfying $\bK_d$.
\end{prop}
\begin{pf}
By \S\ref{sect3-3} we may assume that $\dim(X_s)\leq d$. We have $\DC_{X_\eta} \cong \mupn^{\otimes d+1}[2d]$ by \S\ref{rel:dual}, and $Rj_{X*} \DC_{X_\eta}$ is concentrated in $[-2d,-d+1]$ by Lemma \ref{lem1}. On the other hand, $i_{X*}\DC_{X_s}[-1]$ is concentrated in degree $[-d+1,1]$ by Theorem \ref{thm:jss} and the assumption $\dim(X_s) \leq d$. Hence a morphism $\delta_X:Rj_{X*}\DC_{X_\eta} \to i_{X*}\DC_{X_s}[-1]$ is determined by the map $\cH^{-d+1}(\delta_X)$ of the $(-d+1)$-st cohomology sheaves by \S\ref{sect0-5-4}\,(1). Moreover, for a given $\delta_X$, there is a commutative diagram of sheaves on $X_\et$:
\vspace{-10pt}
\[\xymatrix{
\cH^{-d+1}(Rj_{X*}\DC_{X_\eta}) \ar@{=}[r] \ar[d]_{\cH^{-d+1}(\delta_X)}
& R^{d+1}j_{X*} \mupn^{\otimes d+1} \ar[r]^-{\alpha}
& \displaystyle \bigoplus_{y \in (X_\eta)_d}^{\phantom{|^a}} \ R^{d+1}i_{y*} \mupn^{\otimes d+1} \ar[d]_-{\gamma}\\
\us{\phantom{|^{|^|}}} \cH^{-d+1}(i_{X*}\DC_{X_s}[-1]) \ar@<5pt>@{=}[r]
& \us{\phantom{|^a}}{i_{X*}\cH^{-d}(\DC_{X_s})} \; \ar@<5pt>@{^{(}->}[r]^-{\beta}
& \displaystyle \bigoplus_{x \in (X_s)_d} \ i_{x*}\logwitt x n d \,,\hspace{-5pt}
}\]
where $\alpha$ is the adjunction map, $\beta$ is an inclusion obtained from Theorem \ref{thm:jss} and $\gamma$ is the sum of $\delta_X(y,x)$'s. These facts show the uniqueness of $\delta_X$ satisfying $\bK_d$. Next we prove its existence. For this, let us consider the following diagram of sheaves:
\vspace{-10pt}
\[\xymatrix{
& R^{d+1}j_{X*} \mupn^{\otimes d+1} \ar[r]^-{\alpha} \ar@{.>}[d]_{\varrho}
& \displaystyle
 \bigoplus_{y \in (X_\eta)_d}^{\phantom{|^a}} \ R^{d+1}i_{y*} \mupn^{\otimes d+1} \ar[d]_-{\partial_2} \ar[r]^-{\partial_1}
& \displaystyle \bigoplus_{w \in (X_\eta)_{d-1}}^{\phantom{|^a}} \ R^{d}i_{w*} \mupn^{\otimes d} \ar[d]_-{\partial_3} \\
\us{\phantom{|^a}}0 \ar@<5pt>[r]
& \us{\phantom{|^a}}{i_{X*}\cH^{-d}(\DC_{X_s})} \; \ar@<5pt>[r]^-{\beta}
& \displaystyle \bigoplus_{x \in (X_s)_d} \ i_{x*}\logwitt x n d \ar@<5pt>[r]^-{\partial_4}
& \displaystyle \bigoplus_{z \in (X_s)_{d-1}} \ i_{z*}\logwitt z n {d-1} \,,\hspace{-5pt}
}\]
where $\alpha$ and $\beta$ are the same maps as above, and each $\partial_i$ ($i=1,\dotsc,4$) is the sum of Kato's residue maps. We have the following facts for this diagram: the right square is anti-commutative by \cite{kk:hasse} 1.7 for $X$; the upper row is a complex by \S\ref{rem:lam}; the lower row is exact by Theorem \ref{thm:jss}. Hence $\partial_2$ induces a map $\varrho$ as in the diagram, and we obtain a morphism $\delta_X$ satisfying $\bK_d$ by extending this map (cf.\ \S\ref{sect0-5-4}\,(1)). This completes the proof.
\end{pf}

\smallskip
\subsection{Case of points}\label{sect3-6}
We will prove Theorem \ref{Th.1-1}\,(1) by induction on $\dim(X_\eta)$. We start with:
\begin{lem}\label{Cl.1-1}
Theorem {\rm\ref{Th.1-1}} is true for $X$ with $\dim(X_\eta)=0$.
\end{lem}
\begin{pf}
First we show \ref{Th.1-1}\,(1). By Lemma \ref{Cl.1-3} and Remark \ref{rem1} we may assume that $X$ is integral and proper. Then $f:X \ra S$ is flat and finite by \cite{ega3} 4.4.2, and moreover, $X_s$ is irreducible because $S$ is henselian and $X$ is irreducible. Let $j':\eta' \hra X$ (resp.\ $i':s' \hra X$) be the generic (resp.\ closed) point. Then $\eta'=X_\eta$ and $\eta' \to \eta$ is finite \'etale, because $X$ is integral and $\ch(K)=0$. On the other hand, $s'\to s$ is finite \'etale as well by the perfectness of $k$, and this map factors as the composite of a nilpotent closed immersion $s' \hra X_s$ with $f_s:X_s \ra s$. Therefore we have $\DC_{X_\eta}=\mupn$ and $\DC_{X_s}=\zpn$. Now let
\[ \delta_X:Rj_{X*}\DC_{X_\eta} \lra i_{X*}\DC_{X_s}[-1] \]
be the composite morphism $Rj_{X*}\mupn \ra R^1j_{X*}\mupn[-1] \ra i_{X*}\zpn[-1]$, where the last morphism is given by the map $\dval_{\eta',s'}$. Because $\delta_X$ satisfies ${\bK}_0$ by definition, our task is to show the equality $\delta_X=\delta_X^\Sval (:=Rf^!(\delta_S^\val))$. Moreover, by the finiteness of $f$ and Proposition \ref{Lem.1-1}, we have only to show the commutativity of the diagram
\begin{equation}\label{com1}
\xymatrix{
\H^1(\eta',\mupn) \ar[rd]^-{\dval_{\eta',s}} \ar[d]_{\tr_f} \\
\H^1(\eta,\mupn) \ar[r]_-{\dval_{\eta,s}} & \H^0(s,\zpn)\,,\hspace{-5pt}
}\end{equation}
 assuming that $s=s'=\ol s$ (that is, $k$ is algebraically closed).
We show this commutativity. Let $B_0$ be the affine ring of $X$, let $B$ be the normalization of $B_0$, let $L$ be the fraction field of $B$ and let $x$ be the closed point of $\Spec(B)$. By definition, $\dval_{\eta',s}$ is the composite
\[ \H^1(\eta',\mupn) \lra \H^0(x,\zpn) \isom \H^0(s,\zpn)\,. \]
where the first map is given by the normalized valuation $v_B$ on $L^{\times}$ and the second map is induced by the isomorphism $x \cong  s$. On the other hand, there is a commutative diagram
\[ \xymatrix{
L^{\times}/p^n \ar[r]^-{\cong} \ar[d]_{N_{L/K}} &\H^1(\eta',\mupn) \ar[d]^{\tr_f} \\
K^{\times}/p^n \ar[r]^-{\cong} & \H^1(\eta,\mupn)\,,\hspace{-5pt}
}\]
where $N_{L/K}$ denotes the norm map (cf.\ \cite{sga4} XVIII.2.9 (Var 4)), and the horizontal arrows are boundary maps coming from the Kummer theory for $\eta'$ and $\eta$, respectively. Therefore the commutativity of \eqref{com1} follows from the fact that $v_B = v_A \circ N_{L/K}$. Now we prove the other parts of Theorem \ref{Th.1-1} for $X$. By \S\ref{sect3-3} we may assume that $X$ is reduced. Then, since $\dim(X_\eta)=0$, $X_\eta$ is smooth, and Proposition \ref{Cl.1-2} implies that $\delta_X^\Sval$ is the only morphism satisfying $\bK_0$.
\end{pf}

\smallskip
\subsection{Induction step}\label{sect3-7}
Consider the following situation.
Suppose that $X$ is reduced, separated of finite type over $S$, that $X_\eta$ has dimension $d \ge 1$, and that $X_\eta$ is dense in $X$. Choose a {\it smooth affine} dense open subset $U_\eta \subset X_\eta$. Let $Z_\eta:=X_\eta \ssm U_\eta$ with the reduced structure, let $Z$ be the closure of $Z_\eta$ in $X$, and let $U = X \ssm Z$. Then the composite morphism $f_Z: Z \ra X \ra S$ is flat, and hence we have
\begin{equation}\label{cond1}
(U_s)_d=(X_s)_d\,.
\end{equation}
We name the canonical immersions as follows:
\[ \xymatrix{
U_\eta \, \ar@{^{(}->}[r]^-{j_U} \ar@{_{(}->}[d]_{\phi_\eta} \ar@{}[rd]|{\square}
& \, U_{\phantom.} \ar@{_{(}->}[d]_\phi \ar@{}[rd]|{\square}
& \ar@{_{(}->}[l]_-{i_U} \ar@{_{(}->}[d]_{\phi_s} \, U_s \\ X_\eta \, \ar@{^{(}->}[r]^-{j_X}
& X
& \ar@{_{(}->}[l]_{i_X} \, X_s \\
Z_\eta \phantom{|}\hspace{-2pt} \ar@{^{(}->}[r]^-{j_Z} \ar@{^{(}->}[u]^{\psi_\eta} \ar@{}[ru]|{\square}
& Z \phantom{|}\hspace{-2pt} \ar@{^{(}->}[u]^\psi \ar@{}[ru]|{\square}
& \ar@{_{(}->}[l]_-{i_Z} \ar@{^{(}->}[u]^{\psi_s} \, Z_s\,. \phantom{|}\hspace{-5pt}
} \]
Consider a diagram of the following type in $D^+(X_\et,\zpn)$:
\medskip
\begin{equation}\label{com2}
\xymatrix{
\psi_*Rj_{Z*} \DC_{Z_\eta} \ar[r]^-{\psi_*} \ar[d]_{\delta_Z^S} \ar@{}[rd]|{(1)}
& Rj_{X*} \DC_{X_\eta} \ar[r]^-{\phi^*} \ar@{.>}[d]_{\delta_1} \ar@{}[rd]|{(2)}
& R\phi_* Rj_{U*} \DC_{U_\eta} \ar[r]^-{-\ep_1} \ar@{.>}[d]_{\delta_2} \ar@{}[rd]|{(3)}
& \psi_* Rj_{Z*} \DC_{Z_\eta}[1] \ar[d]_{\delta_Z^S[1]} \\
\psi_* i_{Z*} \DC_{Z_s}[-1] \ar[r]^-{\psi_*} & i_{X*} \DC_{X_s}[-1] \ar[r]^-{\phi^*}
& R\phi_* Ri_{U*} \DC_{U_s}[-1] \ar[r]^-{\ep_2[-1]} & \psi_* i_{Z*} \DC_{Z_s}\,.\hspace{-5pt}
}
\end{equation}
Here we put
\[ \delta_Z^S := \psi_*(\delta_Z^\Sval)=\psi_* Rf_Z^!(\delta_S^\val)\,, \quad
\ep_1 :=\delta_{U,Z}^{\loc}(Rj_{X*} \DC_{X_\eta})\,, \quad \ep_2 :=\delta_{U,Z}^{\loc}(i_{X*} \DC_{X_s})\,, \]
the horizontal rows are the distinguished triangles deduced from the obvious localization triangles (cf.\ \eqref{eq-0.4.2}) and the base-change isomorphisms
\begin{align*}
R\psi^! Rj_{X*} \DC_{X_\eta} = Rj_{Z*} \DC_{Z_\eta}\,,\qquad & \phi^* Rj_{X*} \DC_{X_\eta} = Rj_{U*} \DC_{U_\eta}\,,\\
R\psi^! i_{X*} \DC_{X_s} = i_{Z*} \DC_{Z_s}\,,\qquad & \phi^* i_{X*} \DC_{X_s} = Ri_{U*} \DC_{U_s}\,.
\end{align*}
\addtocounter{thm}{2}
\begin{lem}\label{lem:delta1}
If $\delta_2$ is given, there is at most one morphism $\delta_1$ making the squares {\rm(1)} and {\rm(2)} in \eqref{com2} commutative.
\end{lem}
\begin{pf}
We want to apply Lemma \ref{lem.derived}\,(3). Because $U_\eta$ is smooth and affine, we have $\DC_{U_\eta} \cong \mupn^{\otimes d+1}[2d]$ by \S\ref{rel:dual}, and $C = R (\phi j_U)_* \DC_{U_\eta}$ is concentrated in $[-2d,-d+1]$ by a similar argument as for Lemma \ref{lem1}. On the other hand, because $\dim(X_s) \leq d$, $i_{X*}\DC_{X_s}$ is concentrated in $[-d,0]$ by Theorem \ref{thm:jss} (note that $i_X$ is a closed immersion). Similarly, $A' = \psi_* i_{Z*}\DC_{Z_s}$ is concentrated in $[-d+1,0]$, because we have $\dim(Z_s) \leq d-1$ by the flatness of $f_Z:Z \ra S$. Therefore we get
\[ \Hom_{D(X,\zpn)}(C,A') = 0. \]
On the other hand, for $A = \psi_* Rj_{Z*} \DC_{Z_\eta}$ and $C' = R\phi_* Ri_{U*} \DC_{U_s}[-1]$ we have
{\allowdisplaybreaks
\begin{align*}
\Hom^{-1}_{D(X,\zpn)}(A,C') & = \Hom_{D(X,\zpn)}(\psi_* Rj_{Z*} \DC_{Z_\eta},R\phi_* Ri_{U*} \DC_{U_s}[-2]) \\
 & = \Hom_{D(X,\zpn)}(\phi^* \psi_* Rj_{Z*} \DC_{Z_\eta}, Ri_{U*} \DC_{U_s}[-2]) & \hbox{(adjunction)} \\
 & = 0  &\hbox{($\phi^* \psi_* = 0$)}
\end{align*}
}So the lemma follows from Lemma \ref{lem.derived}\,(3).
\end{pf}

\begin{lem}\label{Cl.1-5}
Consider the diagram \eqref{com2} and assume that $\delta_2 = R\phi_*(\delta_U)$ where $\delta_U: Rj_{U*} \DC_{U_\eta} \to Ri_{U*} \DC_{U_s}$ denotes the morphism obtained by applying Proposition {\rm\ref{Cl.1-2}} to $U$. Assume that $\bK_{d-1}$ holds for $\delta_Z^\Sval$.
\begin{enumerate}
\item[(i)] If $X$ is integral, then the square {\rm(3)} in $\eqref{com2}$ commutes. Consequently, there exists a morphism $\delta_1$ which makes the squares {\rm(1)} and {\rm(2)} in \eqref{com2} commutative at the same time.
\item[(ii)] If $f: X \to S$ is proper, then any morphism $\delta_1$ making the square {\rm(1)} in \eqref{com2} commutative necessarily coincides with $\delta_X^\Sval$.
\end{enumerate}
\end{lem}
\begin{pf*}{\it Proof of Lemma \ref{Cl.1-5}}
(i) As we have seen in the proof of lemma \ref{lem:delta1}, $R(\phi j_U)_* \DC_{U_\eta}$ is concentrated in $[-2d,-d+1]$ and $\psi_* i_{Z*}\DC_{Z_s}$ is concentrated in $[-d+1,0]$.
By these facts, the square (3) commutes if and only if the square of the induced homomorphisms on the $(-d+1)$-st cohomology sheaves commutes. We prove this commutativity on cohomology sheaves. By Theorem \ref{thm:jss}, we have
\[ \cH^{-d+1}(\psi_* i_{Z*}\DC_{Z_s})= \psi_* i_{Z*} \cH^{-d+1} (\DC_{Z_s}) \hra
 \bigoplus_{x\in (Z_s)_{d-1}} \ i_{x*}\logwitt x n {d-1} \,. \]
Hence we may suppose that $(Z_s)_{d-1}$ is not empty, and the problem is local at each point in $(Z_s)_{d-1}$. Now fix a point $x \in (Z_s)_{d-1}$, and define $B$ (resp.\ $C$, $D$) as $\Spec(\O_{X,x})$ (resp.\ $U \times _X B$, $Z \times _X B$), and let $\sigma$ be the open immersion $C_\eta \hra B$. Note that $B$ is integral local of dimension two and that $D_\eta$ and $E:=(C_s)_{\red}$ are finite sets of points in $B^1 \subset X^1$. Our task is to show the commutativity of the following diagram on $B_\et$:
\vspace{-10pt}
\addtocounter{equation}{2}
\begin{equation}\label{com3}
\xymatrix{
R^{d+1}\sigma_*\mupn^{\otimes d+1} \ar[r]^-{\delta_3} \ar[d]_{\delta_4}
& \displaystyle \bigoplus_{z \in D_\eta}^{\phantom{|}} \ R^d i_{z*}\mupn^{\otimes d} \ar[d]^-{\delta_5} \\
\displaystyle \bigoplus_{y \in E} \ i_{y*}\logwitt y n d \ar@<6pt>[r]^-{\delta_6}
& \us{\phantom{l}}{i_{x*}\logwitt x n {d-1}}\,,\hspace{-5pt}
}
\end{equation}
where for a point $v \in B$, we wrote $i_v$ for the map $v \ra B$ and we have used the isomorphisms $\DC_{Z_\eta}\vert_{D_\eta} \cong \mupn^{\otimes d}[2(d-1)]$ (cf.\ \S\ref{rel:dual}) and $\DC_{U_s}\vert_{C_s} \cong \logwitt E n d$ (cf.\ Theorem \ref{thm:jss}). Each $\delta_i$ $(i=3,\dotsc,6)$ denotes the map obtained by restricting the corresponding morphism in the square (3) of \eqref{com2}. Now let $w$ be the generic point of $B$ and let $\alpha$ be the adjunction map $R^{d+1}\sigma_* \mupn^{\otimes d+1} \ra R^{d+1}i_{w*} \mupn^{\otimes d+1}$. We have the following facts for the maps in \eqref{com3}.
\begin{itemize}
\item $\delta_3$ factors, by \S\ref{rem:lam}, as
\[
\delta_3: R^{d+1}\sigma_* \mupn^{\otimes d+1} \os{\alpha}
\lra R^{d+1}i_{w*} \mupn^{\otimes d+1}
\os{\bigoplus_z \, \dKwz}{\llllra} \bigoplus_{z \in D_\eta} \ R^d i_{z*} \mupn^{\otimes d} \,.
\]
\item $\delta_4$ factors as
\[
\delta_4: R^{d+1}\sigma_* \mupn^{\otimes d+1}
\os {\alpha}\lra R^{d+1}i_{w*} \mupn^{\otimes d+1}
\os{\bigoplus_y \, \dKwy}\llllra \bigoplus_{y \in E} \ i_{y*} \logwitt y n d
\]
by the construction of $\delta_U$ (cf.\ Proposition \ref{Cl.1-2}).
\item $\delta_5= \sum_{z \in D_\eta} \ \dKzx$ \; by the assumption of the lemma.
\item $\delta_6= - \sum_{y \in E} \ \dKyx$ \; by Theorem \ref{thm:jss} and the construction of $\cM_{n,X_s}^\bullet$, cf.\ \S\ref{sect2-5}.
\end{itemize}
Therefore we obtain the commutativity of \eqref{com3} from a result of Kato \cite{kk:hasse} 1.7 for $B$, by noting that $B^1=(C_\eta)^1 \coprod D_\eta \coprod E$ and that $\Image(\alpha)$ is contained in the kernel of the map
\[
\bigoplus_v \ \dval_{w,v} :
R^{d+1}i_{w*} \mupn^{\otimes d+1} \lra \bigoplus_{v \in (C_\eta)^1} \ R^d i_{v*} \mupn^{\otimes d}
\]
(cf.\ proof of Proposition \ref{Cl.1-2}). This completes the proof of Lemma \ref{Cl.1-5}\,(i), because its second claim follows with \S\ref{sect0-5-2}.
\par
(ii) By the properness of $f$ and Proposition \ref{Lem.1-1}, we have only to show the commutativity of the right square (1)$''$ of the following diagram, assuming that $s=\ol s$:
\begin{equation}\label{com4}
\xymatrix{
\H^1(Z_\eta,\DC_{Z_\eta}) \ar[r]^-{\psi_*} \ar[d]_{Rf_*(\delta_Z^S)\hspace{-1.5pt}} \ar@{}[rd]|{\text{(1)$'$}}
 & \H^1(X_\eta,\DC_{X_\eta}) \ar[r]^-{f_*} \ar[d]_{Rf_*(\delta_1)\hspace{-1.5pt}} \ar@{}[rd]|{\text{(1)$''$}}
 & \H^1(\eta ,\DC_\eta) \ar[d]_{\delta_S^\val \hspace{-1.5pt}} \\
\H^0(Z_{s},\DC_{Z_s}) \ar[r]^-{\psi_*}
 & \H^0(X_{s},\DC_{X_s}) \ar[r]^-{f_*}
 & \H^0(s,\DC_{s})\,,\hspace{-5pt}
} \end{equation}
where for a proper morphism $g$ of schemes, we wrote $g_*$ for the adjunction map $Rg_*Rg^! \to \id$. The outer square of this diagram commutes, because $\delta_Z^S=\psi_* Rf_Z^!(\delta_S^\val)$ and the composite
\[ Rf_{Z*}Rf_Z^!=Rf_*\psi_*R\psi^!Rf^! \os{\psi_*}\lra Rf_*Rf^! \os{f_*}\lra \id \]
is functorial (in fact, this coincides with $f_{Z*}$). On the other hand, the square (1)$'$ commutes, because $\delta_1$ makes the square (1) in \eqref{com2} commutative. Moreover, in view of the exact sequence \eqref{exact:local}, the upper horizontal arrow $\psi_*$ in (1)$'$ is surjective, because we have
\[ \H^1(U_\eta,\DC_{U_\eta}) \cong \H^{2d+1}(U_\eta,\mupn^{\otimes d+1})=0 \]
by the assumptions that $s=\ol s$ and that $U_\eta$ is smooth affine of dimension $d \geq 1$ (cf.\ Lemma \ref{lem1}).
Therefore (1)$''$ is commutative, and we obtain Lemma \ref{Cl.1-5}\,(ii).
\end{pf*}

\smallskip
\subsection{Proof of Theorem \ref{Th.1-1}}\label{sect3-8}
First consider Theorem \ref{Th.1-1}\,(1). By Lemma \ref{Cl.1-3} and Remark \ref{rem1} it suffices to show:
\par\medskip\noindent
\quad\quad $(\sharp)$ \;
For integral $X$, $\delta_X^\Sval :=Rf^!(\delta_S^\val)$ satisfies $\bK_d$ with $d:=\dim(X_\eta)$.
\par\medskip\noindent
We show this property by induction on $d = \dim(X_\eta)$. The case $d=0$ is settled by Lemma \ref{Cl.1-1}. Now let $\dim(X_\eta)\ge 1$ and choose $U$ and $Z=X \ssm U$ as in \S\ref{sect3-7}. Assume that $\delta_2 = R\phi_*(\delta_U)$ with $\delta_U$ as in Lemma \ref{Cl.1-5}. The assumption of this lemma holds because $(\sharp)$ holds for $Z$ by induction assumption. Therefore there is a morphism $\delta_1$ making \eqref{com2} commutative, and this morphism is $\delta_1 = \delta_X^\Sval$. We conclude that $\delta_U = \phi^*(\delta_X^\Sval) = \delta_U^\Sval$. Hence $\delta_U^\Sval$ satisfies $\bK_q$ (by choice of $\delta_U$), and $\delta_X^\Sval$ satisfies $\bK_q$ as well, because $(X_\eta)_d = (U_\eta)_d$ by density of $U$ in $X$, and $(X_s)_d = (U_s)_d$ as noted in \eqref{cond1}.

Theorem \ref{Th.1-1}\,(3) now follows from Proposition \ref{Cl.1-2}, because $\delta_X^\Sval$ satisfies $\bK_d$.

Theorem \ref{Th.1-1}\,(2) follows once more by induction on $d=\dim(X_\eta)$, the case $d=0$ being given by Lemma \ref{Cl.1-1}. If $d\ge 1$ we may assume that $X$ is reduced and then again choose $U$ and $Z=X \ssm U$ as in \S\ref{sect3-7}. Assume that a morphism
\[ \delta_1: Rj_{X*}\DC_{X_\eta} \lra i_{X*}\DC_{X_s}[-1] \]
satisfies $\bK_q$ for all $q\geq 0$. Then $R\psi^!(\delta_1)$ satisfies $\bK_q$ for all $q\geq 0$ and agrees with $\delta_Z^\Sval$ by induction assumption. On the other hand, $\phi^*(\delta_1)$ satisfies $\bK_d$ and thus coincides with $\delta_U^\Sval$ by Theorem \ref{Th.1-1}\,(3) just proved. The conclusion is that $\delta_1$ makes the square (1) of \eqref{com2} commutative with $\delta_Z^S = \psi_*(\delta_Z^\Sval)$, and the square (2) of \eqref{com2} commutative with $\delta_2 = R\phi_*(\delta_U^\Sval)$. Since obviously $\delta_X^\Sval$ makes these diagrams commutative as well, Lemma \ref{lem:delta1} implies $\delta_1 =\delta_X^\Sval$ as wanted. This concludes the proof of Theorem \ref{Th.1-1}. \qed

\smallskip
\subsection{Dualizing complexes}\label{sect3-9}
We apply our results to the study of dualizing complexes as indicated in part \S\ref{sect0-2} of the introduction. Let $f: X \to S$ be separated of finite type, and define
\[ \DC_X := Rf^!\zpn(1)'_S \in D^+(X_\et,\zpn). \]
Also, let $\DC_{X_\eta} = Rf^!_\eta \mupn$ and $\DC_{X_s} = Rf^!_s\zpn$, as we defined at the beginning of this section. Then, by applying $Rf^!$ to the exact triangle \eqref{eq3-9-2} and using the base-change isomorphisms as in \eqref{dg:proper} we get a canonical isomorphism of exact triangles
\begin{equation}\label{triangles.DC-X}
\xymatrix{
i_{X *}\DC_{X_s}[-2] \ar[r]^-{g_X} \ar@{=}[d]
& \DC_X \ar[r]^-{t_X} \ar@{=}[d]_{\text{(def)}}
& Rj_{X*}\DC_{X_\eta} \ar[r]^-{\delta^\Sval_X} \ar@{=}[d]
& i_{X*}\DC_{X_s}[-1] \ar@{=}[d] \\
Rf^!i_*\zpn[-2] \ar[r]^-{Rf^!(g)}
& Rf^!\zpn(1)'_S \ar[r]^-{Rf^!(t)}
& Rf^!Rj_*\mupn \ar[r]^-{Rf^!(\delta^\val_S)}
& Rf^!i_*\zpn[-1]\,,\hspace{-5pt}
}
\end{equation}
where $g_X$ and $t_X$ are the adjunction maps for $i_X$ and $j_X$, respectively. By Theorem \ref{Th.1-1} the morphism $\delta_X^\Sval$ satisfies the localization property $\bK_q$ for all $q\geq 0$ (i.e., is locally given by Kato's residue maps), and is determined by this property (and just by $\bK_d$ if $X_\eta$ is smooth of dimension $d$). Moreover, by Lemma \ref{lem.deltaloc} below (see also \eqref{eq4-2-1} below), we see that
\begin{equation}\label{eq3-9-5}
\delta^\Sval_X = -\delta^\loc_{X_\eta,X_s}(\DC_X)\,.
\end{equation}
Because the dualizing complex is $\D_{X,p^n} = \DC_X[2]$ by definition (cf.\ \S\ref{sect0-2}), this equality implies the last claim in the part (iv) of \S\ref{sect0-2}. In fact, it is easy to see that the local version treated in this section can be extended to the more global situation described in the introduction.
\addtocounter{thm}{4}
\begin{lem}\label{lem.deltaloc} Consider cartesian squares of schemes
\[ \xymatrix{
X_Z \, \ar@{^{(}->}[r]^-{i} \ar[d] \ar@{}[rd]|{\square} & X \ar[d]_f \ar@{}[rd]|{\square} & \ar@{_{(}->}[l]_-{j} \, X_U \ar[d] \\
Z \, \ar@{^{(}->}[r]^-{i'} & Y & \ar@{_{(}->}[l]_-{j'} \, U\,,\hspace{-5pt}
}\]
where $i'$ is a closed immersion and $j'$ is the open immersion of the complement $U = Y \ssm Z$. Then, for any complex of torsion sheaves $\cK \in D^+(Y_\et)$ the base-change isomorphisms give an identification
\[ Rf^!(\delta^\loc_{U,Z}(\cK)) = \delta^\loc_{X_U,X_Z}(Rf^!\cK)\,. \]
\end{lem}
\begin{pf}
There is a commutative diagram with distinguished rows
\[ \xymatrix{
i_* Ri^!(Rf^!\cK) \ar[r]^-{i_*} \ar[d]_\beta^{\rwr} & Rf^!\cK \ar[r]^-{j^*} \ar@{=}[d] & Rj_* j^*(Rf^!\cK) \ar[rr]^-{-\delta^\loc_{X_U,X_Z}(Rf^!\cK)} \ar@{.>}[d]_\alpha && i_* Ri^!(Rf^!\cK)[1] \ar[d]_{\beta[1]}^{\rwr} \\
Rf^!(i'_* Ri'^!\cK) \ar[r]^-{Rf^!(i'_*)} & Rf^!\cK \ar[r]^-{Rf^!(j'^*)} & Rf^!(Rj'_* j'^*\cK) \ar[rr]^-{-Rf^!(\delta^\loc_{U,Z}(\cK))} && Rf^!(i'_* Ri'^!\cK)[1]\,,\hspace{-5pt}
}\]
where the top row is the localization exact triangle \eqref{eq-0.4.2} for $Rf^!\cK$, the bottom row is obtained by applying $Rf^!$ to a localization exact triangle for $\cK$ and the arrow $\beta$ is a base-change isomorphism. By adjunction, the left hand square commutes. Therefore there exists a morphism $\alpha$ which makes the other squares commute (see \S\ref{sect0-5-2}). By the commutativity of the middle square, $\alpha$ is mapped to the identity under the canonical isomorphisms
\begin{align*}
\Hom_{D(X)}(Rj_* j^*Rf^!\cK, Rf^!Rj'_*j'^*\cK) & \cong \Hom_{D(U)}(j^* Rf^!\cK, j^* Rf^!Rj'_* j'^*\cK) \\
& = \Hom_{D(U)}(j^* Rf^!\cK, j^* Rf^!\cK)\,.
\end{align*}
But this means that $\alpha$ is the base-change isomorphism, and the claim follows.
\end{pf}

\smallskip
\subsection{Bloch-Ogus complexes and Kato complexes}\label{sect3-10}
As an application, used in \cite{js} 2.20, 2.21, we deduce the following result on Kato complexes, analogous to \S\ref{sect1-5} and \S\ref{sect2-11}. As in \cite{js} p.\ 497, we define a homology theory on all separated $S$-schemes $f : X \to S$ of finite type by letting
\[ H_a(X/S,\zpn(-1)) := H^{-a}(X,Rf^!(\zpn(1)'_S))\,, \]
and, following the method of Bloch and Ogus, a niveau spectral sequence
\begin{equation}\label{arspectralsequence}
E^1_{q,t}(X/S,\zpn(-1)) = \bigoplus_{x\in (X/S)_q} H_{q+t}(x/S,\zpn(-1)) \Ra H_{q+t}(X/S,\zpn(-1)),\hspace{-10pt}
\end{equation}
where $H_a(x/S;\zpn(-1))$ is defined as the inductive limit over all $H_a(U/S,\zpn(-1))$, for all non-empty open subschemes $V\subset \overline{\{x\}}$. See \eqref{eq1-5-3'} for the definition of $(X/S)_q$. Then we have
\stepcounter{thm}
\begin{thm} \label{thm:katoc3}
\begin{enumerate}
\item[(1)] For $X=X_\eta$ the spectral sequence \eqref{arspectralsequence} is canonically isomorphic to the spectral sequence
\begin{align*} E^1_{q-1,t+1}(X_\eta/\eta,\zpn(-1)) = & \bigoplus_{x\in (X_\eta)_{q-1}} \ H_{q+t}(x/\eta,\zpn(-1)) \\
 & \qquad \qquad \Lra H_{q+t}(X_\eta/\eta;\zpn(-1))\,. \end{align*}
obtained  from \eqref{eq1-5-4} after shifts in the both degrees.
\item[(2)] For $X=X_s$ the spectral sequence \eqref{arspectralsequence} is canonically isomorphic to the
spectral sequence
\[ E^1_{q,t+2}(X_s/s,\zpn) = \bigoplus_{x\in (X_s)_q} \ H_{q+t+2}(x/s,\zpn) \Lra H_{q+t+2}(X_s/s;\zpn) \]
obtained from \eqref{eq2-11-2} after a shift in the second degree.
\item[(3)] Let $x \in (X/S)_q\cap X_s=(X_s)_q$ and $y\in (X/S)_{q+1} \cap X_\eta=(X_\eta)_q$ with $x\in\overline{\{y\}}$. Then there are canonical purity isomorphisms
\begin{align*}
H_{q+1+t}(y/S,\zpn(-1)) & \cong  H^{q-t-1}(y,\zpn(q+1))\,, \\
H_{q+t}(x/S,\zpn(-1)) & \cong  H^{q-t-2}(x,\zpn(q))\,.
\end{align*}
Via these isomorphisms, the $(y,x)$-component
\[ d^1_{q+1,t}(y,x): H^{q-t-1}(y,\zpn(q+1)) \to H^{q-t-2}(x,\zpn(q)) \]
of the differential $d^1_{q+1,t}$ in \eqref{arspectralsequence} coincides with $-\dKyx$.
\item[(4)]
The isomorphisms in {\rm(1)}, {\rm(2)} and {\rm(3)} induce isomorphisms
\[ E^1_{*,t}(X/S,\zpn(-1)) \cong C^{-t-2,0}_{p^n}(X)^{(-)} \]
between Bloch-Ogus complexes and sign-modified Kato complexes.
\end{enumerate}
\end{thm}
\begin{pf}
(1) and (2) are obvious from the isomorphisms \eqref{isos.DC-S}. The first claim in (3) is clear from the fact that $\overline{\{y\}}$ meets $X_s$, and the isomorphisms then follow from (1) and (2) and the purity isomorphisms \eqref{eq1-5-5} and \eqref{eq2-11-3}, respectively. For the third statement of (3) we recall that the upper exact triangle in \eqref{triangles.DC-X} induces isomorphisms \[ t_X: j_X^*\DC_X \cong \DC_{X_\eta} \quad \hbox{ and } \quad g_X : \DC_{X_s} \cong Ri_X^!\DC_X \] identifying $\delta_X^\Sval$ with the connecting morphism $-\delta^\loc_{X_\eta,X_s}(\DC_X)$, cf.\ \eqref{eq3-9-5}. Since $\delta_X^\Sval$ induces Kato's residue maps, we get the claim. As for (4), the compatibility $d^1(y,x) = -\dKyx$ between the differentials and Kato's residue maps follow from (1) and Theorem \ref{thm:Katocomplexprimetop} for $y,x \in X_\eta$, and from (2) and Theorem \ref{thm:Katocomplexp-case} for $y,x \in X_s$. The remaining case is covered by (3).
\end{pf}
\begin{rem}
It is easy to see that this theorem proves the claims in \cite{js} {\rm 2.20} and {\rm 2.21}, except that the signs needed to be corrected. The reason for this lies in the interpretation of the connecting morphism and the resulting minus sign in \eqref{eq-0.4.2}.
\end{rem}

\smallskip
\subsection{Unicity of the cone}\label{sect3-11}
As a complement we show the following unicity result for $\DC_X=Rf^!(\zpn(1)_S')$. Recall the situation at the beginning of this section
\[ \xymatrix{
X_\eta \, \ar@{^{(}->}[r]^-{j_X}\ar[d]_{f_\eta} \ar@{}[rd]|{\square}
& X \ar[d]_f \ar@{}[rd]|{\square}
& \,X_s \ar@{_{(}->}[l]_-{i_X}\ar[d]_{f_s} \\
\eta \, \ar@{^{(}->}[r]^-j
& S
& \ar@{_{(}->}[l]_-i \, s \,,\hspace{-5pt} }\]
and the associated exact triangle, cf.\ \eqref{triangles.DC-X}
\[\xymatrix{
i_{X_*}\DC_{X_s}[-2]  \ar[r] & \DC_X \ar[r] & Rj_{X*}\DC_{X_\eta} \ar[r]^-{\delta_X^\Sval} & i_{X*}\DC_X[-1]\,.
}\]
\begin{thm}\label{thm.cone}
The object $\mathcal E_X$ is uniquely determined, up to unique isomorphism, as the mapping fiber of $\delta_X^\Sval$.
\end{thm}
\begin{rem}\label{rem.cone}
There is a similar uniqueness claim in \S\ref{sect0-2}\,{\rm(iv)} in the situation where $S$ is the spectrum of the integer ring in a number field. This follows from the same arguments as the proof of Theorem \ref{thm.cone} below by replacing
\begin{align*}
 & j:\eta \hra S \qquad \hbox{ with } \quad S[p^{-1}] \hra S, \\
 & i:s \hra S \qquad \hbox{ with } \quad \varSigma \hra S, \\
 & j_X:X_\eta \hra X \qquad \hbox{ with } \quad X[p^{-1}] \hra X, \\
 & i_X:X_s \hra X \qquad \hbox{ with } \quad X_\varSigma \hra X,
\end{align*}
where $\varSigma$ denotes the set of all closed points on $S$ of characteristic $p$.
\end{rem}
\begin{pf*}{Proof of Theorem \ref{thm.cone}}
By Lemma \ref{lem.derived}\,(3) it suffices to show
\smallskip
\begin{enumerate}
\item[(i)] $\Hom^{-1}_{D(X,\zpn)}(i_{X*}\DC_{X_s}[-2],Rj_{X*}\DC_{X_\eta})=0$. \smallskip
\item[(ii)] $\Hom_{D(X,\zpn)}(Rj_{X*}\DC_{X_\eta},i_{X_*}\DC_{X_s}[-2])=0$.
\end{enumerate}
\smallskip \noindent
(i) follows by adjunction for $j_X$, because $j^*_Xi_{X*}=0$. As for (ii), since
\[ \Hom_{D(X,\zpn)}(Rj_{X*}\DC_{X_\eta},\,i_{X*}\DC_{X_s}[-2])= \Hom_{D(s,\zpn)}(Rf_{s!}i_X^{*}Rj_{X*}\DC_{X_\eta},\, \zpn[-2]) \]
by adjunction for $i_X$ and $f_s$, it suffices to show
\begin{lem}\label{lem.boundmixed}
$Rf_{s!}i^*_XRj_{X*}\DC_{X_\eta}$ is concentrated in $[-2d,1]$, where $d = \dim(X_\eta)$.
\end{lem}
\noindent We first show the following result, which may be of own interest.
\begin{lem}\label{lem.boundfield} Let $k$ be a field, and let $f : X \to \Spec(k)$ be separated of finite type, and let $n$ be a positive integer which is invertible in $k$. Then $Rf_!Rf^!\bZ/n(i)$ is concentrated in $[-2d,0]$, where $d:=\dim(X)$.
\end{lem}
\begin{pf} We proceed by induction on $d=\dim(X)$. We may assume that $k$ is separably closed, that $i=0$ and that $X$ is reduced, and then the case $d=0$ is clear. Choose an affine open subset $U \subset X$ which is smooth of pure $\dim d$ and whose complement $Z := X \ssm U$ has dimension $e \le d-1$. We get a commutative diagram
\[ \xymatrix{
U \,\ar@{^{(}->}[r]^\phi \ar[dr]_{f_U}
& X \ar[d]^f
& \, Z\ar@{_{(}->}[l]_\psi \ar[dl]^{f_Z} \\
& \Spec(k) \,, \hspace{-5pt}} \]
where $\phi$ (resp.\ $\psi$) denotes the natural open (resp.\ closed) immersion, and we defined $f_U:=f \circ \phi$ and $f_Z:=f \circ \psi$. We note that $\phi$ is affine, because $X$ is separated over $k$ (if $V\subset X$ is affine, then $\phi^{-1}(V)=U\cap V$ is affine). There is an exact triangle
\[ Rf_!\psi_* R\psi^!Rf^!\bZ/n \lra Rf_!Rf^!\mathbb Z/n \lra Rf_!R\phi_*\phi^* Rf^!\bZ/n \os{+1}\lra . \]
Since $U$ is smooth of pure dimension $d$, we have
\[ \phi^* Rf^!\bZ/n = Rf^!_U \bZ/n \cong \bZ/n(d)[2d] \,. \]
Moreover we have $R\psi^!Rf^!=Rf_Z^!$. Therefore we can identify the above triangle with
\addtocounter{equation}{3}
\begin{equation}\label{triangle.localization1}
Rf_{Z!} Rf^!_Z \bZ/n \lra Rf_!Rf^!\bZ/n \lra Rf_!R\phi_* \bZ/n(d)[2d] \os{+1}\lra .
\end{equation}
Since $Rf_{Z!}Rf^!_Z\bZ/n$ is concentrated in $[-2d+2,0]$ by induction, it is enough to show that $A:=Rf_! R\phi_* \bZ/n(d)[2d]$ is concentrated in $[-2d,0]$. Obviously $A$ is concentrated in degrees $\ge -2d$, because this holds for $\bZ/n(d)[2d]$. On the other hand, we note that $\bZ/n(d)[d]$ is a perverse sheaf on $U$ (\cite{pervers} p.\ 102), so that $R\phi_*\bZ/n(d)[d]$ is perverse, because $\phi$ is an affine open immersion (and hence $t$-exact for the perverse $t$-structure loc.\ cit.\ 4.1.10\,(i)), and that $A=Rf_!R\phi_*\,\bZ/n(d)[d]$ is of perversity $\le d$ (loc.\ cit.\ 4.2.4), i.e., lies in $D^{p \le d}_c(k,\bZ/n)=D_c^{p\le 0}(k,\bZ/n)[-d]$. This means that
\[ A \in D^{p\leq 0}_c(k,\mathbb Z/n). \]
Since the perverse $t$-structure is the classical $t$-structure on $\Spec(k)$, we get that $A$ is concentrated in degrees $\le 0$. Thus we obtain Lemma \ref{lem.boundfield}.
\end{pf}

\noindent{\bf Proof of Lemma \ref{lem.boundmixed}.}
We may assume that $X$ is reduced and the closure of $X_\eta$. Then we prove the lemma by induction on $d=\dim X_\eta$. The case $d=0$ is easy and left to the reader. Suppose $d>1$. Then there is a commutative diagram
\[ \xymatrix{U \, \ar@{^{(}->}[r]^-\phi \ar[dr] & X \ar[d]^f & \, Z\ar@{_{(}->}[l]_\psi \ar[dl]^g \\ & S\,,\hspace{-5pt}} \]
where $\phi$ is an open immersion, $U_\eta$ is affine, smooth over $\eta$ and has pure dimension $d$, $\psi$ is the closed immersion of the complement $Z = X \ssm U$ (with reduced subscheme structure), and $\dim Z_\eta\leq d-1$. We get an exact triangle
\begin{equation}\label{triangle.localization2}
i^*_X Rj_{X*}\psi_{\eta*}R\psi_\eta^!\DC_{X_\eta} \lra i^*_X Rj_{X*}\DC_{X_\eta}
\lra i^*_X Rj_{X*}R\phi_{\eta*}\phi_\eta^* \DC_{X_\eta} \os{+1}\lra\,,
\end{equation}
where we used morphisms in the following diagram:
\[ \xymatrix{
U_\eta \, \ar@{^{(}->}[r]^-{j_U} \ar@{_{(}->}[d]_{\phi_\eta} \ar@{}[rd]|{\square}
& \, U_{\phantom.} \ar@{_{(}->}[d]_\phi \ar@{}[rd]|{\square}
&\ar@{_{(}->}[l]_-{i_U} \ar@{_{(}->}[d]_{\phi_s} \, U_s \ar[rrd] \\
X_\eta \, \ar@{^{(}->}[r]^-{j_X}
& X
& \ar@{_{(}->}[l]_{i_X} \, X_s \ar[rr]^{f_s}
&& s \, . \\
Z_\eta \phantom{|}\hspace{-2pt} \ar@{^{(}->}[r]^-{j_Z} \ar@{^{(}->}[u]^{\psi_\eta} \ar@{}[ru]|{\square}
& Z \phantom{|}\hspace{-2pt} \ar@{^{(}->}[u]^\psi \ar@{}[ru]|{\square}
&\ar@{_{(}->}[l]_-{i_Z} \ar@{^{(}->}[u]^{\psi_s} \ar[rru]_{g_s} \, Z_s \phantom{|}\hspace{-2pt}
}\]
By the proper base-change theorem for $\psi$ we identify
\[ i^*_X Rj_{X*}\psi_{\eta*}R\psi_\eta^!\DC_{X_\eta} = i^*_X \psi_{*}Rj_{Z*}\DC_{Z_\eta} = \psi_{s*}i^*_ZRj_{Z*}\DC_{Z_\eta} \,. \]
Because $\phi$ is \'etale and $U_\eta$ is smooth of pure dimension $d$, we have
\[ R\phi_{\eta*}\phi^*_\eta \DC_{X_\eta}=R\phi_{\eta*}\DC_{U_\eta}=R\phi_{\eta*}\mupn^{\otimes d+1}[2d]\,.\]
Therefore triangle \eqref{triangle.localization2}, after application of $Rf_{s!}$, leads to an exact triangle
\[ Rg_{s!}i^*_ZRj_{Z*}\DC_{Z_\eta} \lra Rf_{s!}i^*_X Rj_{X*}\DC_{X_\eta}
\lra Rf_{s!}i^*_XRj_{X*} R\phi_{\eta*}\mupn^{\otimes d+1}[2d] \os{+1}\lra. \]
Since the first term is concentrated in $[-2d+2,1]$ by induction, it suffices to show that
\begin{equation}\label{eq3-11-6}
A:=Rf_{s!}i^*_XRj_{X*} R\phi_{\eta*}\mupn^{\otimes d+1}[2d]
\end{equation}
is concentrated in $[-2d,1]$. It is clearly concentrated in degrees $\ge -2d$, because this holds for $\mupn^{\otimes d+1}[2d]$. We prove that $A$ is concentrated in degrees $\le 1$ in what follows. By the proof of Lemma \ref{lem.boundfield}, $R\phi_{\eta*}\mupn^{\otimes d+1}[d]$ is a perverse sheaf, i.e.,
\[ \cP^q :=\cH^q(R\phi_{\eta*}\mupn^{\otimes d+1}[d])=R^{q+d}\phi_{\eta*}\mupn^{\otimes d+1} \]
has support in dimension $\le -q$. In particular, it is non-zero only for $-d \le q \le 0$. We will prove
\par\bigskip\noindent
{\it Claim. The sheaf $i^*_XR^mj_{X*}\cP^q$ is zero for $m+q>1$.}
\par\bigskip\noindent
We see that $i^*_XRj_{X*}R\phi_{\eta*}\mupn^{\otimes d+1}[d]$ is concentrated in degrees $\le 1$ by the claim and the Leray spectral sequence
\[ E_2^{a,b}=i^*_XR^aj_{X*}\cP^b \Lra \cH^{a+b}(i^*_X Rj_{X*}R\phi_{\eta*}\mupn^{\otimes d+1}[d])\,. \]
Moreover since $\dim X_s \le d$ (see the beginning of proof of Lemma \ref{lem.boundmixed}) and $\ch(s)=p$,
we see that \[ A[d]=Rf_{s!}i^*_XRj_{X*}R\phi_{\eta*}\mupn^{\otimes d+1}[d] \] is concentrated in degrees $\le d+1$, so that $A$ is concentrated in degrees $\le 1$. Thus it remains to show the above claim. By the remark before the claim, it suffices to prove
\addtocounter{thm}{3}
\begin{lem}\label{lem.specialization}
If $\cF$ is an \'etale sheaf on $X_\eta$ with $\dim(\Supp\,\cF)\le b$, then we have
\[ \dim(\Supp \; i^*_XR^mj_{X*}\cF) \le b \quad \hbox{ for \; $m \ge 0$,} \]
and $i^*_XR^mj_{X*}\cF=0$ for $m>b+1$.
\end{lem}
\begin{pf} By assumption, there is a closed subset $V \os{\iota}\hra X_\eta$ of dimension $\le b$ such that $\cF=\iota_* \cG$ with $\cG=\iota^*\cF$. Let $Y=\ol V$, the closure of $V$ in $X$ endowed with the reduced subscheme structure. Then $V=Y_\eta$, and $Y_s$ has dimension $\le b$. We get cartesian squares
\[\xymatrix{
Y_\eta \, \ar@{^{(}->}[r]^{j_Y} \ar@{_{(}->}[d]_{\iota=\kappa_\eta} \ar@{}[rd]|{\square}
& Y{\phantom{,}}\hspace{-2pt} \ar@{_{(}->}[d]_{\kappa} \ar@{}[rd]|{\square}
& \,Y_s\ar@{_{(}->}[l]_{i_Y} \ar@{_{(}->}[d]_{\kappa_s} \\
X_\eta \, \ar@{^{(}->}[r]^{j_X}
& X
& \,X_s \ar@{_{(}->}[l]_{i_X}
}\]
with $\kappa,\kappa_\eta$ and $\kappa_s$ being closed immersions. Since $\cF=\kappa_{\eta*}\cG$, we get
\[ i^*_XR^aj_{X*}\cF = i^*_XR^aj_{X*} \kappa_{\eta*}\cG = i^*_X\kappa_* R^qj_{Y*}\cG = \kappa_{s*}i^*_Y R^qj_{Y*}\cG\,,\]
where the last equality is a base-change isomorphism. This shows that $i^*_XR^mj_{X*}\cF$ has support in $Y_s$, i.e., in dimension $\le b$. Finallly, since $R^mj_{Y*}\cG=0$ for $m>b+1$ by Lemma \ref{lem1}, we have $i^*_XR^mj_{X*}\cF=\kappa_{s*}i_Y^*R^mj_{Y*}\cG=0$ for $m>b+1$.
\end{pf}
\noindent
This completes the proof of Lemma \ref{lem.boundmixed} and Theorem \ref{thm.cone}.
\end{pf*}
\par\medskip
By the above results, we obtain the following bounds for $\DC_X = Rf^!(\zpn(1)'_S)$.
\begin{cor}\label{cor.boundglobal}
Put $d:=\max(\dim X_\eta,\dim X_s)$. Then{\rm:}
\begin{enumerate}
\item[(1)] $i^* Rf_!\DC_X$ is concentrated in $[-2d,2]$.
\item[(2)] $j^* Rf_!\DC_X$ is concentrated in $[-2d,0]$.
\end{enumerate}
In particular, $Rf_!\DC_X$ is concentrated in $[-2d,2]$.
\end{cor}
\begin{pf}
(1) Consider the exact triangle
\[ \xymatrix{
i^{*}Rf_!i_{X*}\DC_{X_s}[-2] \ar[r] \ar@{=}[d] & i^* Rf_!\DC_X \ar[r] \ar@{=}[d] &
i^* Rf_!Rj_{X*}\DC_{X_\eta} \ar[r]^-{+1} \ar@{=}[d] & \,. \\
A & B & C }\]
Here
\[ A=Rf_{s!}\DC_{X_s}[-2] \quad \hbox{ and }\quad C=Rf_{s!}i_X^*Rj_{X*}\DC_{X_\eta} \]
by the proper base-change theorem. Since $C$ is concentrated in $[-2d,1]$ by Lemma \ref{lem.boundmixed}, it is enough to show that $A$ is concentrated in $[-d+2,2]$. Since $\DC_{X_s} \cong \cM_{n,X_s}$ by Theorem \ref{thm:jss}, the complex $\DC_{X_s}$ is concentrated in $[-d,0]$ and any non-zero section of $\cH^q(\DC_{X_s})$ has support of dimension $\le -q$.
This implies that
\[ R^mf_{s!}\cH^q(\DC_{X_s})=0 \quad\hbox{ for }\; m+q> 0. \]
Indeed, $Rf_{s!}$ commutes with inductive limits of sheaves, and for any separated of finite type morphism $g : Z \to s$ with $\dim(Z)=e$ and any $p$-primary torsion sheaf $\cF$ on $Z_\et$, the complex $Rg_!\cF$ is concentrated in $[0,e]$. Therefore $Rf_{s!}\DC_{X_s}$ is concentrated in $[-d,0]$, and $A$ is concentrated in $[-d+2,2]$.
\par
(2) Since
\[ j^* Rf_!\DC_X=Rf_{\eta!}\DC_{X_\eta} = Rf_{\eta!}Rf_\eta^!\mupn \,, \]
the assertion follows from Lemma \ref{lem.boundfield}.
\end{pf}

\newpage
\appendix
\section{Trace maps for logarithmic Hodge-Witt cohomology}\label{sectA}
\medskip
For the definition of the Kato complexes one needs corestriction maps
\setcounter{equation}{0}
\begin{equation}\label{eqA-0-1}
 \Cor_{L/K} : H^j(L,\logwitt L n r) \lra H^j(K,\logwitt K n r)
\end{equation}
for logarithmic Hodge-Witt cohomology and finite extensions $L/K$ of fields of characteristic $p>0$ (cf.\ \S\ref{sect0-6} or \cite{kk:hasse}). These are not defined explicitly in \cite{kk:hasse}, but Kato constructed such maps in earlier papers and referred to results in these papers. In this appendix we discuss Kato's construction and some alternative descriptions used in the main body of the paper (cf.\ Lemmas \ref{lem:norm0}, \ref{lem:norm1} and Corollary \ref{cor:norm1} below). Recall that the groups above are non-zero only for $j=0, 1$.
\subsection{The case $\bs{j=0}$}\label{sectA-1}
First we consider the case $j=0$. Here the definition \eqref{eq-0.6.6} works for arbitrary extensions $L/K$. But in the following situation this corestriction map coincides with Gros' Gysin maps.
\begin{lem}\label{lem:norm0}
Let $k$ be a perfect field of characteristic $p>0$, and let $r \ge 0$ be an integer. For a finite extension $L/K$ of finitely generated fields over $k$, the following diagram is commutative{\rm:}
\[\xymatrix{
K^M_r(L)/p^n \ar[r]^-{N_{L/K}} \ar[d]_{h^r=\dlog} & K^M_r(K)/p^n \ar[d]^{h^r=\dlog} \\
\H^0(z,\logwitt z n r) \ar[r]^{\gys_f} & \H^0(x,\logwitt x n r)\,,\hspace{-5pt}
}\]
where $z:=\Spec(L)$, $x:=\Spec(K)$ and $\gys_f$ denotes Gros' Gysin map for the finite morphism $f : z \to x$, cf.\ {\rm\S\ref{sect2-1}}. The vertical arrows are the differential symbol maps, and $N_{L/K}$ denotes the norm map of Milnor $K$-groups. In other words, the corestriction map $\Cor_{L/K}$ of \eqref{eq-0.6.6} coincides with $\gys_f$.
\end{lem}
This property was first shown by Shiho under the assumption $r=[K:K^p]$ (unpublished). Later he gave a proof for general $r$ but under the assumption $n=1$ (\cite{shm} p.\ 624 Claim 2). We include a simplified proof of Lemma \ref{lem:norm0} to extend his result to general $r$ and $n$, which will be useful in \S\ref{sectA-2} below.
\par
\begin{pf}
By the transitivity properties of Gros' Gysin maps (cf.\ (P2) in \S\ref{sect2-1}) and the norm maps (\cite{kk:cft} p.626 Proposition 5), we may suppose that $L/K$ is a simple extension, i.e., $L=K(\alpha)$ for some $\alpha \in L$. Fix an $K$-rational point $\infty$ on $\P^1:=\P^1_K$ and an affine coordinate $t$ on $\P^1 \ssm \{\infty \}$. We regard $z = \Spec(L)$ as the closed point on $\P^1 \ssm \{\infty \}$ corresponding to the minimal polynomial (in $t$) of $\alpha$ over $K$. By a result of Bass and Tate \cite{BT} p.\ 379 (7), there is an exact sequence
\addtocounter{equation}{1}
\begin{equation}\label{eqA-1-2}
K^M_{r+1}(K(t)) \os{\partial}\lra \bigoplus_{v \in (\P^1)_0} \ K^M_r(\kappa(v)) \os{N}\lra K^M_r(K) \lra 0\,,
\end{equation}
where $N$ denotes the sum of the norm maps $(N_{v/x})_{v \in (\P^1)_0}$ of Milnor $K$-groups. This sequence yields the upper exact row in the following commutative diagram:
\vspace{-15pt}
\begin{equation}\label{eqA-1-3}
\xymatrix{
K^M_{r+1}(K(t))/p^n \ar[r]^-{\partial} \ar[d]_-{\dlog}
& \displaystyle \bigoplus_{v \in (\P^1)_0}^{\phantom{v \in (\P^1)_0}} \ K^M_r(\kappa(v))/p^n \ar@{->>}[r]^-{N} \ar[d]_-{\dlog}
& K^M_r(K)/p^n \\
\us{\phantom{|^{|^|}}} \H^0(\eta,\logwitt \eta n {r+1}) \ar@<6pt>[r]^-{\partial'}
& \displaystyle \bigoplus_{v \in (\P^1)_0} \ \H^0(v,\logwitt v n r) \ar@<6pt>[r]^-{G}
& \us{\phantom{|^{|^|}}} \H^0(x,\logwitt x n r)\,.\hspace{-5pt}
}
\end{equation}
Here we put $\eta:=\Spec(K(t))$, and $\partial'$ is induced by $(\dval_{\eta,v})_{v \in (\P^1)_0}$. The square commutes by the definition of $\dval_{\eta,v}$'s. The arrow $G$ denotes the sum of the maps $(G_v)_{v \in (\P^1)_0}$, where $G_v$ is Gros' Gysin map for the morphism $v \to x$. We will show
\par\vspace{8pt} \noindent
{\it Claim. The lower row of \eqref{eqA-1-3} is a complex.}
\par\vspace{8pt} \noindent
This claim implies that the above commutative square induces a map
\[ K^M_m(K)/p^n \lra \H^0(x,\logwitt x n r)\,. \]
Because the $\infty$-components of $N$ and $G$ are identity maps, this induced map must be $\dlog$. In particular, $N_{z/x}$ commutes with $G_z=\gys_f$ via the $\dlog$ maps. Therefore it remains to show the claim.
\par\vspace{8pt} \noindent
{\it Proof of Claim.}
Let $g : \P^1 \to x$ be the structure map, and consider the following diagram:
\vspace{-20pt}
\[\xymatrix{
\H^0(\eta,\logwitt \eta n {r+1}) \ar[r]^-{\partial'} \ar@{=}[d]
& \displaystyle \bigoplus_{v \in (\P^1)_0}^{\phantom{(\P^1)_0}} \ \H^0(v,\logwitt v n r) \ar[r]^-{G} \ar[d]_-{G'}
& \H^0(x,\logwitt x n r) \\
\us{\phantom{|^{|^|}}}\H^0(\eta,\logwitt \eta n {r+1}) \ar@<6pt>[r]
&\displaystyle \bigoplus_{v \in (\P^1)_0} \ \H^1_v(\P^1,\logwitt {\P^1} n {r+1}) \ar@<6pt>[r]
& \us{\phantom{|^{|^|}}} \H^1(\P^1,\logwitt {\P^1} n {r+1}) \ar[u]^-{\gys_g} \,,\hspace{-5pt}}\]
where $G'$ is induced by the Gros' Gysin maps and the lower row is the localization exact sequence. By the results in \S\ref{sect2-2}, which does not use this lemma, the left square commutes up to a sign. The right square commutes by the transitivity of Gros' Gysin maps ((P2) in \S\ref{sect2-1}). Hence the upper row is a complex. \qed
\par \vspace{8pt} \noindent
This completes the proof of Lemma \ref{lem:norm0}.
\end{pf}
\subsection{The case $j=1$}\label{sectA-2}
Now we consider the case $j=1$ of \eqref{eqA-0-1}. Kato again used a symbol map to define a corestriction map in this case for an arbitrary finite field extension $L/K$ of fields of characteristic $p$. Recall that one has an exact sequence of \'etale sheaves on $x=\Spec(K)$
\[ 0 \lra \logwitt x n r \lra \witt x n r \os{1-\fF}{-\hspace{-7pt}\lra} \witt x n r/dV^{n-1}\Omega_x^{r-1} \lra 0\,,\]
where $\fF$ denotes the Frobenius operator and $V$ denotes the Verschiebung operator. We get an associated `long' exact cohomology sequence
\[ 0 \to \H^0(x,\logwitt x n r) \to \witt K n r \os{1-\fF}{-\hspace{-7pt}\lra} \witt K n r/dV^{n-1}\Omega_K^{r-1}
  \os{\delta}\to \H^1(x,\logwitt x n r) \to 0\,. \]
This induces an isomorphism
\begin{equation}\label{eqA-2-1}
\Coker\big(\witt K n r \os{1-\fF}\lra \witt K n r/dV^{n-1}\Omega_K^{r-1}\big) \cong H^1(x,\logwitt x n r)\,.
\end{equation}
We adapt the definitions in \cite{kk:cft} (where a discrete
valuation field with residue field $K$ was treated) to directly
define a symbol map for $H^1(x,\logwitt x n r)$.
\stepcounter{thm}
\begin{defn}[\cite{kk:cft}]\label{def:P} Define the group $P^r_n(K)$ as
\[ P^r_n(K) := W_n(K) \otimes (K^\times)^{\otimes r}/J\,, \]
where $J$ is the subgroup of $W_n(K) \otimes (K^\times)^{\otimes r}$ generated by elements of the following forms{\rm:}
\begin{enumerate}
\item[(i)] $(\overbrace{0,\dotsc,0}^{i\text{ times}},a,0,\dotsc,0) \otimes a \otimes b_1 \otimes \dotsb \otimes b_{r-1}$ {\rm(}$0 \le i \le r-1$, $a, b_1,\dotsc,b_{r-1} \in K^\times${\rm)}.
\item[(ii)] $(\fF(w)-w) \otimes b_1 \otimes \dotsb \otimes b_r$ {\rm(}$w \in W_n(K),\ b_1,\dotsc,b_r \in K^\times${\rm)}. Here $\fF$ denotes the Frobenius operator on $W_n(K)$.
\item[(iii)] $w \otimes b_1 \otimes \dotsb \otimes b_r$ {\rm(}$w \in W_n(K),\ b_1,\dotsc,b_r \in K^\times$ with $b_i=b_j$ for some $i \ne j${\rm)}.
\end{enumerate}
\end{defn}
\bigskip\noindent
We will construct a map $h^r :  P^r_n(K) \to \H^1(x,\logwitt x n r)$, and show that it is bijective. First of all, there is a natural map
\begin{eqnarray*}
g^r : W_n(K) \otimes (K^\times)^{\otimes r}
& \lra
& \witt K n r/dV^{n-1}\Omega_K^{r-1} \\
w \otimes b_1 \otimes \dotsb \otimes b_r
&\lmt
& w \,\dlog(\ul{b_1}) \cdot \, \dotsb \, \cdot \dlog(\ul{b_r}) \; \mod dV^{n-1}\Omega_K^{r-1}
\end{eqnarray*}
($w \in W_n(K),\ b_1,\dotsc,b_r \in K^\times$). For $a \in K$, we wrote $\ul a \in W_n(K)$ for its Teichm\"uller representative. This map $g^r$ annihilates the elements of $J$ of the form (iii).
\begin{lem}\label{lem:VF}
Let $\omega$ be an element of $J$ of the form {\rm (i)} or {\rm (ii)}. Then $g^r(\omega)$ is contained in the image of $1-\fF$.
\end{lem}
\begin{pf}
The assertion is obvious for $\omega$ of the form (ii). We show the case that $\omega$ is of the form (i). Let $a, b_1,\dotsc,b_{r-1}$ be elements of $K^\times$ and let $i$ be an integer with $0 \le i \le n-1$. Put
\begin{align*}
\omega_i & := (V^i\ul a)\,\dlog(\ul a)\, \dlog(\ul{b_1}) \cdot \, \dotsb \, \cdot \dlog(\ul{b_{r-1}}) \; \in \; \witt K n r\,, \\
\tau_i & := d\ul a \, \dlog(\ul{b_1}) \cdot \, \dotsb \, \cdot \dlog(\ul{b_{r-1}}) \; \in \; \witt K {n-i} r\,.
\end{align*}
We will prove
\par\vspace{8pt} \noindent
{\it Claim. \; We have $\omega_i = V^i\tau_i$ in $\witt K n r$.}
\par\vspace{8pt} \noindent
We first finish the proof of Lemma \ref{lem:VF}, admitting this claim. By the proof of \cite{il} I.3.26, $\tau_i$ is contained in the image of $1-\fF: \witt K {n-i+1} r \ra \witt K {n-i} r$. Hence $\omega_i$ is contained in the image of $1-\fF: \witt K {n+1} r \ra \witt K n r$ by the claim and the equality $V\fF=\fF V$. The lemma immediately follows from this fact. Thus it remains to show the claim.
\par\vspace{8pt} \noindent
{\it Proof of Claim.} Since $\fF\, \dlog(\ul b) = \dlog(\ul b)$ ($b \in K^\times$) and
\[ Vx \cdot y = V(x \cdot \fF\, y)  \quad  \hbox{($x \in \witt K {n-1} r, \; y \in \witt K n {r'}$)}\]
by \cite{il} I.2.18.4, we have only to show the equality
\[ (V^i \ul a)\,\dlog(\ul a) = V^i d \ul a \quad \hbox{ in } \; \witt K n 1\,. \]
The case $i=0$ is clear. The case $i = 1$ follows from the equalities
\[ (V \ul a) \, \dlog(\ul a) = (V {\ul a}^{1-p}) \, d \ul a \os{\text{(V3)}}{=}
V({\ul a}^{p-1} \cdot {\ul a}^{1-p})\, dV \ul a = (V \ul{1}) \, dV \ul a \os{\text{(V2)}}{=} V d \ul a \,, \]
where the first equality follows from loc.\ cit.\ 0.1.1.9, and the equalities (V2) and (V3) mean those in loc.\ cit.\ I.1.1. Finally for $i \ge 2$, we have
\[ (V^i \ul a) \, \dlog(\ul a) = (V^{i-1} (V \ul a)) \, \dlog(\ul a) =
V^{i-1} ((V \ul a) \, \dlog(\ul a)) = V^{i-1} (V d \ul a) = V^i d \ul a \,. \]
This completes the proof of the claim and Lemma \ref{lem:VF}.
\end{pf}
\medskip
By the above, we get an induced map
\[ g^r: P^r_n(K) \to \Coker\big(\witt K n r \os{1-\fF}\lra \witt K n r/dV^{n-1}\Omega_K^{r-1}\big)\,, \]
and, by composition with the isomorphism \eqref{eqA-2-1}, the wanted symbol map
\[ h^r :  P^r_n(K) \to \H^1(x,\logwitt x n r)\,. \]
\begin{prop}
$h^r$ is bijective.
\end{prop}
\begin{pf}
The case $n=1$ follows from \cite{kk:cft} p.\ 616 Corollary. For the case $n \ge 2$, consider the commutative diagram with exact rows
\[ \xymatrix{
& P^r_{n-1}(K) \ar[r] \ar[d]_{h^r}
& P^r_n(K) \ar[r] \ar[d]_{h^r}
& P^r_1(K) \ar[r] \ar[d]_{h^r}
& 0 \\
0 \ar[r]
& \H^1(x,\logwitt x {n-1} r) \ar[r]
& \H^1(x,\logwitt x n r) \ar[r]
& \H^1(x,\Omega_{x,\log}^r)\ar[r]
& 0\,,\hspace{-5pt}}\]
where we put $x:=\Spec(K)$. The exactness of the lower row follows from \cite{CTSS} p.\ 779 Lemma 3 and the Bloch-Gabber-Kato theorem \cite{Bloch-Kato} 2.1. The exactness of the upper row is obtained from the natural isomorphisms
\[ P^r_n(K) \otimes \Z/p^i \cong  P^r_i(K) \qquad \hbox{ for } \;\; 1 \le i \le n. \]
Therefore the map $h^r$ is bijective by induction on $n \ge 1$.
\end{pf}
Now we come to the corestriction map defined by Kato. In \cite{kk:cft} p.\ 637 Corollary 4, it is shown that there is an exact sequence
\[ C^r_n(K) \lra C^r_n(K)/\{C_n^{r-1}(K),T \} \lra P^r_n(K) \lra 0\,, \]
where $C^r_n(K)$ is a group defined in terms of the group $T\wh{C}K_{r+1}(K)$ considered by Bloch \cite{B} and $T$ is an indeterminate used in defining $T\wh{C}K_{r+1}(K)$.
See \cite{kk:cft} p.\ 636 for the precise definition of $C^r_n(K)$. By this exact sequence $P^r_n(K)$ is expressed by algebraic $K$-groups, and in loc.\ cit.\ p. 637 Proposition 3\,(1),\,(2) and p.\ 658, Kato defined a transfer map
\[ \Tr_{L/K}:P^r_n(L) \lra P^r_n(K) \]
using the transfer map in algebraic $K$-theory. The crucial claim \cite{kk:hasse} 1.9 then relies on a result in \cite{kk:res} and the corestriction map $\Cor_{L/K}$ defined as the composite
\addtocounter{equation}{3}
\begin{equation}\label{eq-def-cor}
\xymatrix{
\Cor_{L/K} : \H^1(z,\logwitt z n r) \ar[r]^-{(h^r)^{-1}}
& P^r_n(L) \ar[r]^-{\Tr_{L/K}}
& P^r_n(K) \ar[r]^-{h^r}
&\H^1(x,\logwitt x n r)\,,
}
\end{equation}
where $z=\Spec(L)$. We show that this definition agrees with the one given in \eqref{eq-0.6.5}:
\stepcounter{thm}
\begin{lem}\label{lem:norm1}
The following diagram is commutative{\rm:}
\[ \xymatrix{
P^r_n(L) \ar[r]^-{\Tr_{L/K}} \ar[d]_{h^r}^{\hspace{-1.5pt}\wr}
& P^r_n(K) \ar[d]^{h^r}_{\wr\hspace{-1.5pt}} \\
\H^1(z,\logwitt z n r) \ar[r]^-{\tr_{z/x}}
& \H^1(x,\logwitt x n r)\,,\hspace{-5pt}
}\]
where $\tr_{z/x}$ denotes the corestriction map in the sense of \eqref{eq-0.6.8}.
\end{lem}
\begin{pf}
We prove this lemma in a similar way as for Lemma \ref{lem:norm0}. By the transitivity properties of the two transfer maps, we may suppose that $L/K$ is a simple extension, i.e., $z=\Spec(L)$ is a closed point on $\P^1:=\P^1_x$. Let $\eta$ be the generic point of $\P^1$, and consider a commutative diagram
\vspace{-13pt}
\stepcounter{equation}
\begin{equation}\label{eqA-2-2}
\xymatrix{
P^{r+1}_n(\kappa(\eta)) \ar[r]^-{\partial}
& \displaystyle \bigoplus_{v \in (\P^1)_0}^{\phantom{(\P^1)_0}} \ P^r_n(\kappa(v)) \ar[r]^-{\Tr}
& P^r_n(\kappa(x)) \\
\us{\phantom{|^{|^|}}}\H^1(\eta,\logwitt \eta n {r+1}) \ar[u]^{(h^r)^{-1}}_{\hspace{-1.5pt}\wr} \ar@<6pt>[r]^-{\partial'}
& \displaystyle \bigoplus_{v \in (\P^1)_0} \ \H^1(v,\logwitt v n r) \ar[u]^-{(h^r)^{-1}}_-{\hspace{-1.5pt}\wr} \ar@<6pt>[r]^-{\tr}
& \us{\phantom{|^{|^|}}} \H^1(x,\logwitt x n r) \ar@<6pt>[r]
& \us{\phantom{|^{|^|}}}0\,.\hspace{-5pt}}
\end{equation}
Here $\partial'$ is defined as $(\dval_{\eta,v})_{v \in (\P^1)_0}$, $\tr$ is the sum of the maps $\tr_{v/x}$, and $\Tr$ is the sum of the maps $\Tr_{v/x} = \Tr_{\kappa(v)/\kappa(x)}$. The arrow $\partial$ is a residue map induced by the residue maps of algebraic $K$-groups (cf.\ \cite{kk:cft} \S2.1 and p.\ 637 Proposition 3) and the upper row is a complex obtained from the localization theory in algebraic $K$-theory. The square commutes up to a sign (loc.\ cit.\ p.\ 660 Proof of Lemma 3).
By a similar argument as for Lemma \ref{lem:norm0}, we have only to show that the lower row of \eqref{eqA-2-2} is exact.
Consider the following diagram:
\vspace{-12pt}
\[\xymatrix{
K^M_{r+1}(\kappa(\eta))/p^n \ar[r]^-{\partial} \ar[d]_-{\dlog}^-{\hspace{-1.5pt}\wr}
& \displaystyle
\bigoplus_{v \in (\P^1)_0}^{\phantom{(\P^1)_0}} \ K^M_r(\kappa(v))/p^n \ar[r]^-{N} \ar[d]_-{\dlog}^-{\hspace{-1.5pt}\wr}
& K^M_r(\kappa(x))/p^n \ar[d]_-{\dlog}^-{\hspace{-1.5pt}\wr} \ar[r]
& 0 \\
\us{\phantom{|^{|^|}}}\H^0(\eta,\logwitt \eta n {r+1}) \ar@<6pt>[r]^-{\partial'}
& \displaystyle \bigoplus_{v \in (\P^1)_0} \ \H^0(v,\logwitt v n r) \ar@<6pt>[r]^-{\Cor}
& \us{\phantom{|^{|^|}}} \H^0(x,\logwitt x n r) \ar@<6pt>[r]
& \us{\phantom{|^{|^|}}} 0\,.\hspace{-5pt}} \]
where the maps are defined as in \eqref{eqA-1-3}, except that at the bottom we now have the map $\Cor$, the sum of the corestriction maps $\Cor_{v/x}$.
Then the left square commutes by the definition of the residue maps, and the right square commutes by the definition \eqref{eq-0.6.6} of the corestriction maps. The upper row is exact as we have seen in \eqref{eqA-1-3} and the vertical maps are all isomorphisms. Therefore the lower sequence is exact as well. Sheafifying in the \'etale topology for $K$ we obtain an exact sequence of sheaves whose stalks at the separable closure $\ol K$ of $K$ are
\[ \logwitt {\eta'} n {r+1} \os{\partial'}\lra \bigoplus_{v \in (\P^1)_0} \ \logwitt {v'} n r \os{\Cor}\lra \logwitt {x'} n r \lra 0\,,\]
with $x'=\Spec(\ol K)$, $\eta' = \eta\times_x x'$ and $v'=v\times_x x'$. By taking cohomology $\H^1(x,-)$ (which is a right exact functor on $p$-primary torsion sheaves) and applying Shapiro's lemma, we obtain an exact sequence
\[ \H^1(\eta,\logwitt \eta n {r+1}) \os{\partial'}\lra \bigoplus_{v \in (\P^1)_0} \ \H^1(v,\logwitt v n r)
\os{\tr}\lra \H^1(x,\logwitt x n r) \lra 0\,, \]
where $\tr$ is the sum of the race maps $\tr_{v/x}$ and $\partial'$ coincides with the map $\partial'$ in \eqref{eqA-2-2}, by the definition of the maps $\dval_{\eta,v}$. Therefore this sequence coincides with the lower row of \eqref{eqA-2-2}, which shows the exactness of the latter.
\end{pf}
By Lemmas \ref{lem:norm0} and \ref{lem:norm1} we immediately obtain:
\stepcounter{thm}
\begin{cor}\label{cor:norm1}
Under the same setting as in Lemma {\rm\ref{lem:norm0}}, the following diagram commutes{\rm:}
\[ \xymatrix{
P^r_n(L) \ar[r]^-{\Tr_{L/K}} \ar[d]_{h^r}^{\hspace{-1.5pt}\wr}
& P^r_n(K) \ar[d]^{h^r}_{\wr\hspace{-1.5pt}} \\
\H^1(z,\logwitt z n r) \ar[r]^-{\gys_f}
& \H^1(x,\logwitt x n r)\,.\hspace{-5pt}
}\]
In other words, the corestriction map in the sense of \eqref{eq-def-cor} coincides with $\gys_f$.
\end{cor}
\noindent This property was first shown by Shiho in the case that $[K:K^p] = p^r$ and $n=1$ (\cite{shm} p.\ 630 Claim 3).

\bigskip\noindent
    Uwe Jannsen:\\
    Fakult\"at f\"ur Mathematik, Universit\"at Regensburg\\
    Universit\"atsstr.\ 31, 93040 Regensburg\\
    GERMANY\\
    e-mail: uwe.jannsen@mathematik.uni-regensburg.de\\

    \noindent
    Shuji Saito:\\
    Department of Mathematical Sciences, University of Tokyo\\
    Komaba, Meguro-ku, Tokyo 153-8914\\
    JAPAN\\
    e-mail: sshuji@msb.biglobe.ne.jp\\

    \noindent
    Kanetomo Sato:\\
    Graduate School of Mathematics, Nagoya University\\
    Furo-cho, Chikusa-ku, Nagoya 464-8602\\
    JAPAN\\
    email: kanetomo@math.nagoya-u.ac.jp
\end{document}